\documentclass[a4paper,10pt, preprint,twocolumn]{article}
\usepackage{graphicx} % Required for inserting images
\usepackage[left=2cm, right=2cm, top=2cm, bottom=3cm]{geometry}
\usepackage{algorithm}
\usepackage{algpseudocode}
\usepackage{amsmath,amssymb}
\usepackage{comment}
\usepackage{hyperref}
\usepackage{stmaryrd}
\usepackage{bm}
\usepackage{authblk}
\usepackage{multicol}
\usepackage{indentfirst}
\usepackage{comment}
\usepackage{abstract}
\providecommand{\keywords}[1]
{
  \small	
  \textbf{\textit{Keywords: \quad}} #1
}

\title{Application of Sensitivity Analysis Methods for Studying Neural Network Models}
\author[1]{Jiaxuan Miao\footnote{jiaxuanmiao176@gmail.com}}
\affil[1]{Lomonosov MSU, Faculty of Computational Mathematics and Cybernetics, Moscow, Russia}
\author[1,2]{Sergey Matveev}
\affil[2]{Marchuk Institute of Numerical Mathematics, RAS Moscow, Russia}
\date{}

\begin{document}

\twocolumn[
\maketitle
\begin{onecolabstract}
This study demonstrates the capabilities of several methods for analyzing the sensitivity of neural networks to perturbations of the input data and interpreting their underlying mechanisms. The investigated approaches include the Sobol global sensitivity analysis, the local sensitivity method for input pixel perturbations and the activation maximization technique. As examples, in this study we consider a small feedforward neural network for analyzing an open tabular dataset of clinical diabetes data, as well as two classical convolutional architectures, VGG-16 and ResNet-18, which are widely used in image processing and classification. Utilization of the global sensitivity analysis allows us to identify the leading input parameters of the chosen tiny neural network and reduce their number without significant loss of the accuracy. As far as global sensitivity analysis is not applicable to larger models we try the local sensitivity analysis and activation maximization method in application to the convolutional neural networks. These methods show interesting patterns for the convolutional models solving the image classification problem. All in all, we compare the results of the activation maximization method with popular Grad-CAM technique in the context of ultrasound data analysis.
\end{onecolabstract}

  \\
 \saythanks{$^{*}$ Corresponding author: jiaxuanmiao176@gmail.com}
 \\
]

\section{Introduction}

At present, it is impossible to deny the pervasive influence of advances in machine learning and artificial intelligence on our lives. Even in 2024, the achievements of Nobel Prize winners in physics \cite{hopfield1982neural} and chemistry \cite{jumper2021highly} have been linked to discoveries in the development and application of artificial neural networks (ANNs). Research on artificial neural networks (ANNs) primarily focuses on developing new training rules, exploring novel neural network architectures, and expanding their application domains \cite{antsiferova2022video, liu2019neural, vaswani2017attention, kharyuk2018employing}. At the same time, less attention is given to the development of methods that help understand the nature of the internal representations for the networks solving the specific tasks \cite{matveev2021overview}. ANNs are often perceived by users as a kind of ``black box'', whose highly complex operation transforms the input data into the predefined outputs. In other words, it is not immediately clear how the weights of the ANN or the activation values of the hidden neurons are related to the processed dataset. Thus, unlike the classical statistical models, there is no straightforward approach determining the influence of each explanatory variable on the dependent variable in a neural network.

Gaining a deep understanding of their functioning and interpretability is becoming increasingly important and significant  with the extensive use and reliance on the neural networks \cite{ramsauer2020hopfield}. Numerous examples have long been known where relatively small perturbations in input data, known as adversarial attacks, impact the performance of neural networks in classification tasks \cite{khrulkov2018art, mopuri2017fast}.

In this context, the sensitivity analysis of neural networks enhances their interpretability \cite{angelov2021explainable}, helping to understand the mechanisms how input features influence the output results. This may improves model transparency and trustworthiness. Sensitivity analysis facilitates model structure optimization, eliminates redundant features, and enhances training and inference efficiency. Additionally, it enables an assessment of model robustness to input perturbations and strengthens defenses against adversarial attacks.

Ultimately, the sensitivity analysis can be highly beneficial in practical applications related to data collection and feature selection, which is particularly crucial in fields such as medicine. It also allows evaluating whether the model adequately captures causal relationships between input and output data and helps uncover hidden patterns, leading to valuable insights \cite{kharyuk2025exploring}. Therefore, sensitivity analysis serves as a key tool for understanding, optimizing, and validating neural networks.

\section{Global Sensitivity Analysis of a Feedforward Neural Network}
\subsection{Introduction to Data and Model}

At first, we make a global sensitivity analysis of a small and simple neural network. Sensitivity analysis has a great importance in the medical field \cite{smith1988using}. It also helps to verify whether the model adequately reflects the causal relationships between input and output data. Therefore, we have selected an open dataset \cite{akturk2020diabetes}, which includes factors such as the number of Pregnancies, Glucose, BloodPressure, SkinThickness, Insulin, BMI, DiabetesPedigreeFunction (an estimate of hereditary predisposition), and Age.

The patients in the dataset are divided into two classes based on the final diagnosis of whether they have diabetes or not. The task of the neural network is to predict the final diagnosis based on the available patient data.

\begin{table*}[ht]
 \centering
 \resizebox{\textwidth}{!}{%
 \begin{tabular}{@{}ccccccccc@{}}
 \textbf{Pregnancies} & \textbf{Glucose} & \textbf{BloodPressure} & \textbf{SkinThickness} & \textbf{Insulin} & \textbf{BMI} & \textbf{DiabetesPedigreeFunction} & \textbf{Age} & \textbf{Outcome} \\
 6                   & 148              & 72                    & 35                     & 0               & 33.6         & 0.627                          & 50           & 1               \\
 1                   & 85               & 66                    & 29                     & 0               & 26.6         & 0.351                          & 31           & 0               \\
 8                   & 183              & 64                    & 0                      & 0               & 23.3         & 0.672                          & 32           & 1               \\
 1                   & 89               & 66                    & 23                     & 94              & 28.1         & 0.167                          & 21           & 0               \\
 0                   & 137              & 40                    & 35                     & 168             & 43.1         & 2.288                          & 33           & 1               \\
 \vdots               & \vdots           & \vdots                & \vdots                  & \vdots           & \vdots       & \vdots                          & \vdots       & \vdots           \\ 
 \end{tabular}%
 }
 \caption{Diabetes data}
 \end{table*}

Due to the limitations of the dataset size and model scale, after parameter tuning and optimization, we finalize the model architecture with sufficiently good performance. The training dataset consists of 500 samples, while the test dataset contains 200 samples. The model accuracy on the training dataset is 77.9\%, and on the test dataset, it reaches 82.5\%, indicating that the model does not exhibit clear signs of overfitting. Typically, overfitting occurs when a model performs exceptionally well on the training data but poorly on the test data. However, in this case, the model achieves similar results on both training and test datasets, with slightly higher accuracy on the test set. This suggests that the model generalizes well and adapts effectively to new, unseen data. Therefore, it can be assumed that the model does not overemphasize the noise or specific patterns in the training data but instead captures general trends efficiently, avoiding overfitting.

The selected neural network consists of the two fully connected layers. The first layer is a fully connected layer with an input dimension of 8 and an output dimension of 10, followed by batch normalization and the \texttt{ReLU} activation function for nonlinear transformation. The second layer transforms the 10 features of the hidden layer into a single output, also using batch normalization. In the final stage, the model applies the nonlinear sigmoid function:
 $$
 \sigma(x) = \frac{1}{1 + e^{-x}},
 $$
 which compresses the output value into the range $[0, 1]$, which is suitable for a binary classification task. During training, binary cross-entropy is used as the loss function (\texttt{BCELoss}), where
$$\mathcal{L}(y, \hat{y}) = - \frac{1}{N} \sum_{i=1}^{N} \left( y_i \log(\hat{y}_i) + (1 - y_i) \log(1 - \hat{y}_i) \right),$$
 The batch normalization function in PyTorch for one-dimensional data is computed as follows:
\begin{enumerate}
     \item For the input data \( x = \{x_1, x_2, \dots, x_N\} \) the mean and variance are calculated for the input data:
$$
 \mu = \frac{1}{N} \sum_{i=1}^{N} x_i, \quad \quad \sigma^2 = \frac{1}{N} \sum_{i=1}^{N} (x_i - \mu)^2.
$$
 \item  We normalize the data ($\varepsilon = 10^{-5}$):
$$
 \hat{x}_i = \frac{x_i - \mu}{\sqrt{\sigma^2 + \varepsilon}}
$$

 \item We apply scaling parameters $\gamma$ and shifting $\beta$:
$$
 y_i = \gamma \hat{x}_i + \beta.
$$
 \end{enumerate}

 The stochastic gradient descent (SGD) method with $L_2$-regularization of the error function is used to optimize the model's parameters in order to prevent overfitting \cite{spall2005introduction}. The formula for updating the parameters with $L_2$-regularization is as follows:
 \begin{align*} 
 w_{n+1}=
 w_{n}-\eta \cdot \left. \frac{\partial \left(L(w) + \frac{\lambda}{2} ||w||_2^2\right)}{\partial w} \right|_{w_n}\\
= w_{n}-\eta \cdot \left(\frac{\partial L\left(w_{n}\right)}{\partial w}+\lambda w_{n}\right)
 \end{align*}
 $w_{n}$ is the current weight parameter, $\eta$ is the learning rate, $\dfrac{\partial L\left(w_{n}\right)}{\partial w}$ is the gradient of the loss function with respect to the weights, and $\lambda$ is the regularization coefficient. In this neural network, the regularization coefficient for $L_2$ regularization is set to $\lambda = 0.001$, and the learning rate is $\eta = 0.01$.
%\newpage

\subsection{Sobol Method for Global Sensitivity Analysis}

Sobol sensitivity analysis \cite{sobol2001global, rosolem2012fully} evaluates, using Monte Carlo numerical integration, the proportion of the variability of a scalar function 
$f(x_1, \ldots, x_k)$ that can be attributed to changes in the values of its input $x_i$. Let us consider the problem formulation:
 \begin{equation*}
 y=f(x)=f(x_{1},x_{2},\cdot \cdot \cdot,x_{k}),
 \end{equation*}
where $y$ is the value used to assess the performance of the model $f(x_1, \ldots, x_k)$ (for example, the mean squared error or simply a numerical output), and
$X=\left \{ x_{1},x_{2},\cdot \cdot \cdot ,x_{k} \right \}$ is the vector of $k$ model factors (parameters) that are assumed to influence $y$.

Our interest lies in understanding which portion of the total variance $V(y)$ in $y$ can be explained by the variability of the factors $x$. The Sobol method computes this by decomposing the function $y=f(x)=f(x_{1},x_{2},\cdot \cdot \cdot,x_{k})$ into components of increasing dimensionality, where each subsequent dimension represents higher-order interactions between the parameters
\begin{equation*}
 \begin{aligned}
 & f\left(x_{1}, x_{2}, \ldots, x_{k}\right)= f_{0}+\sum_{i=1}^{k} f_{i}\left(x_{i}\right) \\
 & + \sum_{1 \leq i \leq j \leq k} f_{i j}\left(x_{i}, x_{j}\right)+\ldots  +f_{1,2, \cdots, k}\left(x_{1}, x_{2}, \ldots, x_{k}\right).
 \end{aligned}
 \end{equation*}
Based on this decomposition, it can be shown that the total variance $V(y)$ in $y$ consists of these components
\begin{equation*}
 V(y)=\sum_{i} V_{i}+\sum_{i<j} V_{i j}+\ldots+V_{1,2, \cdots, k},
 \end{equation*}
where $V_{i}$ is the portion of the influence of the parameter $x_{i}$ on $V(y)$, while $V_{ij}$ is the portion contributed by the interaction between the parameters $x_{i}$ and $x_{j}$. The contribution of the $i$-th factor to the total variance of the output $V(y)$ is called the first-order sensitivity index $S_{i}$.
\begin{equation*}
 S_{i}=\frac{V_{i}}{V(y)} 
\end{equation*}
 The second-order indices $S_{ij}$ and the total Sobol indices $S_{Ti}$ are defined as:
\begin{equation*}
 S_{ij}=\frac{V_{ij}}{V(y)}, \quad S_{Ti}=1-\frac{V_{\sim i}}{V(y)},
\end{equation*}
where $V_{\sim i}$ is the variance obtained when all factors except $x_{i}$ are allowed to vary (with $x_{i}$ held at a fixed value). Thus, the total-order index $S_{Ti}$ represents the overall contribution of the factor $x_{i}$ to the total variance $V(y)$ through both direct and indirect effects.
 
\subsection{Result for Tiny Neural Network}
\begin{figure}[h]
    \centering
    \includegraphics[width=0.4\textwidth]{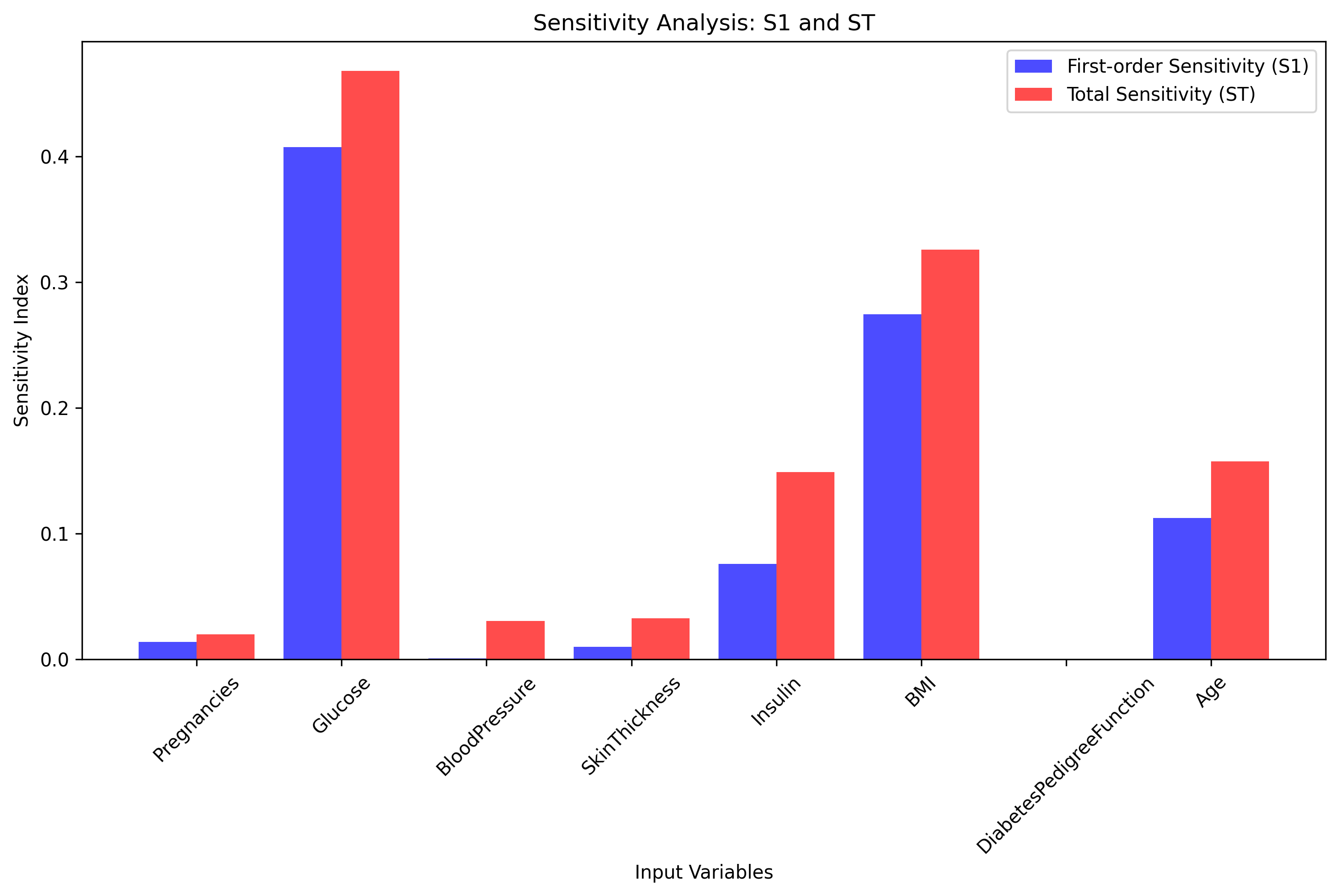}
    \caption{The Sobol indices $S_1$ and  $S_T$ for the output.}
    \label{fig:S1_ST}
\end{figure}

\begin{figure}[h]
    \centering
    \includegraphics[width=0.4\textwidth]{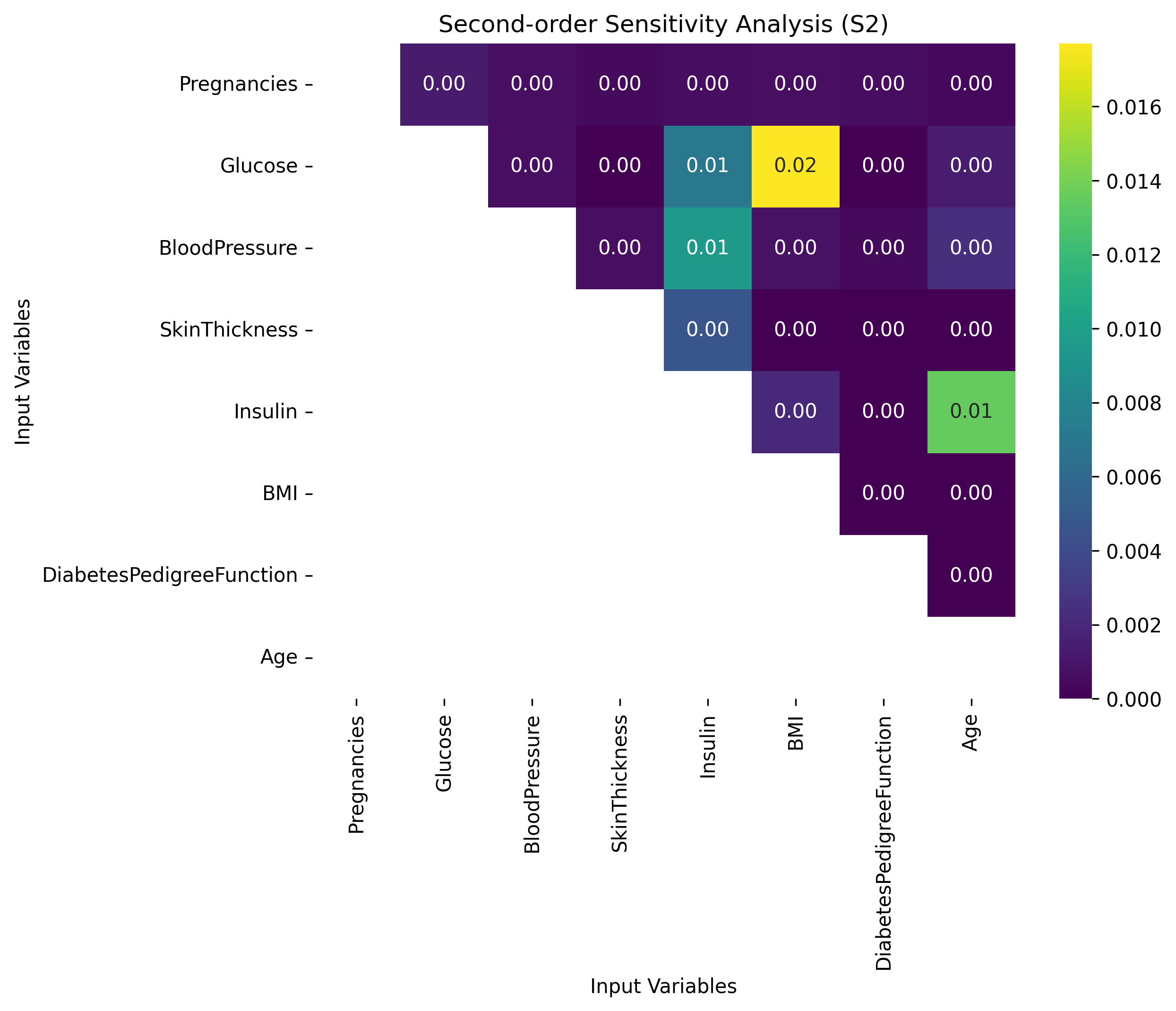}  
    \caption{The second-order Sobol indices $S_2$ for the output.}
    \label{fig:S2}
\end{figure}

\begin{figure}[h]
    \centering
    \includegraphics[width=0.4\textwidth]{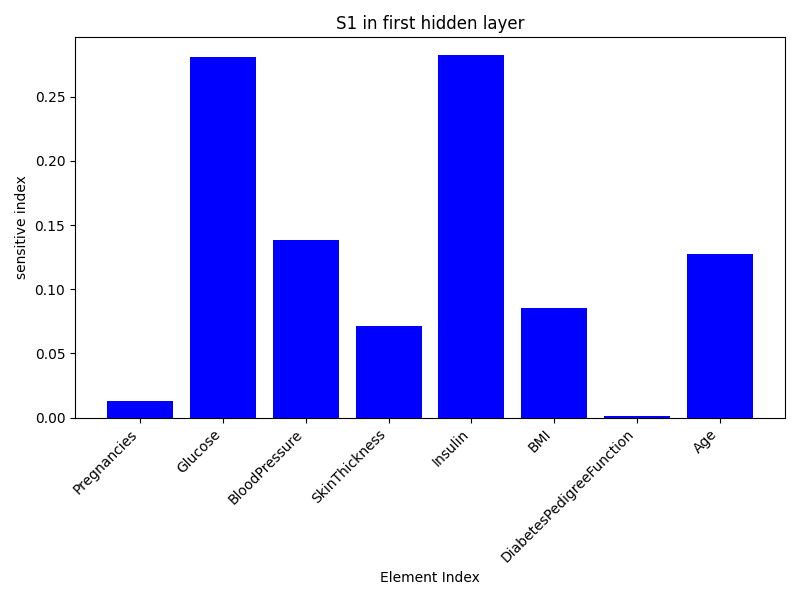}  
    \caption{$S_1$ for the first hidden layer.}
    \label{fig:S1_hidden}
\end{figure}

To better understand the behavior of the neural network, we analyzed the sensitivity of the input data to the first hidden layer and the output. The results are presented in Figures \ref{fig:S1_ST}, \ref{fig:S2}, \ref{fig:S1_hidden}. We used the SALib library in Python \cite{herman2017salib} to compute first-order, second-order, and total sensitivity for the output. The sample size in the presented calculations is 65536 quasi-random points \cite{sobol1998quasi}. The software implementation of the experiments is available \footnote{\url{https://github.com/LUCI1a/diabetes_prediction.git}}.

We can observe that for first-order sensitivity, the most significant variables are Glucose, BMI, and Age. After the first hidden layer, the influence of BloodPressure and Insulin becomes very large but then weakens, while the influence of Glucose, BMI, and Age remains significant. Thus, it can be concluded that for such a simple fully connected feedforward neural network, the sensitivity of input data gradually decreases as they propagate through the network and becomes more precise.

\begin{figure}[h]
    \centering
    \includegraphics[width=0.5\textwidth]{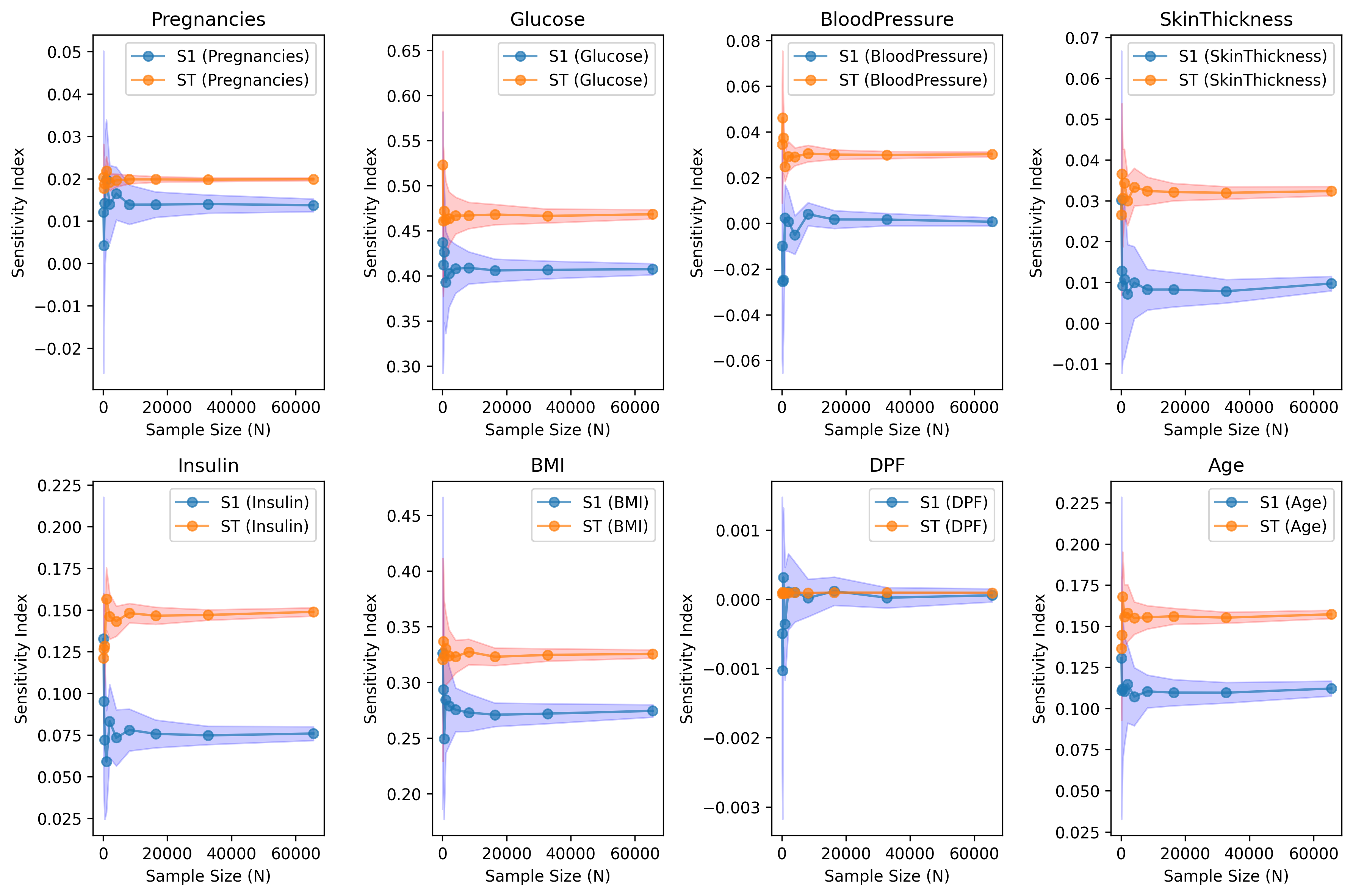}  
    \caption{$S_1$ and $S_T$ for various physiological indicators and sample sizes $N$. The shaded areas represent the corresponding confidence intervals.}
    \label{fig:convergence_Sobol}
\end{figure}

It is evident that increasing the sample size significantly improves the research results, making the sensitivity indices clearer and more distinct \cite{sobol1998quasi, gasanov2020sensitivity}. Our sample data ranged from 128 to 65536, and we observed a gradual convergence of the first and total sensitivity indices (Fig. \ref{fig:convergence_Sobol}). Ultimately, three key factors were identified as influencing diabetes: glucose level, body mass index, and age.

Based on the above, it can be concluded that the Sobol method is suitable for working with neural networks that have a small number of features \cite{sysoev2019sensitivity, zhang2020sensitivity}. These neural networks are suitable for tabular data and relatively small training samples. As the volume of the data for sampling increases, the width of the confidence intervals gradually converges to zero in accordance with theory. Therefore, even for small neural networks, a large number of quasi-random points is required to obtain high-quality approximations of sensitivity indices and to separate the confidence intervals.

One may wonder if nice results could be achieved by using only the most important features to build the neural network. Therefore, while keeping the original parameters and network architecture, we conduct training using only the four most influential features: Glucose, BMI, Age, and Insulin. As a result, the average accuracy on the test set was about 80\%, which closely matched the results on the original data. Additionally, when training the model on the remaining four less significant features — Pregnancies, BloodPressure, SkinThickness, and DiabetesPedigreeFunction — the accuracy on the test set was about 69\%. This suggests that the features Glucose, BMI, Age, and Insulin are sufficiently important and allow for a general description of the probability of diabetes development, while the information contained in Pregnancies, BloodPressure, SkinThickness, and DiabetesPedigreeFunction is less significant. If we had more data, we believe that the model's generalization ability would be higher, and the results would be more apparent.

Overall, this study confirmed the correctness and significance of using the Sobol method for sensitivity analysis, as well as its practical utility for independent testing and early diabetes screening.

\subsection{Comparison with PCA}
We analyzed the clinical data of all patients with diabetes and selected four principal components, whose weights account for more than 70\%. 

After applying principal component analysis (PCA), the most important indicators are Glucose, BMI, Insulin, and BloodPressure. In the global sensitivity analysis, the Glucose, BMI, and Insulin also turn out to be leading factors. However, the differences between the influences of various indicators after using PCA are not as large, and their impact on diabetes is not linear. Therefore, the visual effect of the principal component analysis (PCA) is much smaller than that of neural networks. Hence, conducting the sensitivity analysis for the neural networks is a very appropriate approach.

Although the global sensitivity analysis using the Sobol method is quite appealing, it is practically impossible to apply to larger convolutional models that have become mainstream in computer vision. However, for smaller and simpler neural networks, using the Sobol method for sensitivity analysis proves to be quite significant.

%\newpage

\section{Local sensitivity analysis for convolutional neural networks}

We have already conducted a general analysis of the behavior of a neural network using a fully connected feedforward neural network, and in this section we move to the analysis of the classical, yet more complex, VGG-16 convolutional neural network model \cite{simonyan2014very}. VGG-16 consists of 5 convolutional blocks, each containing several convolutional layers and one subsampling layer. We select the standard CIFAR-10 dataset consisting of color RGB images of size $32 \times 32 \times 3$ \cite{alex2009learning}. During the training, we used a following equation for normalization:
$$
x_{channel}^{\prime}=\frac{x_{channel}-\mu_{channel}}{\sigma_{channel}},
$$
Here, $\mu_{channel}$ and $\sigma_{channel}$ take the value 0.5, and the normalization parameter $\epsilon$ remains set to $10^{-5}$. After 15 epochs of training with a learning rate of 0.001, our model achieves an accuracy of 92\% on the training dataset and 84.9\% on the test dataset. Our goal is to investigate the behavior and role of different convolutional blocks. The software implementation of the experiments related to CIFAR-10 is available \footnote{Use the following link: \url{https://github.com/LUCI1a/CIFAR-10_projects.git}}

\begin{comment}
\begin{figure}[H]
    \centering
    \includegraphics[width=0.32 \textwidth]{images/vgg16.png} 
    \caption{The diagram of the VGG-16 model.}
    \label{fig:example_VGG} 
\end{figure}
\end{comment}

The global sensitivity analysis in the case of this architecture proves to be difficult due to the significantly larger dimensionality of the pixel space. Therefore, we focus on activation visualization techniques for the hidden layers and perform the local sensitivity analysis. We ust the following method use for sensitivity calculation:
\begin{eqnarray}
\notag
   & s_{{ij}_{C}} =  \left \| B_{q}\begin{pmatrix}
  x_{11} &  x_{12}&    \dots& x_{1n}\\
  x_{21} &  x_{22}&    \dots& x_{2n}\\
  \dots&  \dots&  x_{ij}+\varepsilon &  \dots& \\
  x_{m1} &  x_{m2}&    \dots& x_{mn}\\
\end{pmatrix}_{C} \right \| _{2 }  -
\\
\notag
& \left \| B_{q}\begin{pmatrix}
  x_{11} &  x_{12}&    \dots& x_{1n}\\
  x_{21} &  x_{22}&    \dots& x_{2n}\\

  \dots&  \dots&  \dots&  \dots& \\
  x_{m1} &  x_{m2}&    \dots& x_{mn}\\
\end{pmatrix}_{C}  \right \|_{2} 
\end{eqnarray}
$s_{{ij}_{C}}$ is the sensitivity of the pixel in a specific color channel,
$B{q}$ is the operator that describes the behavior of the $q$-th block, and $x_{ij}$ is the pixel value. Using this method, we obtain a sensitivity matrix for each block across the different color channels. As a generalization of this approach, partial derivatives, averages, and variances can be used when perturbing individual variables. Implementations of such generalizations are available in the \texttt{NeuralSens} package \cite{JSSv102i07}.

%\newpage

\subsection{Numerical Results}

\begin{figure}[htbp]
    \centering
%    \begin{subfigure*}{\textwidth}
        \includegraphics[width=0.092\linewidth]{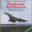}
        \includegraphics[width=0.092\linewidth]{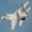}
        \includegraphics[width=0.092\linewidth]{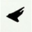}
        \includegraphics[width=0.092\linewidth]{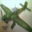}
        \includegraphics[width=0.092\linewidth]{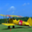}
        \includegraphics[width=0.092\linewidth]{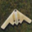}
        \includegraphics[width=0.092\linewidth]{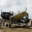}
        \includegraphics[width=0.092\linewidth]{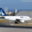}
        \includegraphics[width=0.092\linewidth]{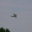}
        \\
%    \end{subfigure*}
    \vspace{0.2cm} % 图片间距
    % 第二行
%    \begin{subfigure*}\textwidth}
        \includegraphics[width=0.092\linewidth]{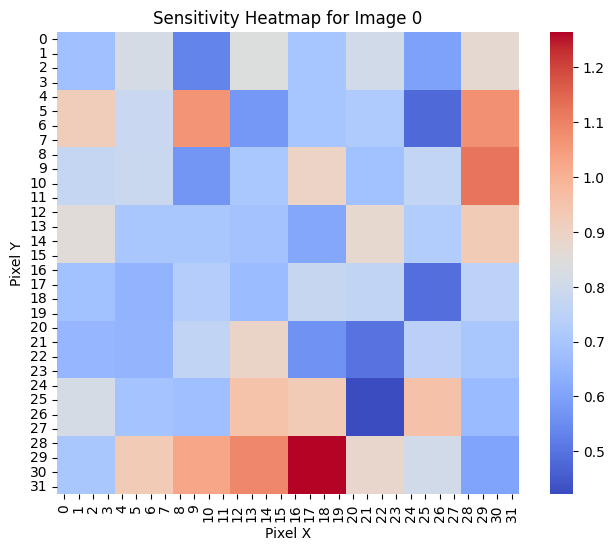}
        \includegraphics[width=0.092\linewidth]{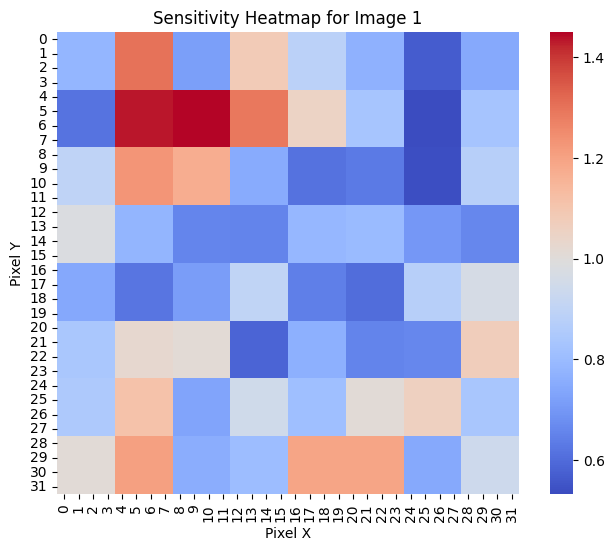}
        \includegraphics[width=0.092\linewidth]{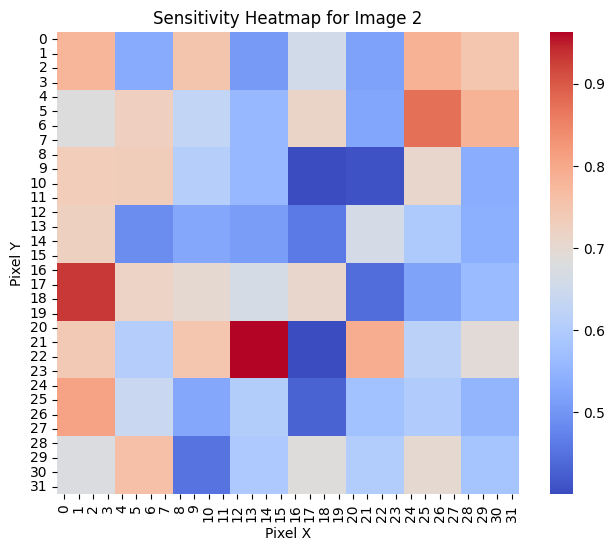}
        \includegraphics[width=0.092\linewidth]{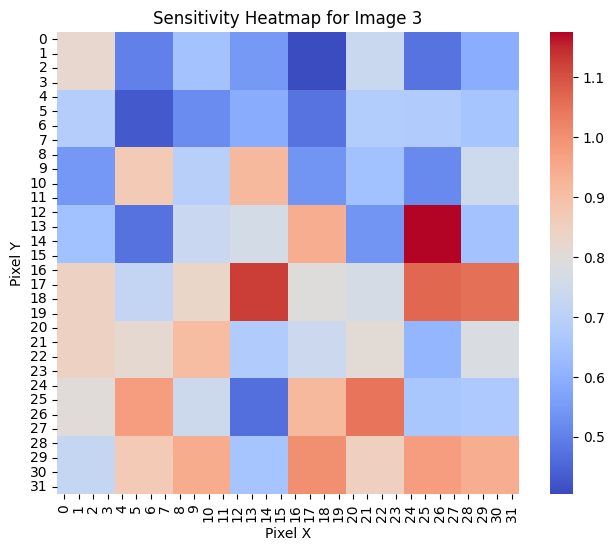}
        \includegraphics[width=0.092\linewidth]{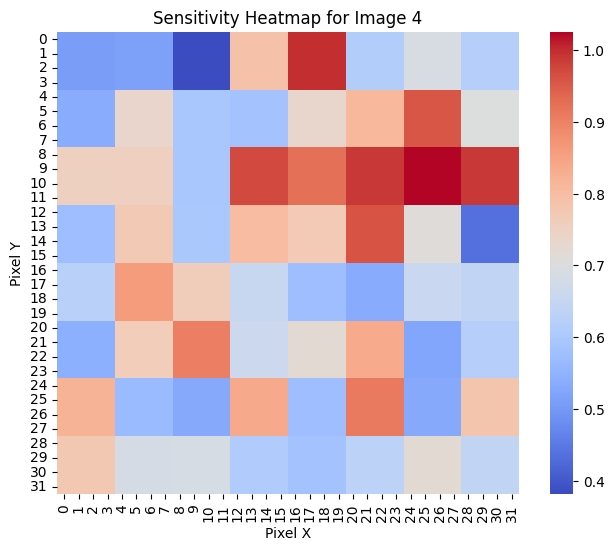}
        \includegraphics[width=0.092\linewidth]{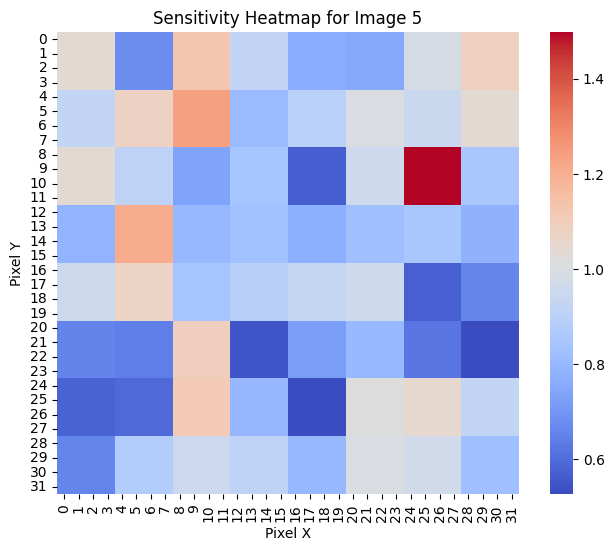}
        \includegraphics[width=0.092\linewidth]{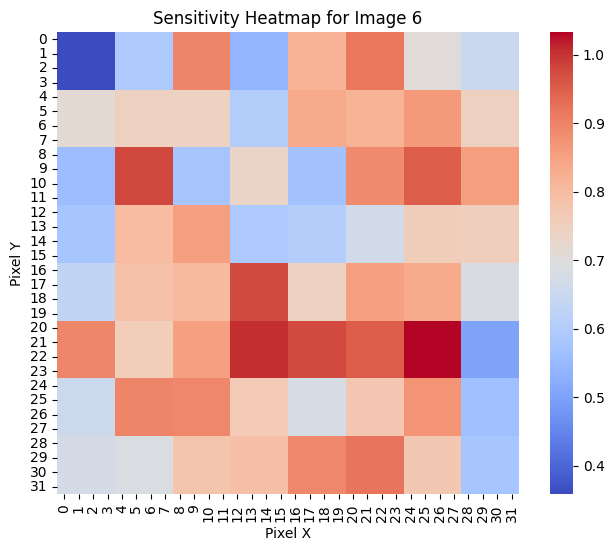}
        \includegraphics[width=0.092\linewidth]{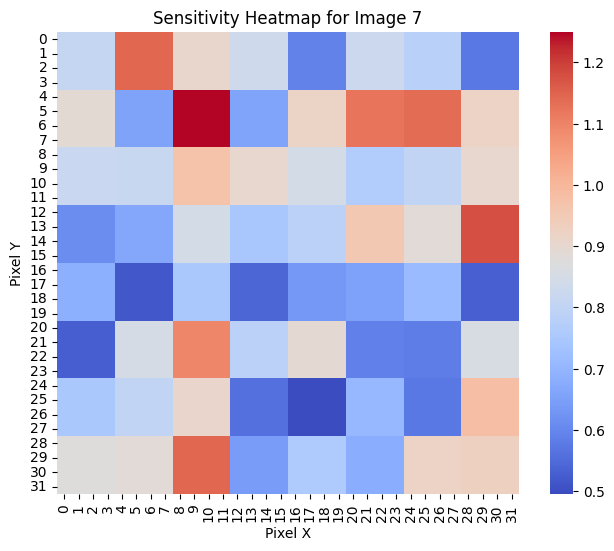}
        \includegraphics[width=0.092\linewidth]{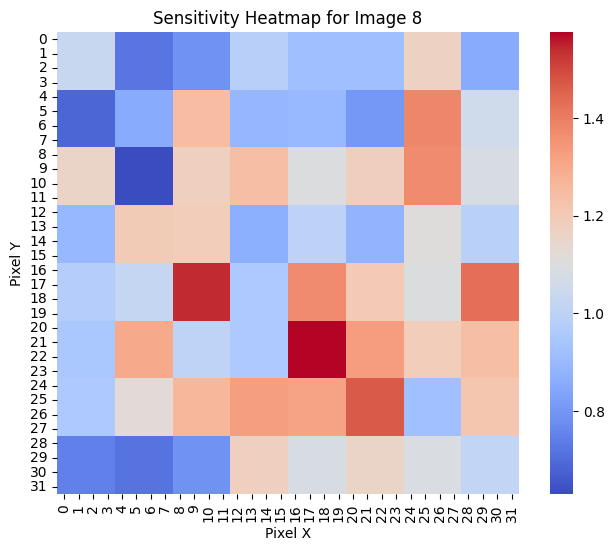}
        \\
%    \end{subfigure*}
%        \begin{subfigure*}{\textwidth}
        \includegraphics[width=0.092\linewidth]{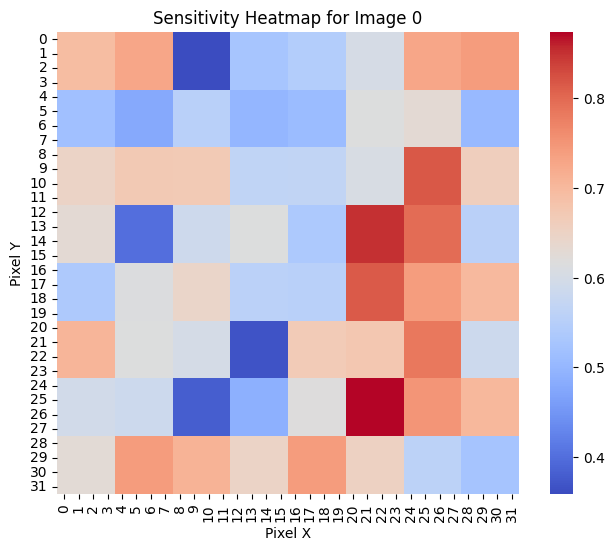}
        \includegraphics[width=0.092\linewidth]{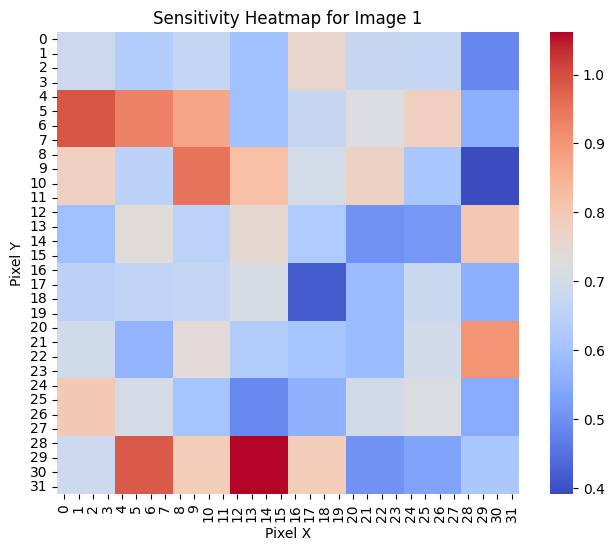}
        \includegraphics[width=0.092\linewidth]{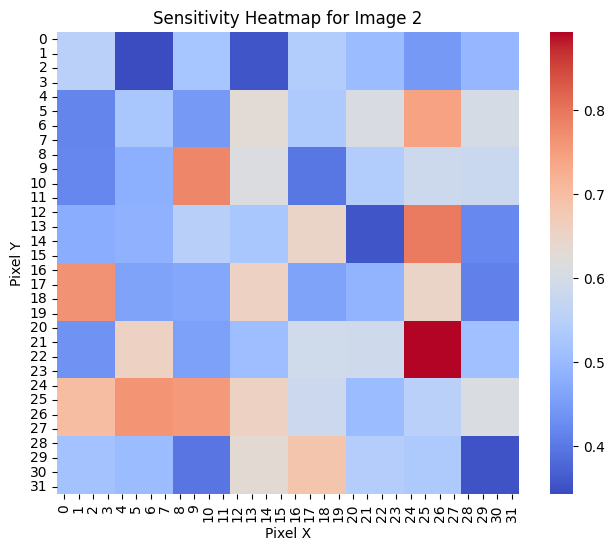}
        \includegraphics[width=0.092\linewidth]{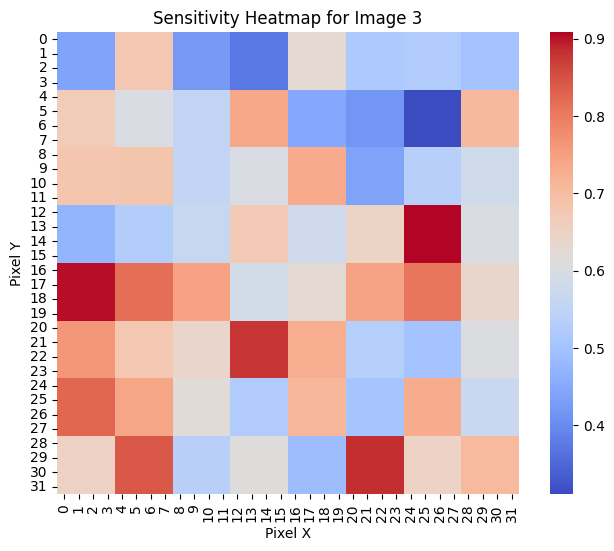}
        \includegraphics[width=0.092\linewidth]{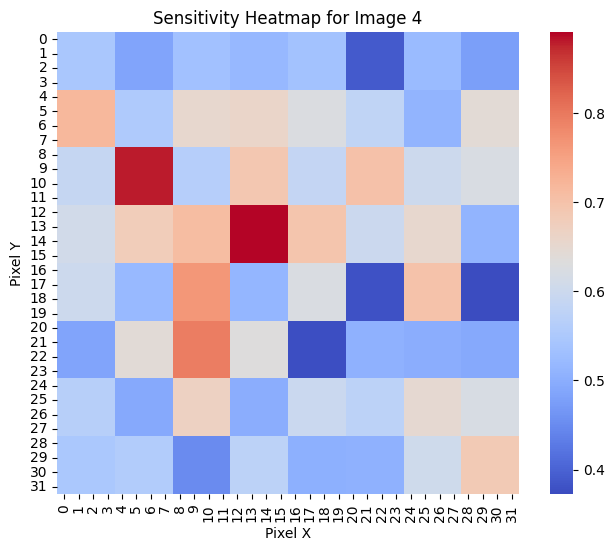}
        \includegraphics[width=0.092\linewidth]{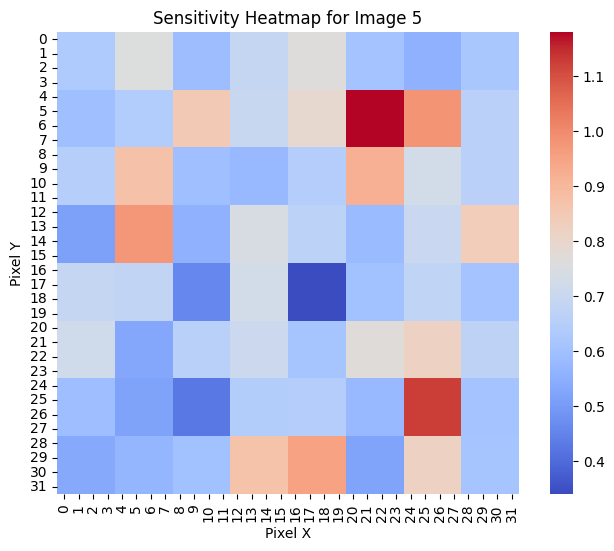}
        \includegraphics[width=0.092\linewidth]{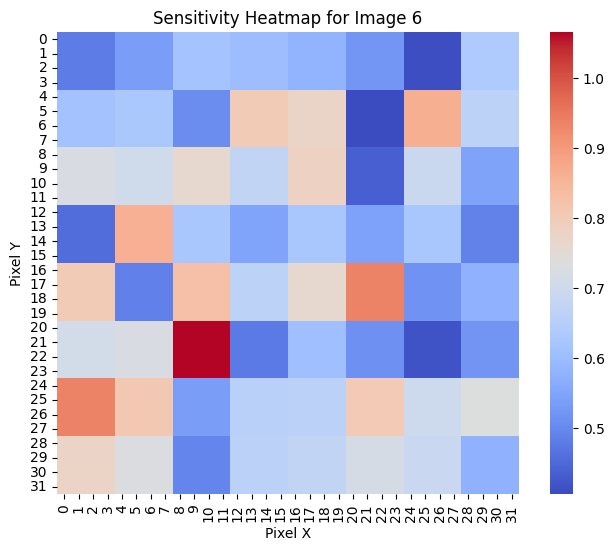}
        \includegraphics[width=0.092\linewidth]{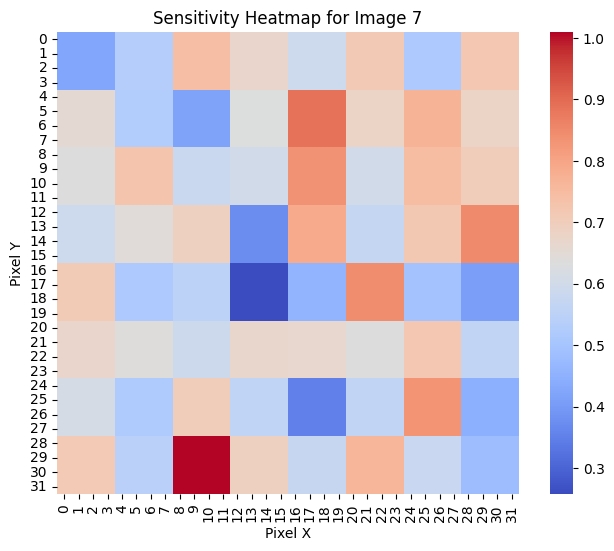}
        \includegraphics[width=0.092\linewidth]{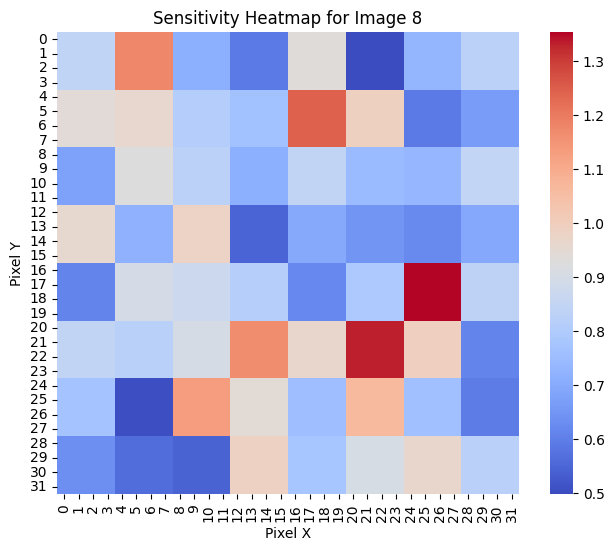}\\
%    \end{subfigure*}\\
%        \begin{subfigure*}{\textwidth}
        \includegraphics[width=0.092\linewidth]{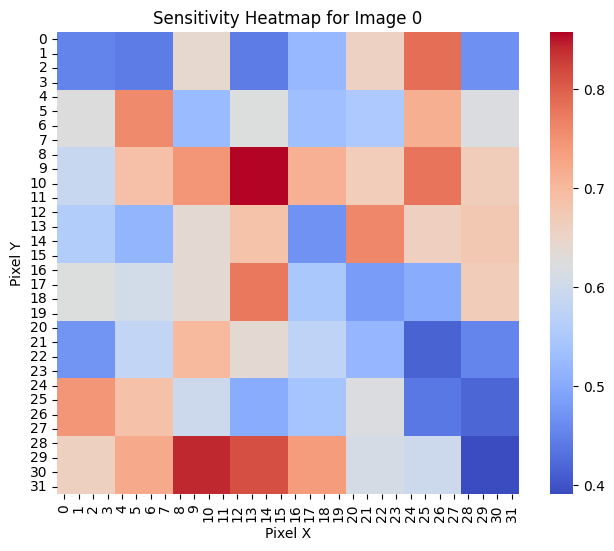}
        \includegraphics[width=0.092\linewidth]{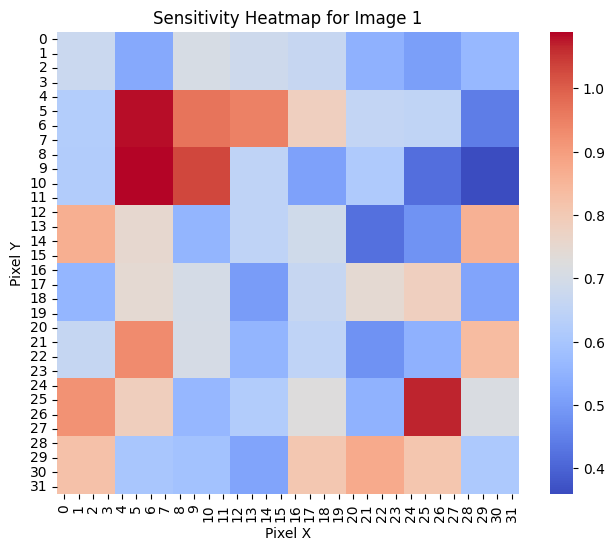}
        \includegraphics[width=0.092\linewidth]{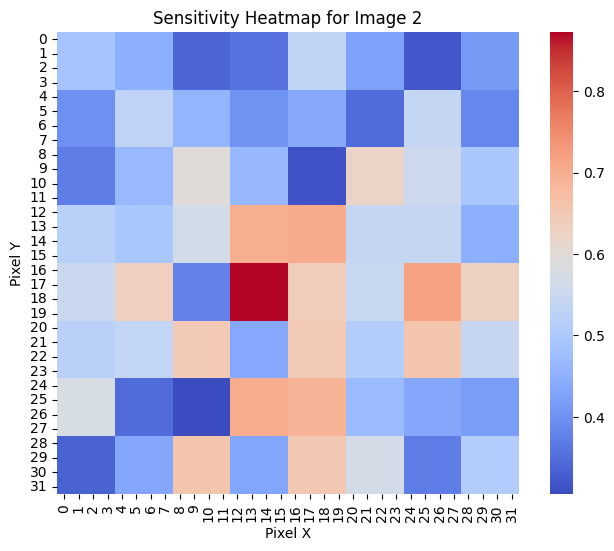}
        \includegraphics[width=0.092\linewidth]{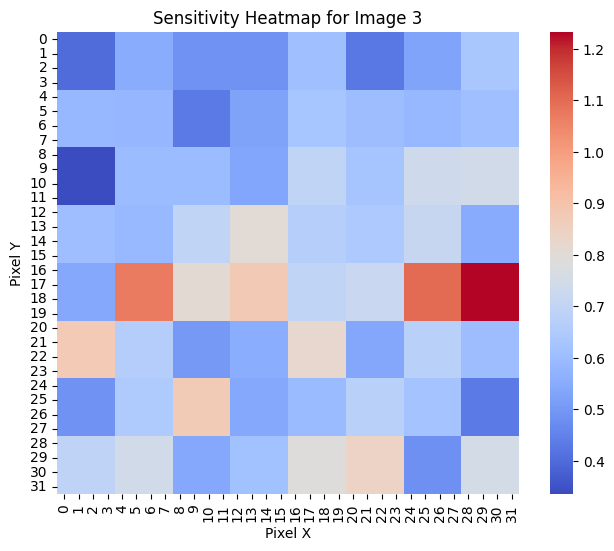}
        \includegraphics[width=0.092\linewidth]{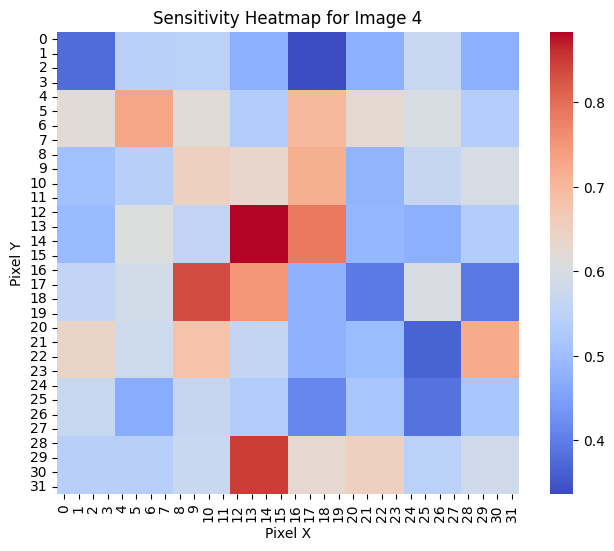}
        \includegraphics[width=0.092\linewidth]{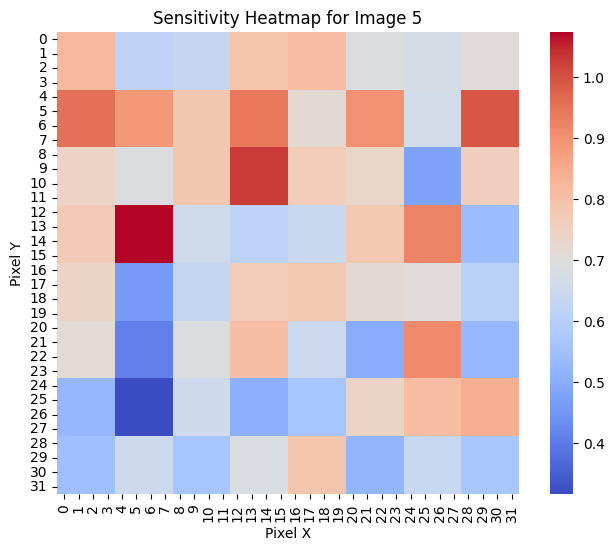}
        \includegraphics[width=0.092\linewidth]{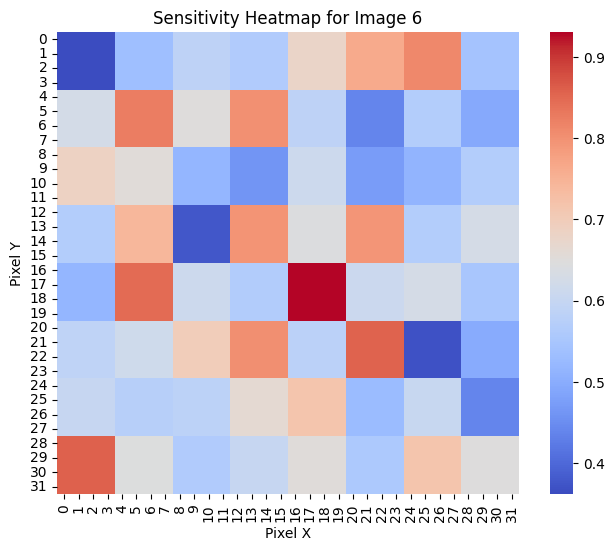}
        \includegraphics[width=0.092\linewidth]{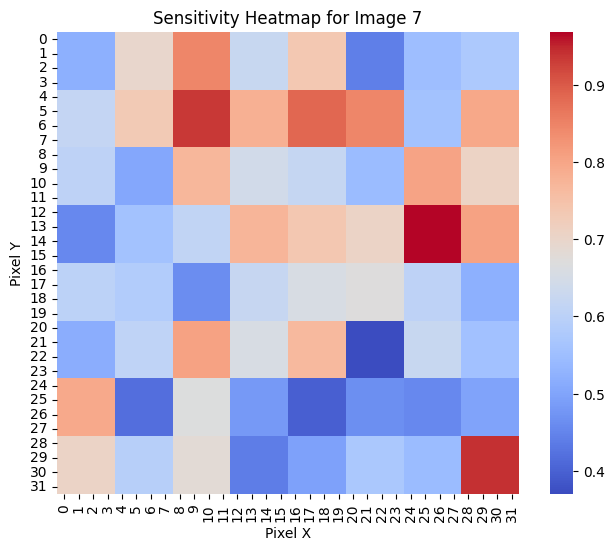}
        \includegraphics[width=0.092\linewidth]{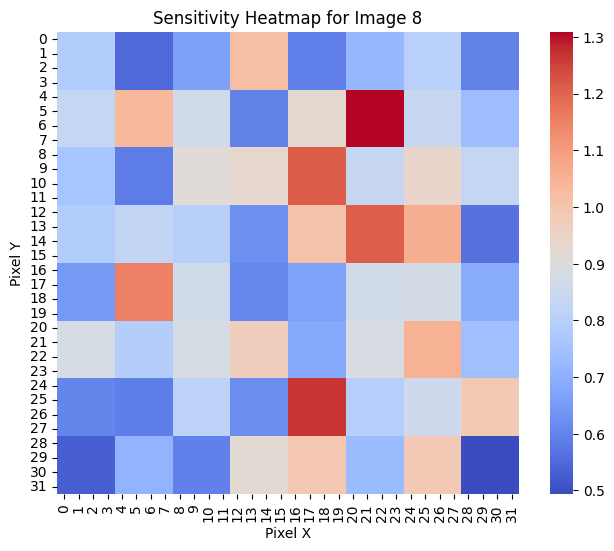}\\
%    \end{subfigure*}
    \vspace{0.2cm} 
    \caption{The sensitivity heatmap for the three color channels of block 1. The scale is $10^{0}$.}
    \label{fig:SensitHeatmapBlock1}
\end{figure}

We demonstrate the images with a fixed class number (zero class — airplanes) from the CIFAR-10 dataset. For ease of observation, we apply the pixelization to the sensitivity maps: each section represents a larger block of $4\times4$ pixels. In Fig. \ref{fig:SensitHeatmapBlock1} we show the heatmaps of the sensitivity for the red, green, and blue channels.

We can see that in the first block, the sensitivity of the different color channels varies significantly. It can be confidently assumed that the first block of the VGG-16 model is sensitive to color and is responsible for classifying and generalizing the color characteristics of the image.

\begin{figure}[htbp]
    \centering
%    \begin{subfigure}{\textwidth}
        \includegraphics[width=0.092\linewidth]{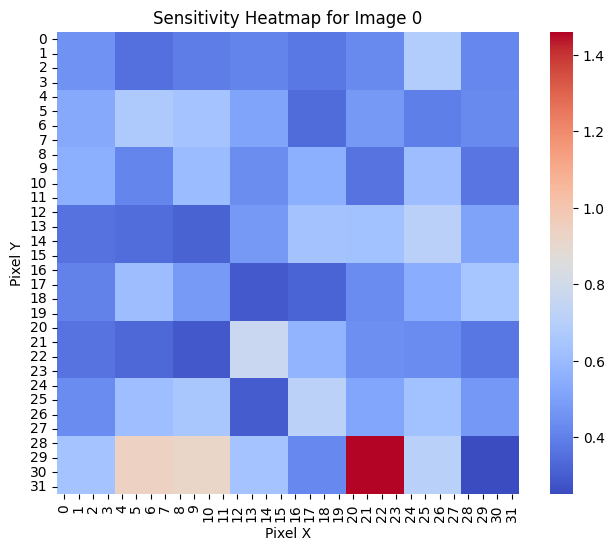}
        \includegraphics[width=0.092\linewidth]{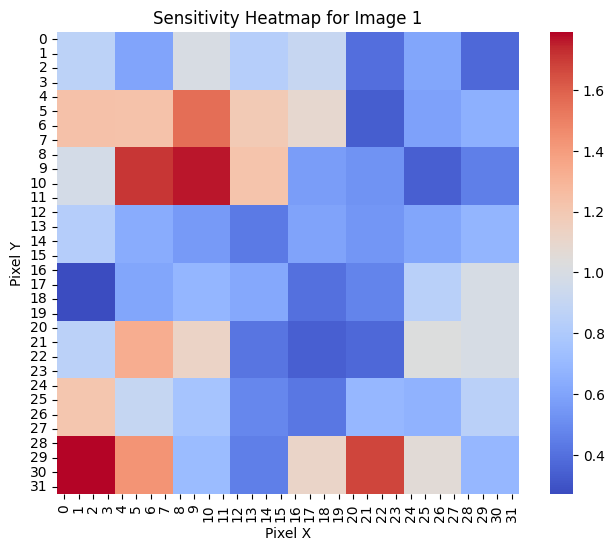}
        \includegraphics[width=0.092\linewidth]{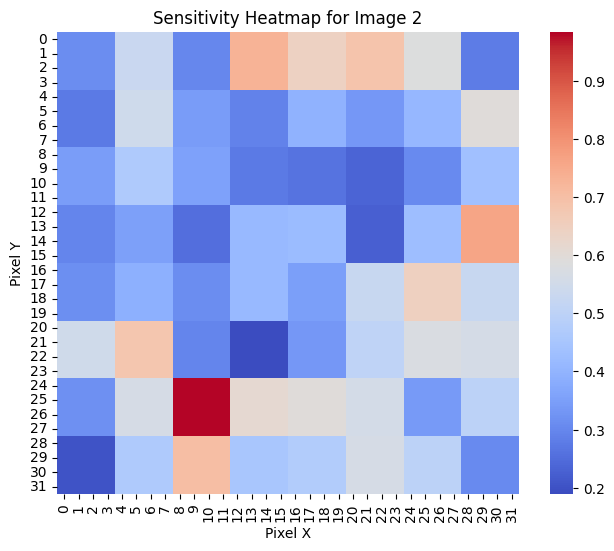}
        \includegraphics[width=0.092\linewidth]{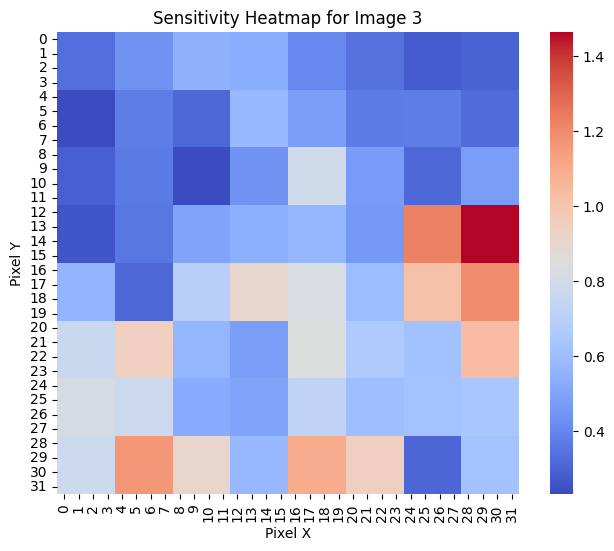}
        \includegraphics[width=0.092\linewidth]{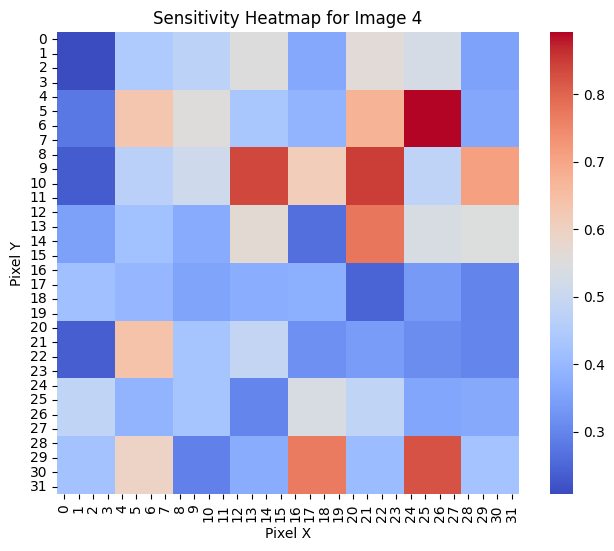}
        \includegraphics[width=0.092\linewidth]{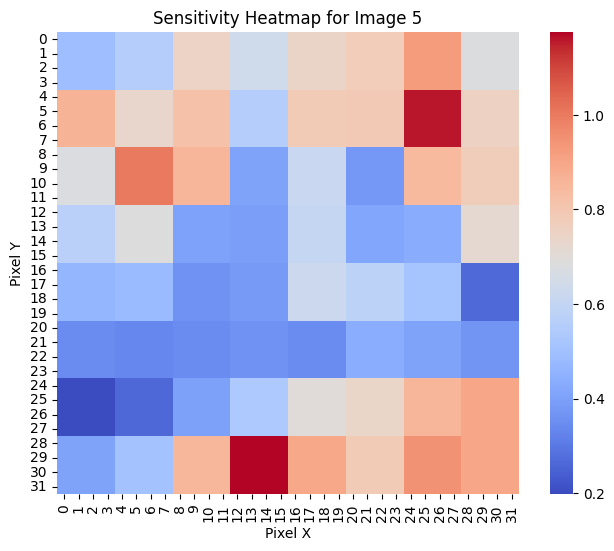}
        \includegraphics[width=0.092\linewidth]{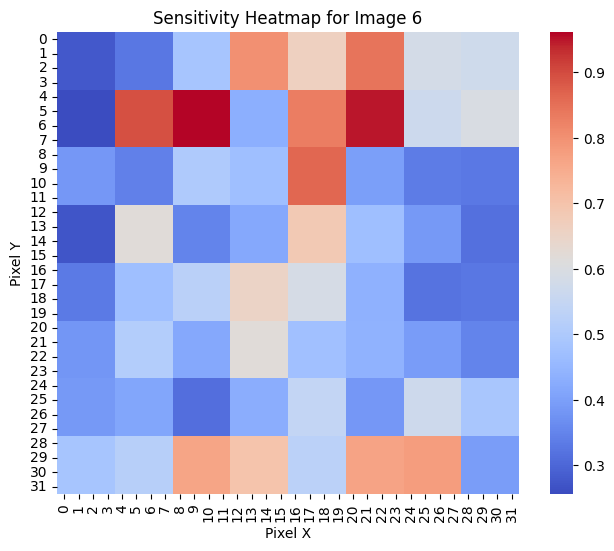}
        \includegraphics[width=0.092\linewidth]{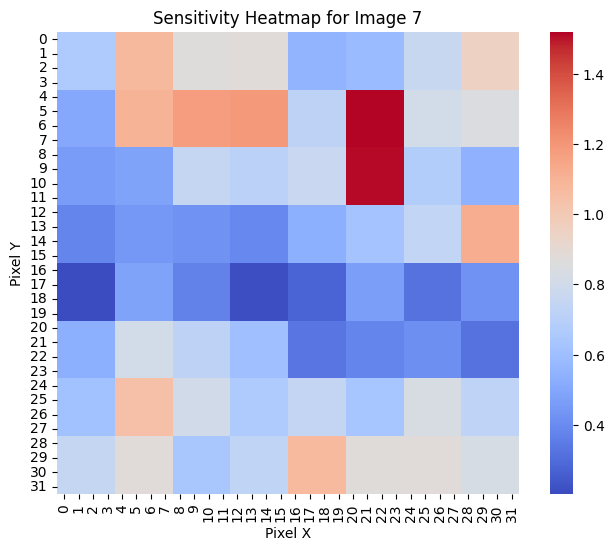}
        \includegraphics[width=0.092\linewidth]{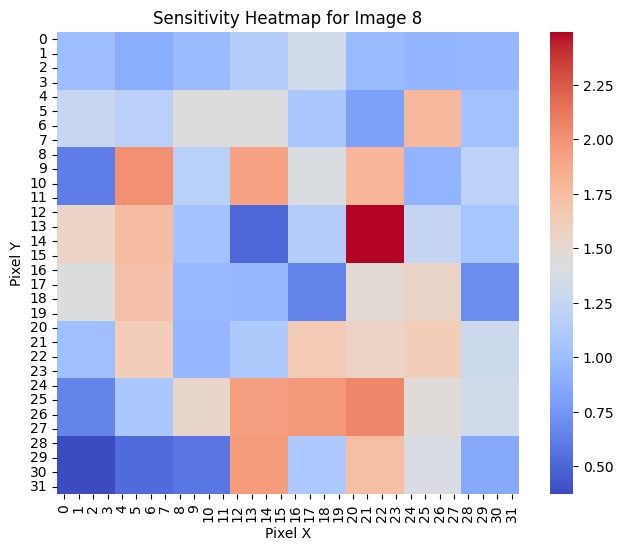}\\
%    \end{subfigure}
%        \begin{subfigure}{\textwidth}
        \includegraphics[width=0.092\linewidth]{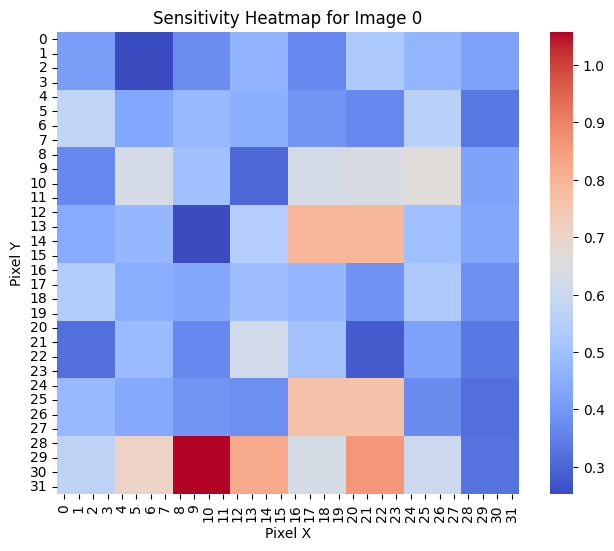}
        \includegraphics[width=0.092\linewidth]{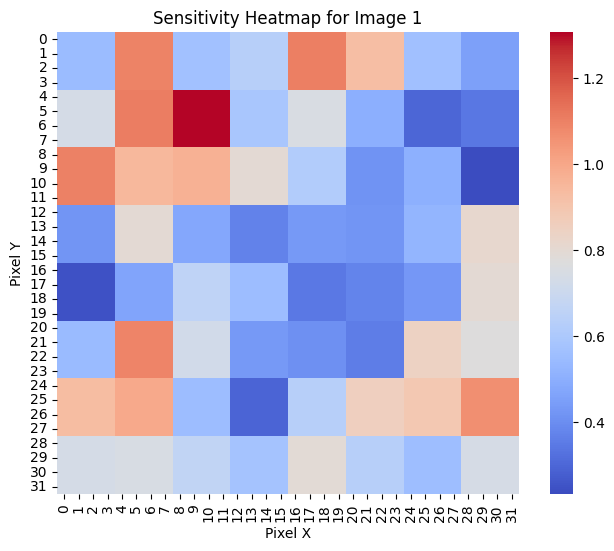}
        \includegraphics[width=0.092\linewidth]{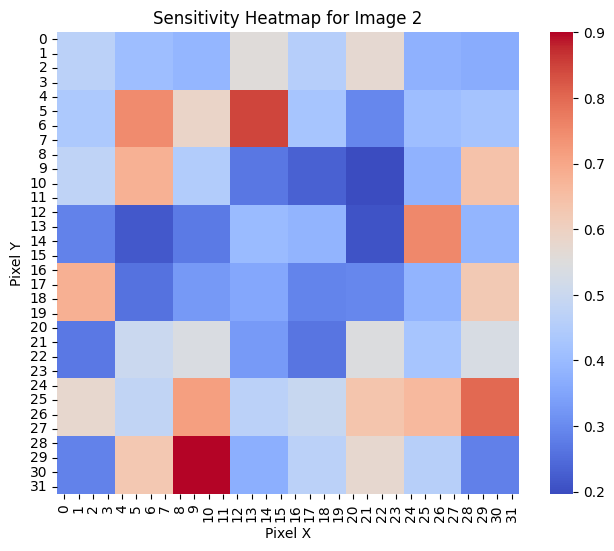}
        \includegraphics[width=0.092\linewidth]{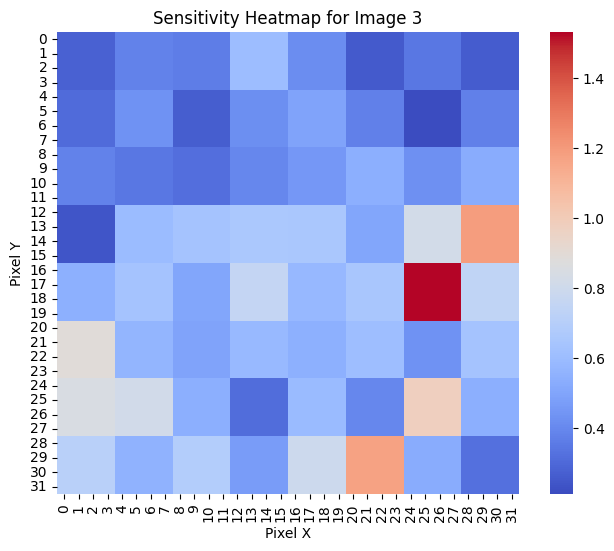}
        \includegraphics[width=0.092\linewidth]{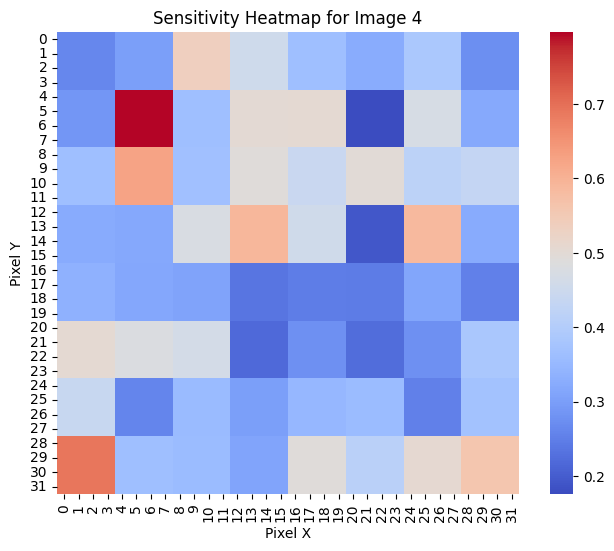}
        \includegraphics[width=0.092\linewidth]{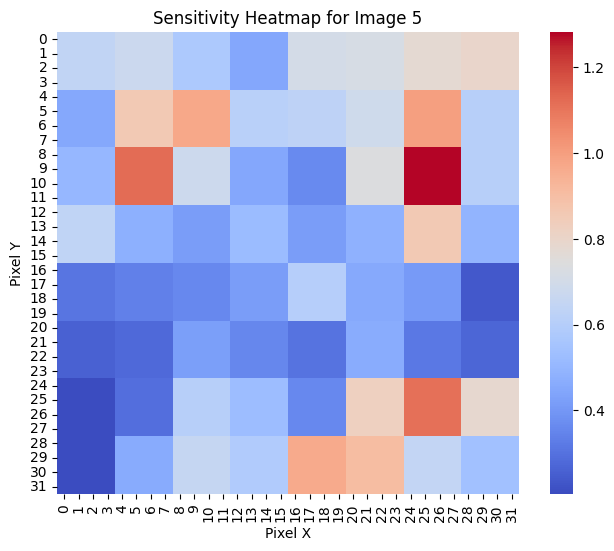}
        \includegraphics[width=0.092\linewidth]{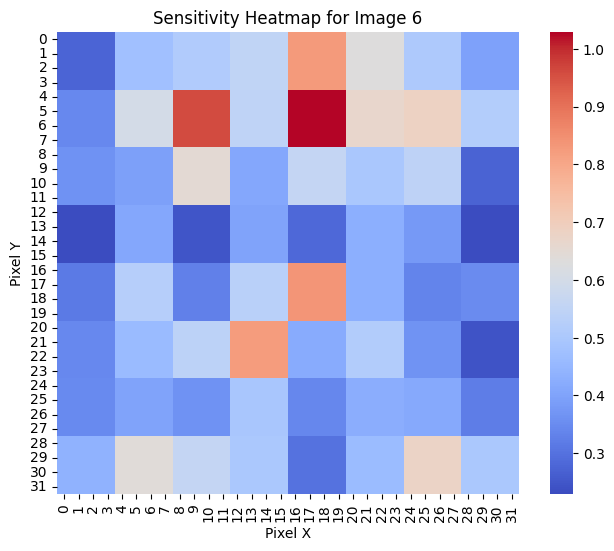}
        \includegraphics[width=0.092\linewidth]{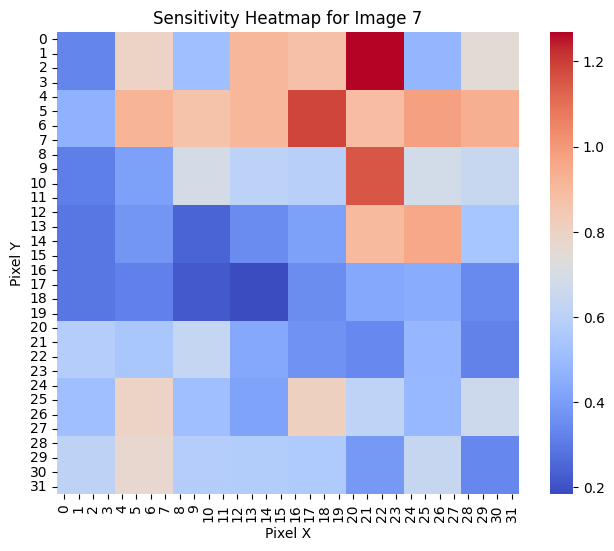}
        \includegraphics[width=0.092\linewidth]{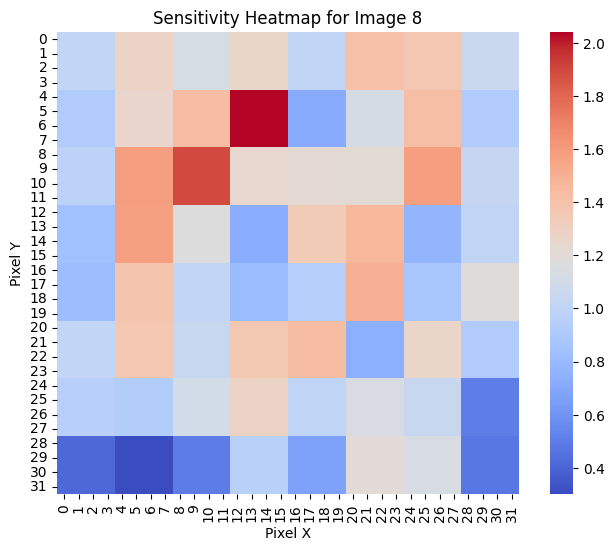}
%    \end{subfigure}
\\
%    \begin{subfigure}{\textwidth}
        \includegraphics[width=0.092\linewidth]{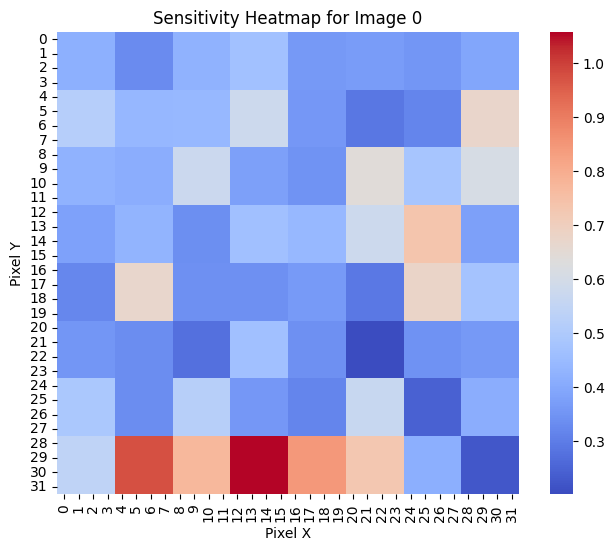}
        \includegraphics[width=0.092\linewidth]{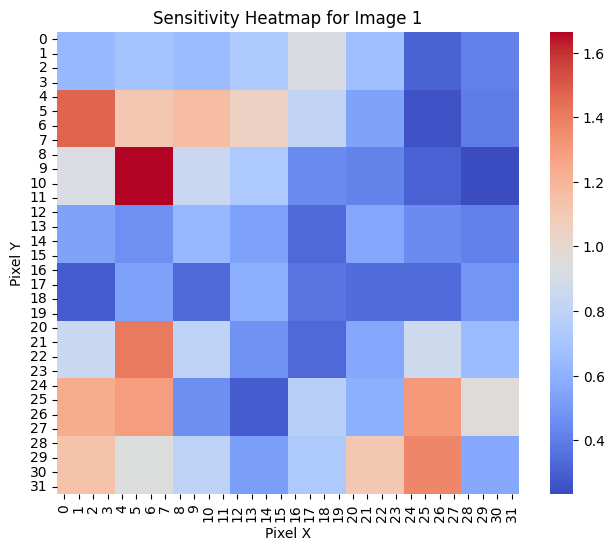}
        \includegraphics[width=0.092\linewidth]{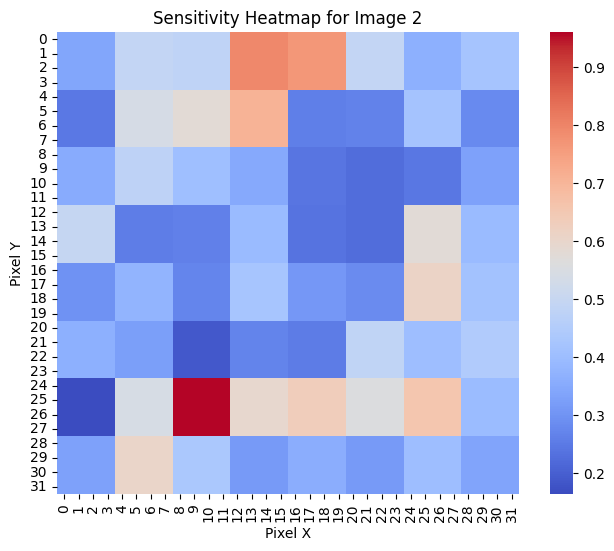}
        \includegraphics[width=0.092\linewidth]{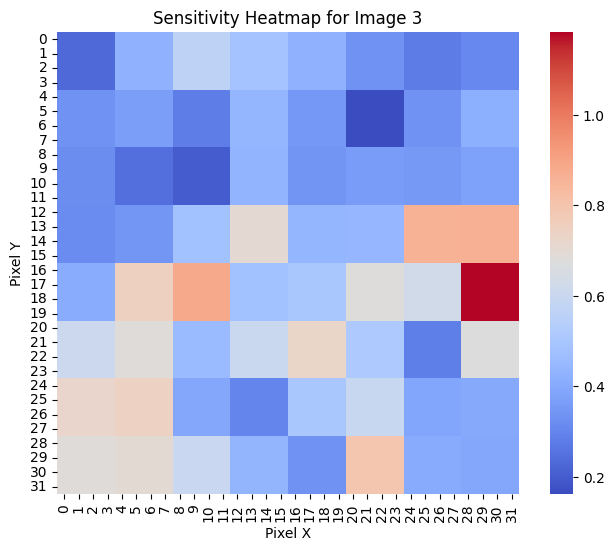}
        \includegraphics[width=0.092\linewidth]{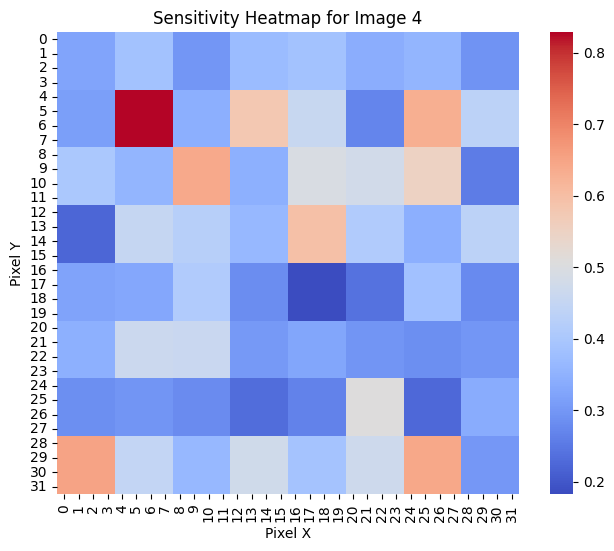}
        \includegraphics[width=0.092\linewidth]{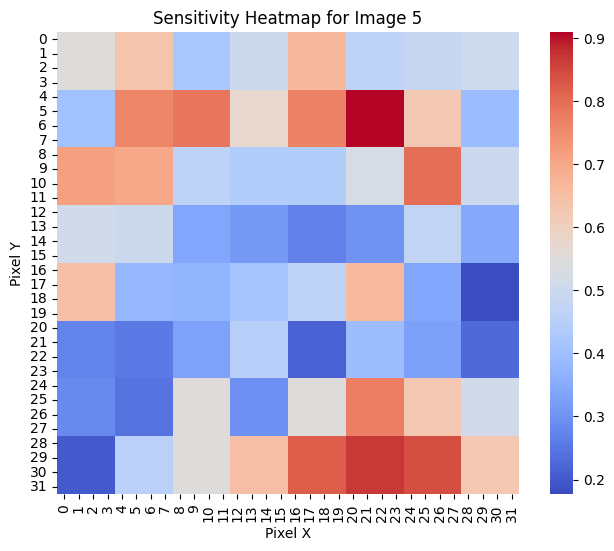}
        \includegraphics[width=0.092\linewidth]{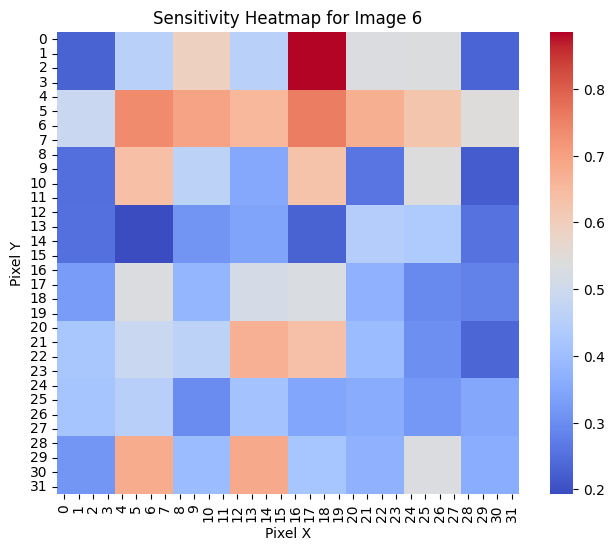}
        \includegraphics[width=0.092\linewidth]{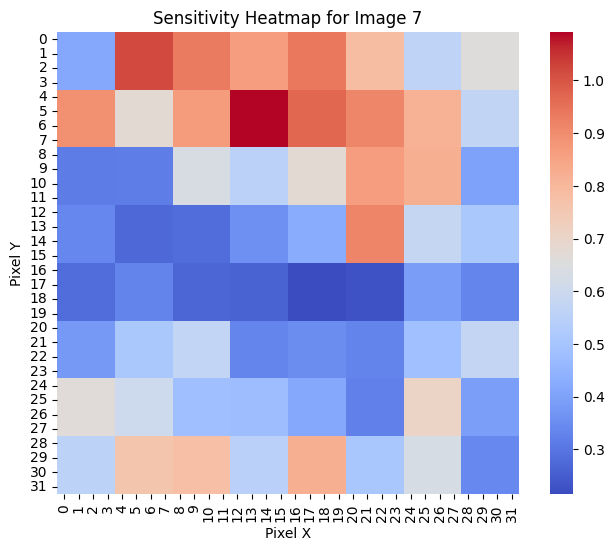}
        \includegraphics[width=0.092\linewidth]{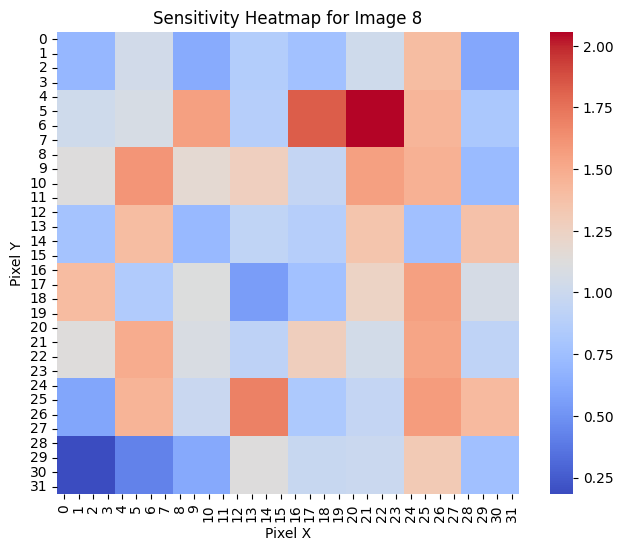}
%    \end{subfigure}
    \vspace{0.2cm} 
    \caption{The sensitivity heatmap for the three color channels of block 2. The scale is $10^{0}$.}
    \label{fig:SensitHeatmapBlock2}
\end{figure}

\begin{figure}[htbp]
    \centering  
%    \begin{subfigure}{\textwidth}
        \includegraphics[width=0.092\linewidth]{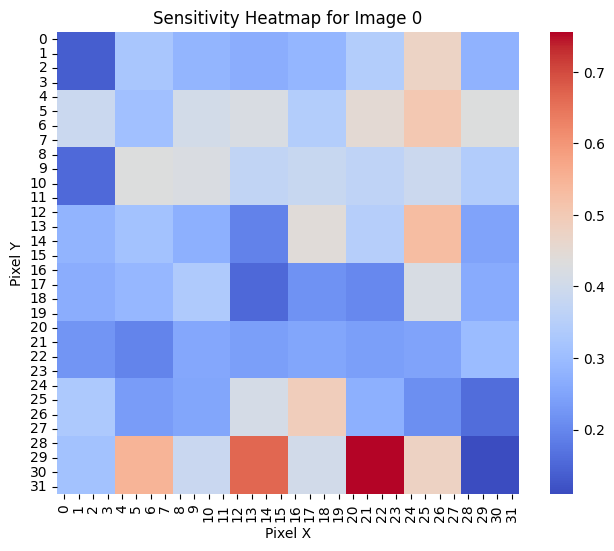}
        \includegraphics[width=0.092\linewidth]{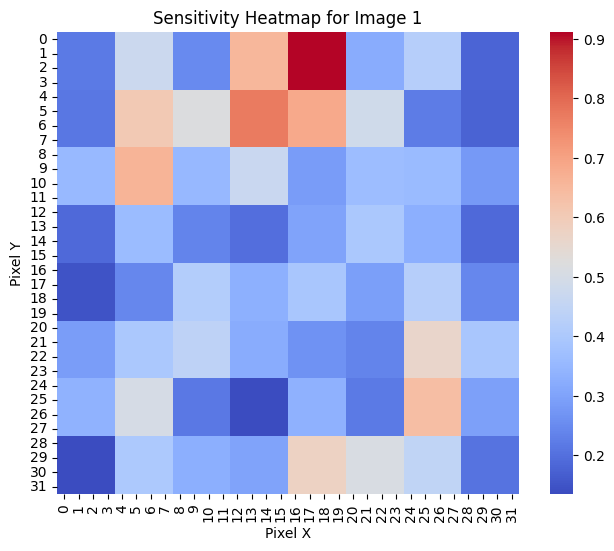}
        \includegraphics[width=0.092\linewidth]{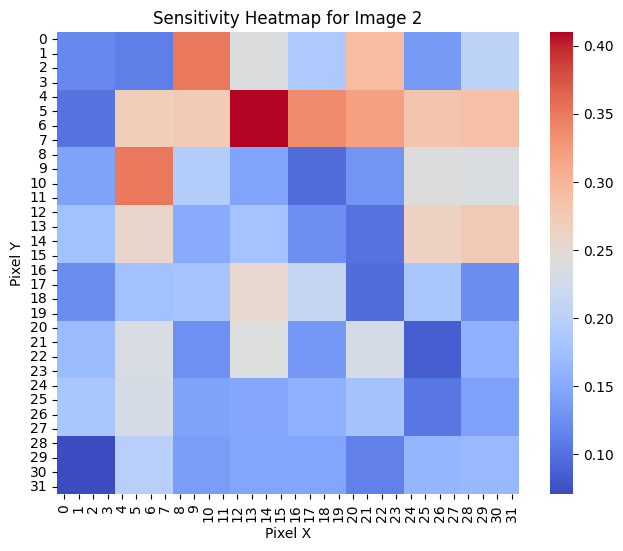}
        \includegraphics[width=0.092\linewidth]{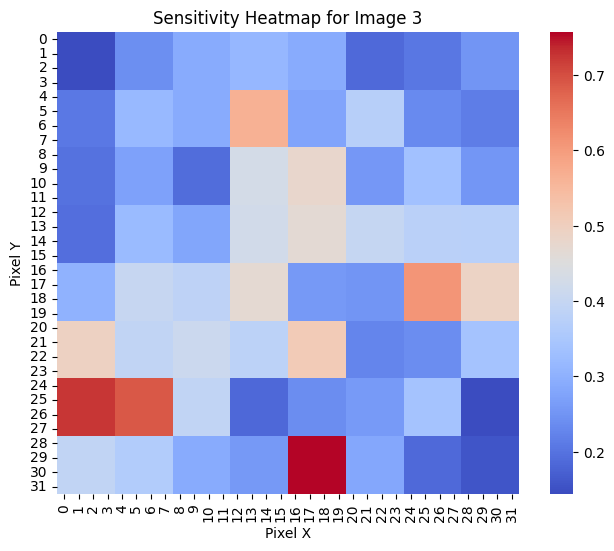}
        \includegraphics[width=0.092\linewidth]{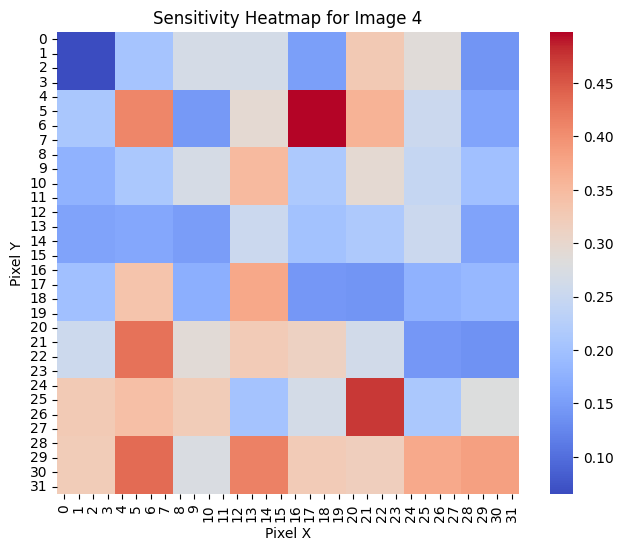}
        \includegraphics[width=0.092\linewidth]{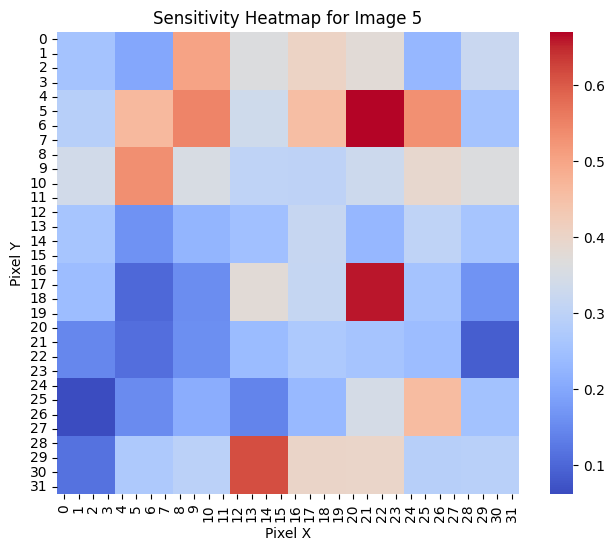}
        \includegraphics[width=0.092\linewidth]{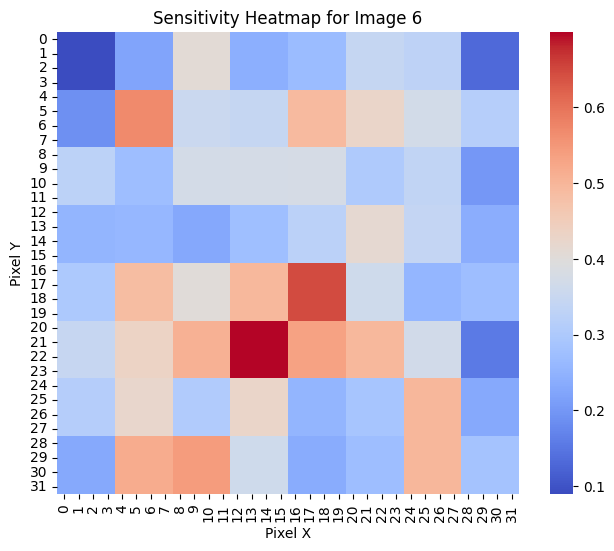}
        \includegraphics[width=0.092\linewidth]{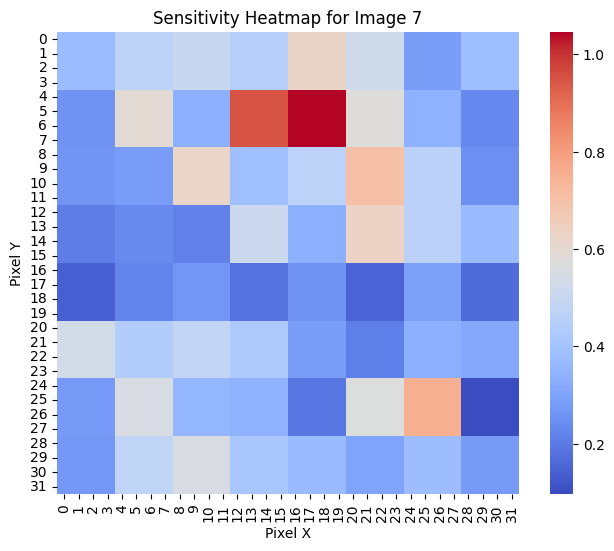}
        \includegraphics[width=0.092\linewidth]{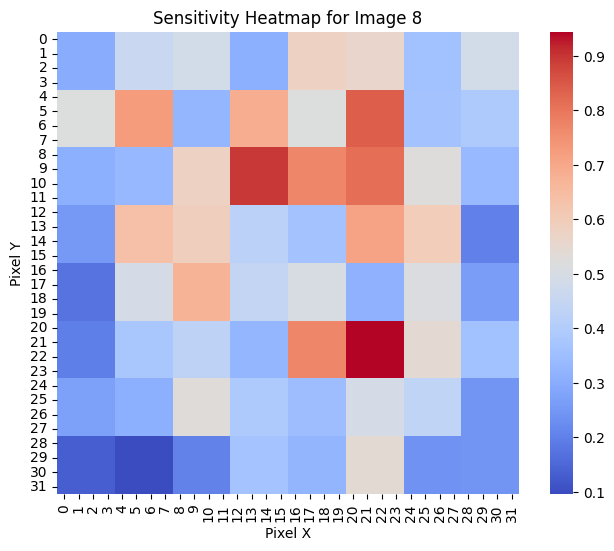}
%    \end{subfigure}
\\
%        \begin{subfigure}{\linewidth}
        \includegraphics[width=0.092\linewidth]{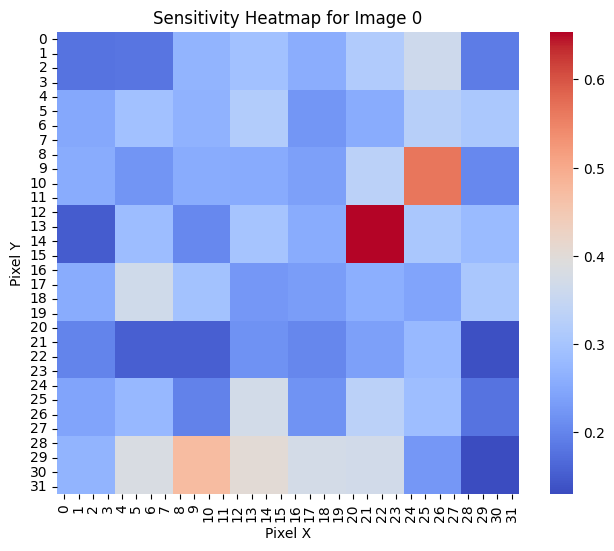}
        \includegraphics[width=0.092\linewidth]{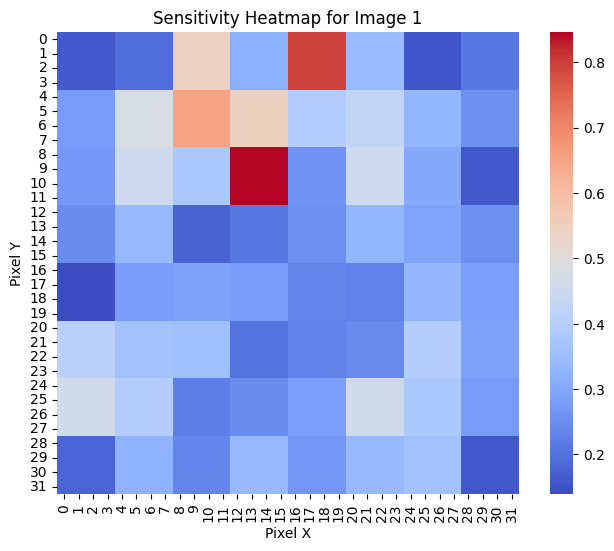}
        \includegraphics[width=0.092\linewidth]{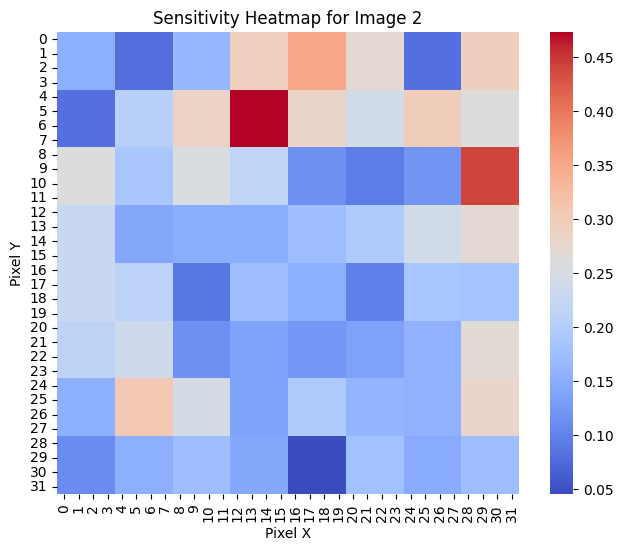}
        \includegraphics[width=0.092\linewidth]{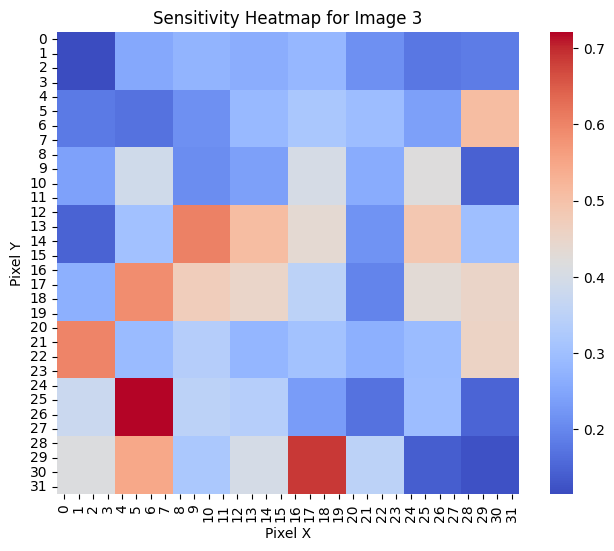}
        \includegraphics[width=0.092\linewidth]{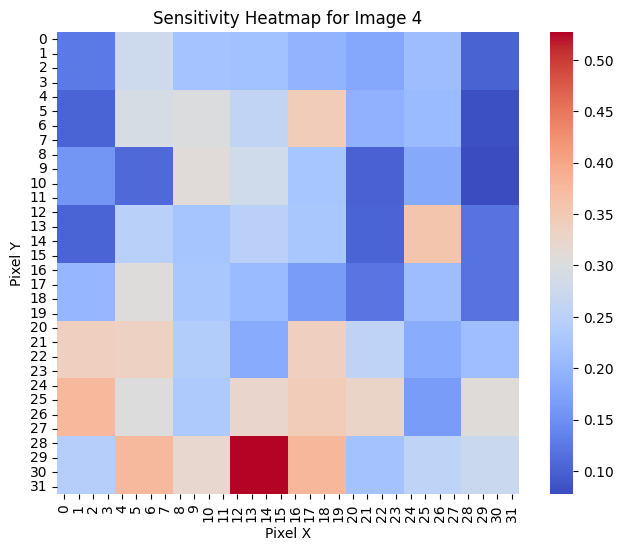}
        \includegraphics[width=0.092\linewidth]{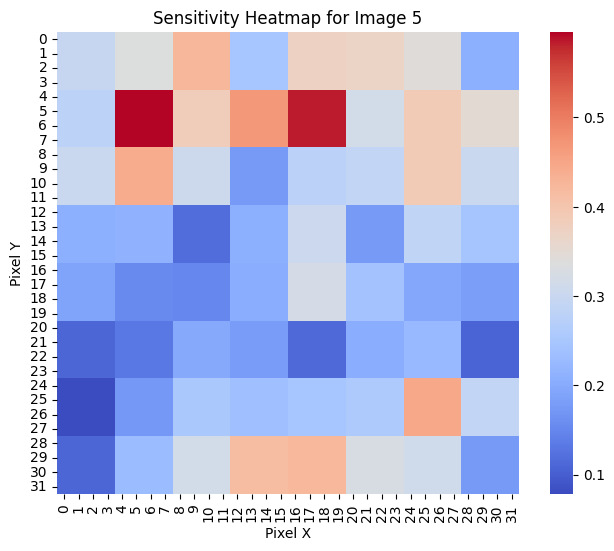}
        \includegraphics[width=0.092\linewidth]{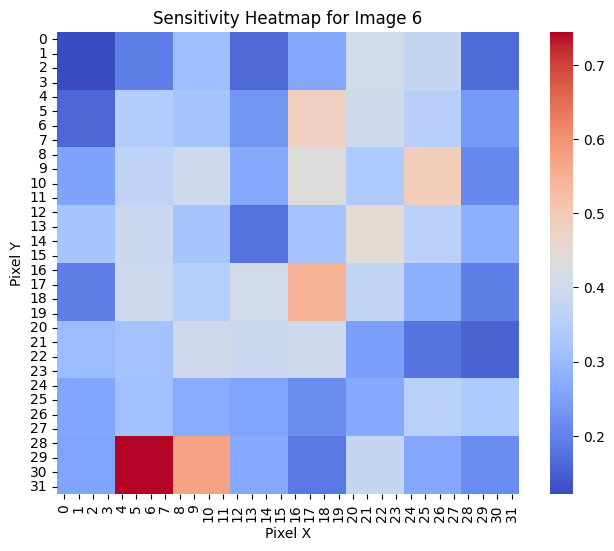}
        \includegraphics[width=0.092\linewidth]{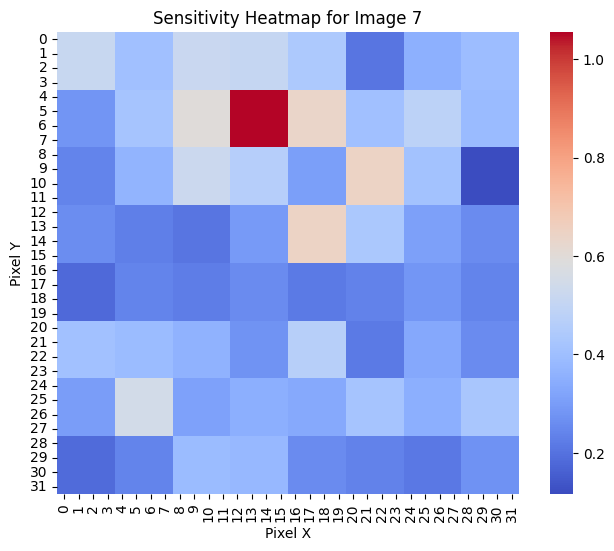}
        \includegraphics[width=0.092\linewidth]{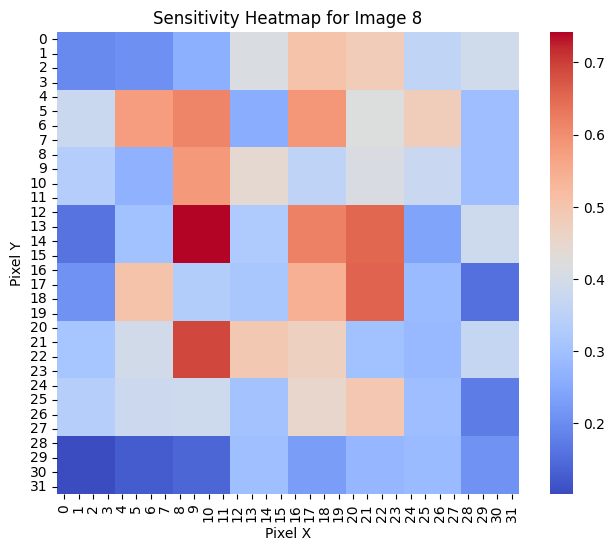}
%    \end{subfigure}
\\
%        \begin{subfigure}{\linewidth}
        \includegraphics[width=0.092\linewidth]{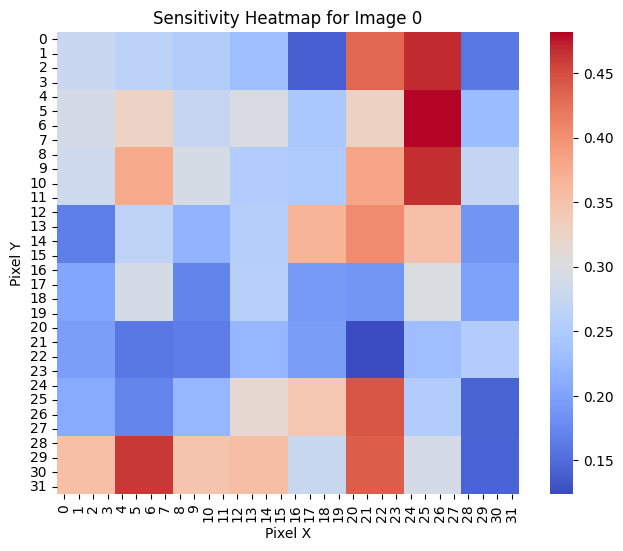}
        \includegraphics[width=0.092\linewidth]{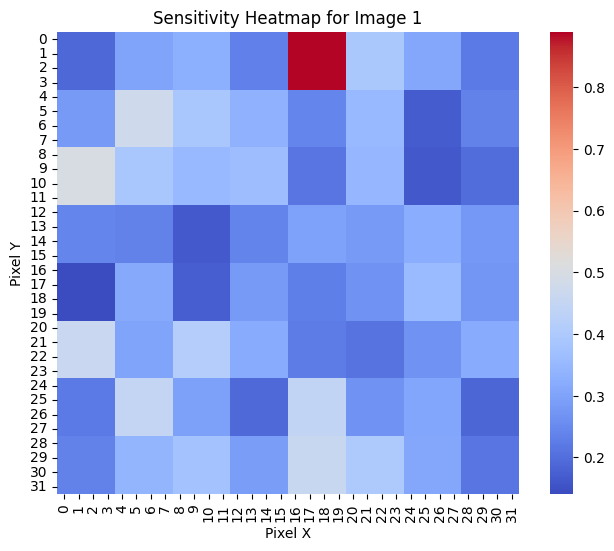}
        \includegraphics[width=0.092\linewidth]{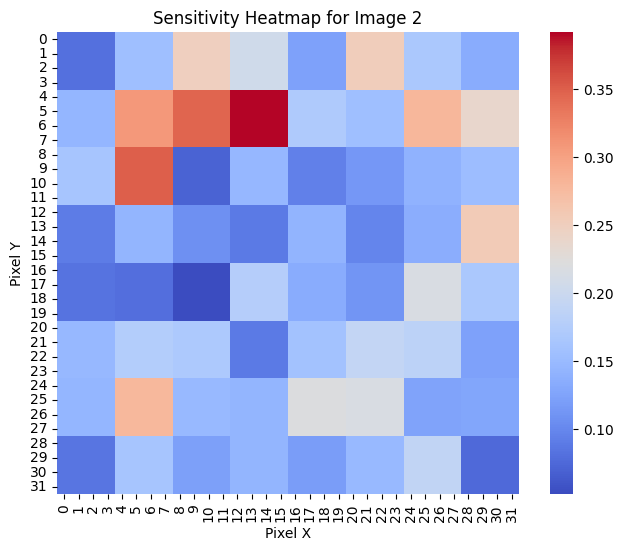}
        \includegraphics[width=0.092\linewidth]{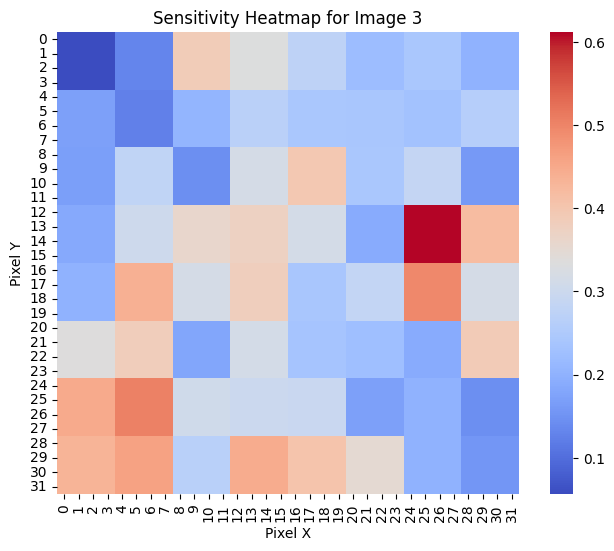}
        \includegraphics[width=0.092\linewidth]{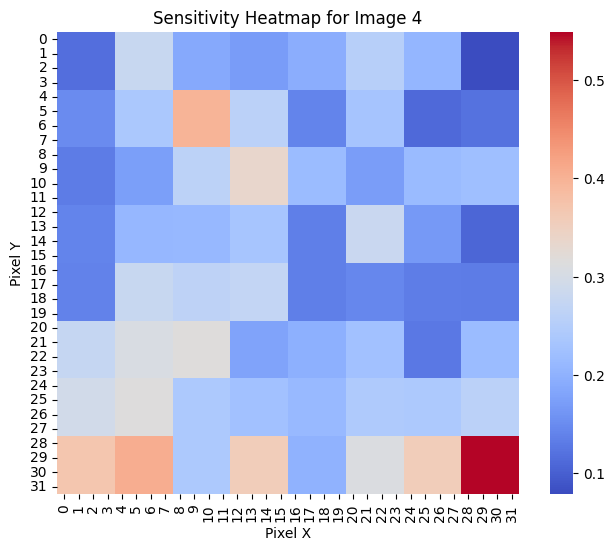}
        \includegraphics[width=0.092\linewidth]{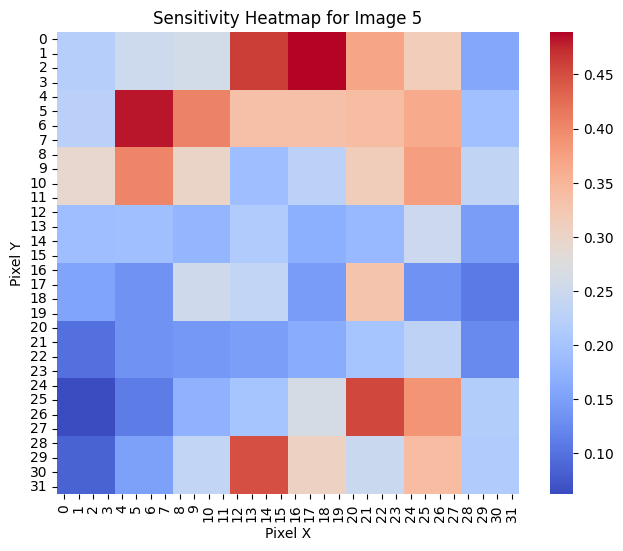}
        \includegraphics[width=0.092\linewidth]{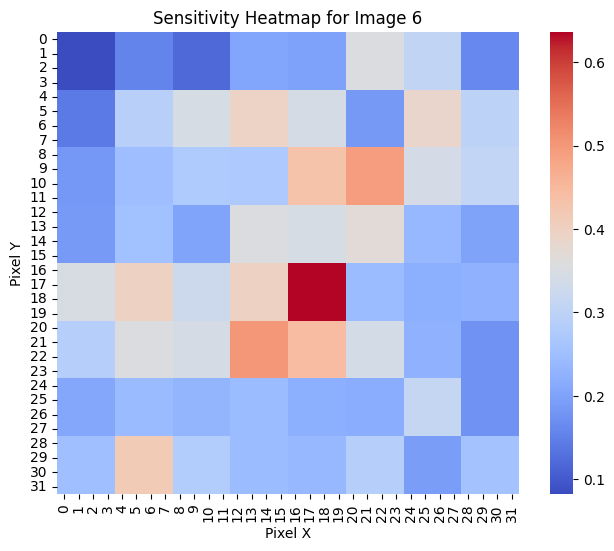}
        \includegraphics[width=0.092\linewidth]{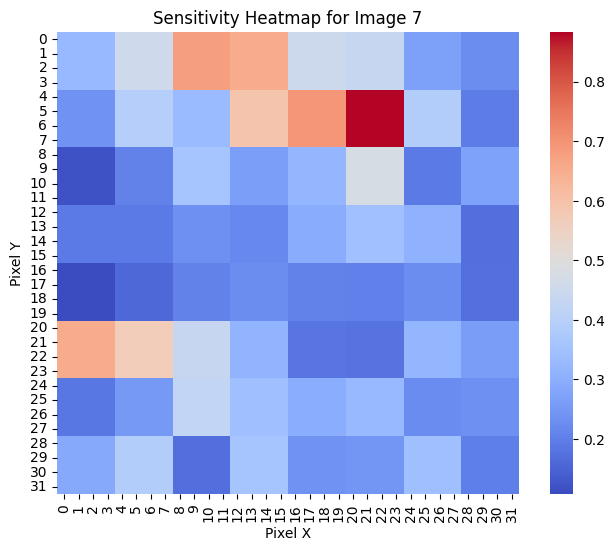}
        \includegraphics[width=0.092\linewidth]{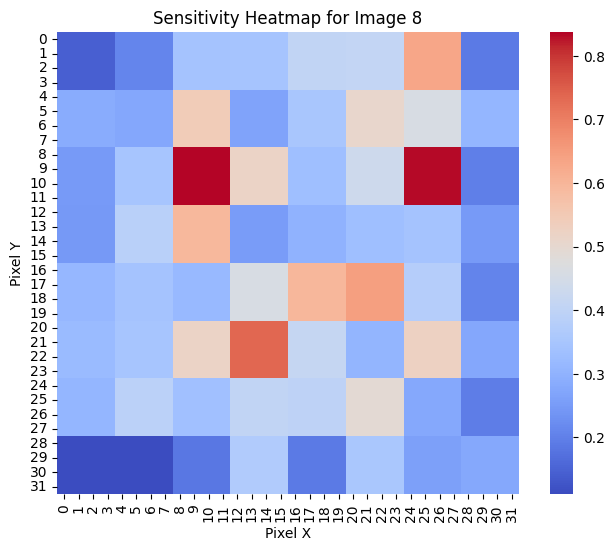}
%    \end{subfigure}
    \vspace{0.2cm} 
    \caption{The sensitivity heatmap for the three color channels of block 3. The scale is $10^{-1}$.}
    \label{fig:SensitHeatmapBlock3}
\end{figure}

\begin{figure}[htbp]
    \centering
%    \begin{subfigure}{\linewidth}
        \includegraphics[width=0.092\linewidth]{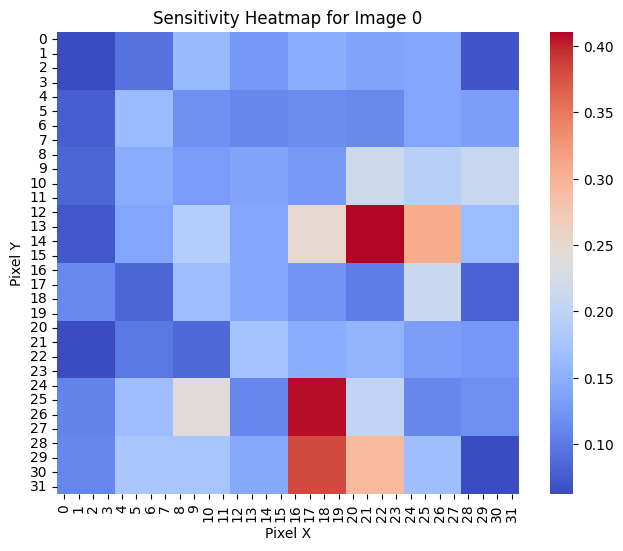}
        \includegraphics[width=0.092\linewidth]{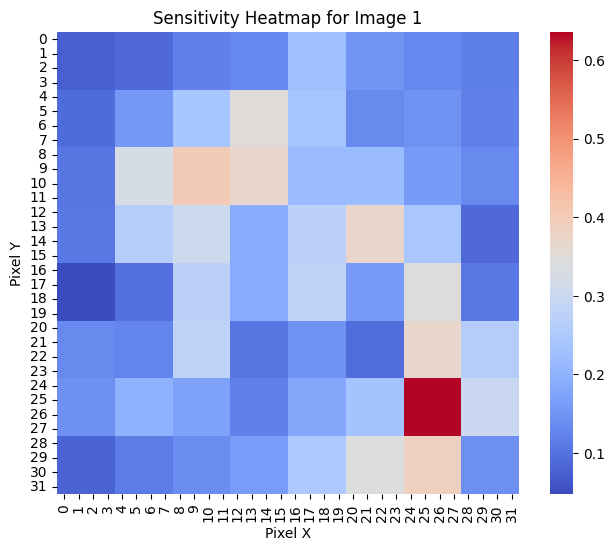}
        \includegraphics[width=0.092\linewidth]{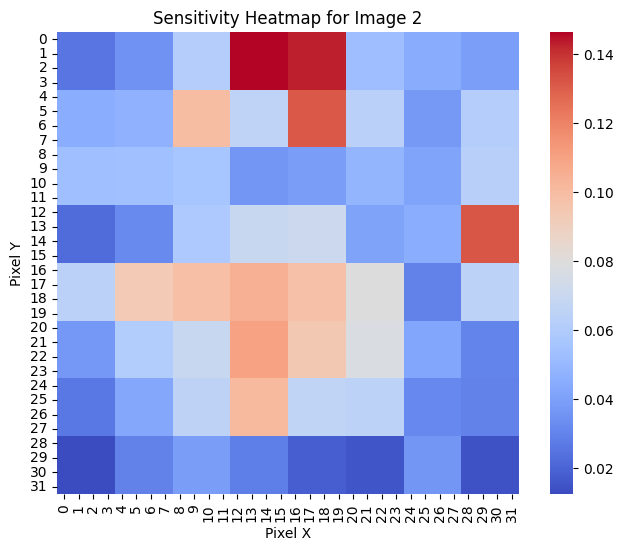}
        \includegraphics[width=0.092\linewidth]{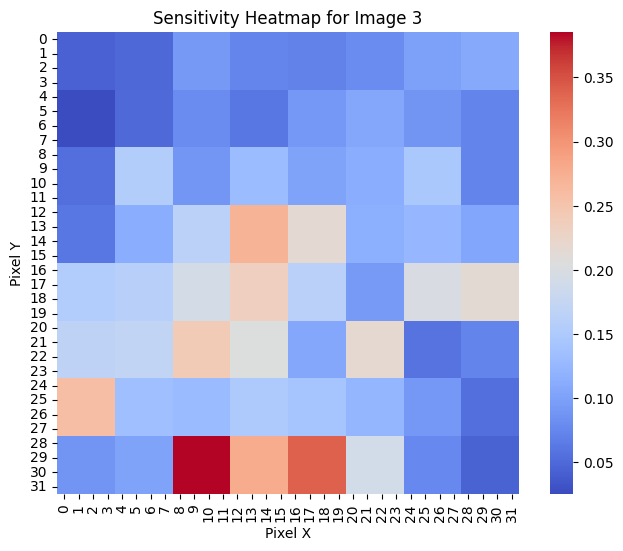}
        \includegraphics[width=0.092\linewidth]{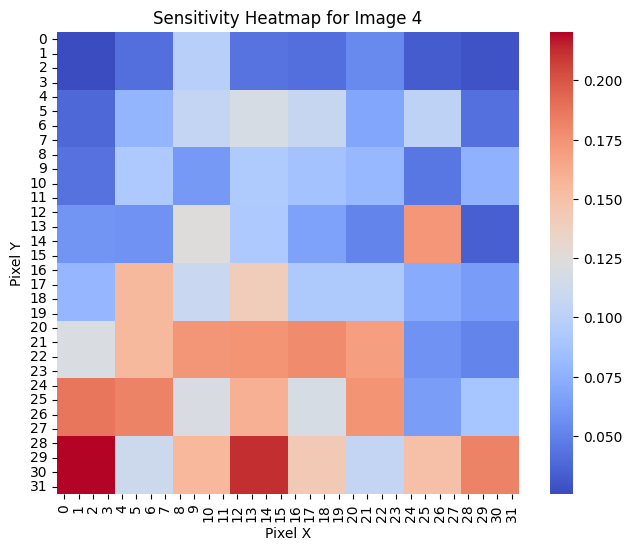}
        \includegraphics[width=0.092\linewidth]{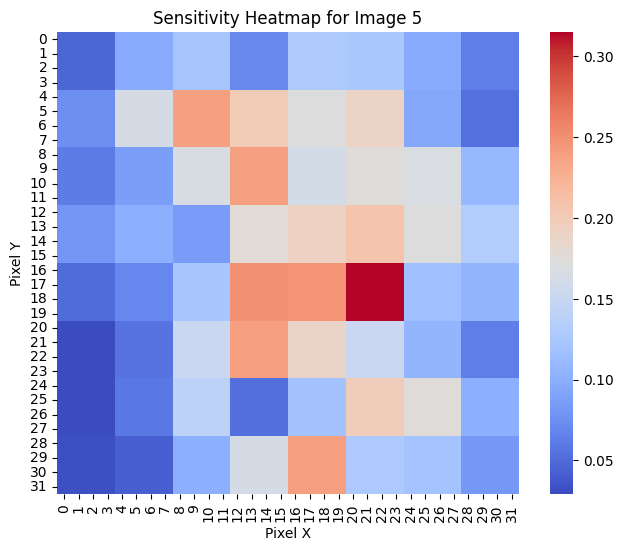}
        \includegraphics[width=0.092\linewidth]{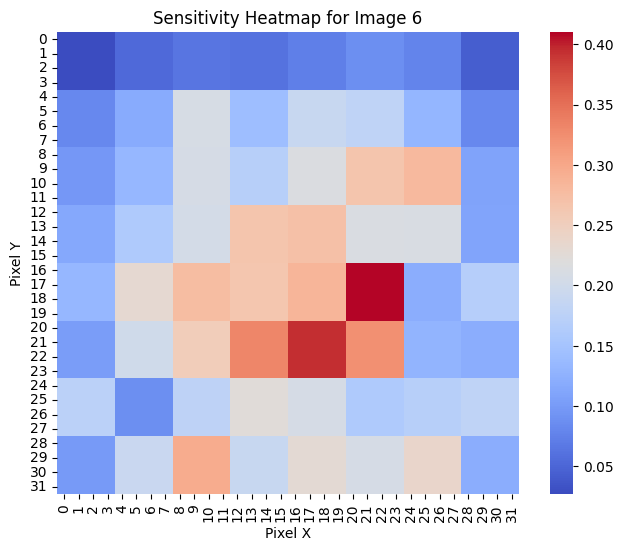}
        \includegraphics[width=0.092\linewidth]{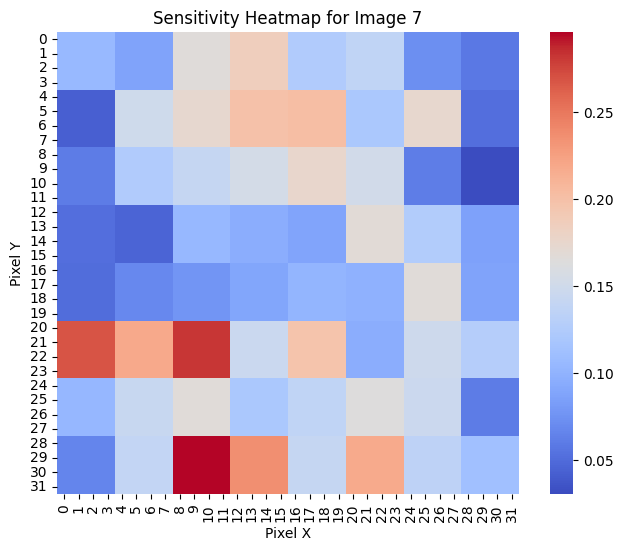}
        \includegraphics[width=0.092\linewidth]{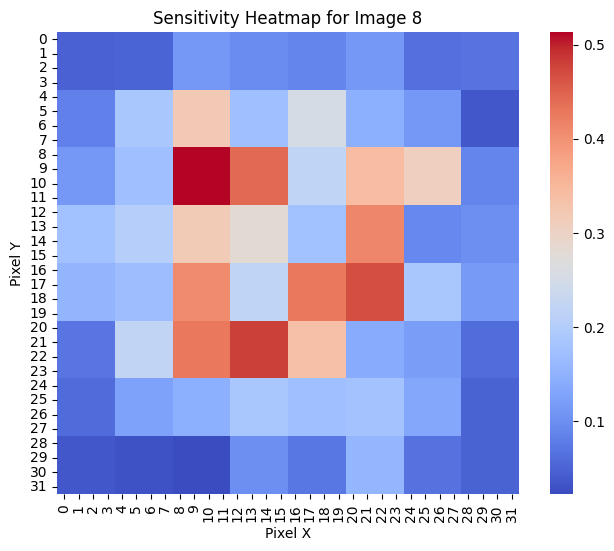}
%    \end{subfigure}
\\
%        \begin{subfigure}{\linewidth}
        \includegraphics[width=0.092\linewidth]{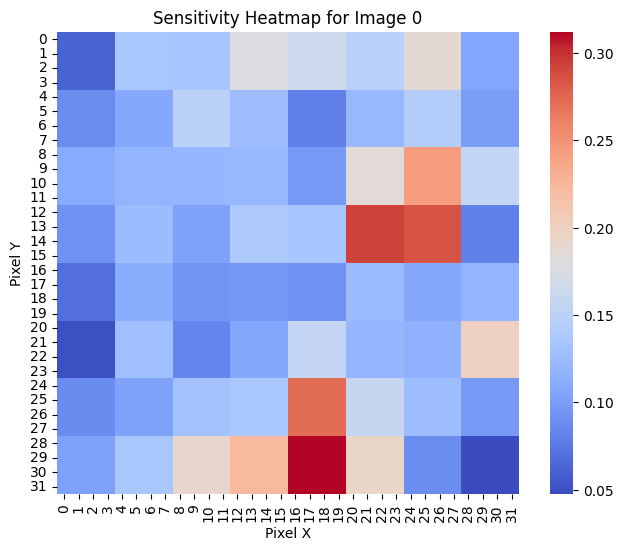}
        \includegraphics[width=0.092\linewidth]{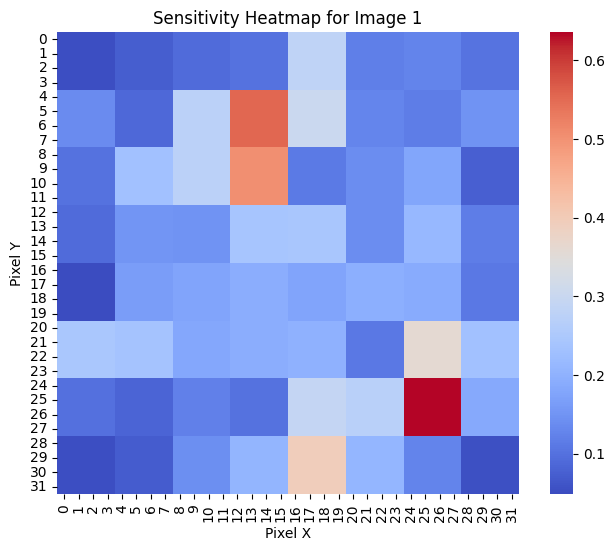}
        \includegraphics[width=0.092\linewidth]{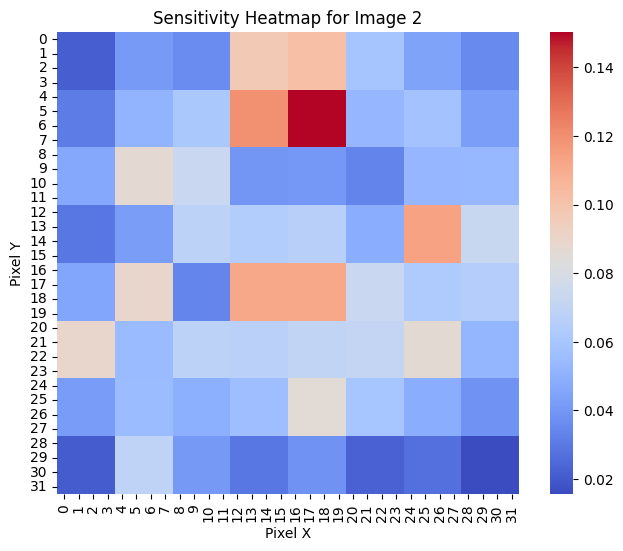}
        \includegraphics[width=0.092\linewidth]{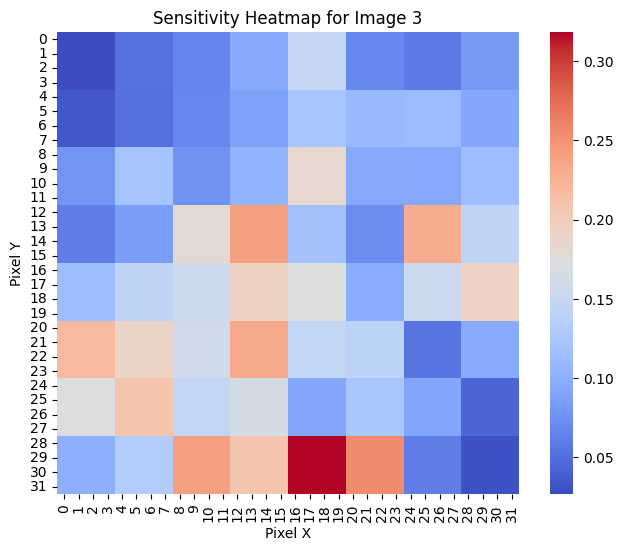}
        \includegraphics[width=0.092\linewidth]{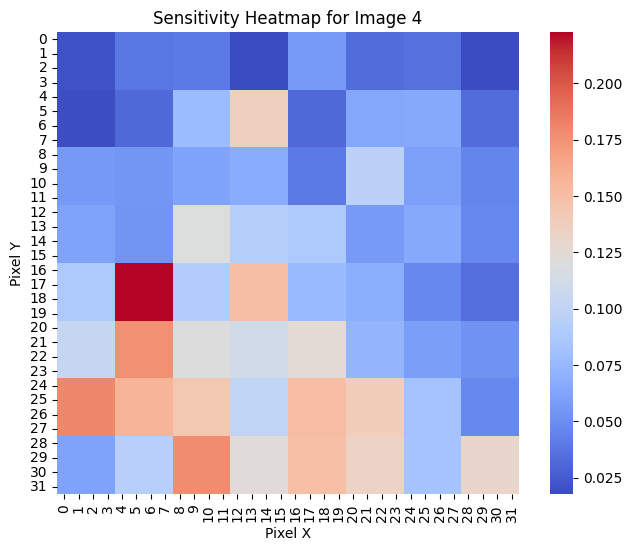}
        \includegraphics[width=0.092\linewidth]{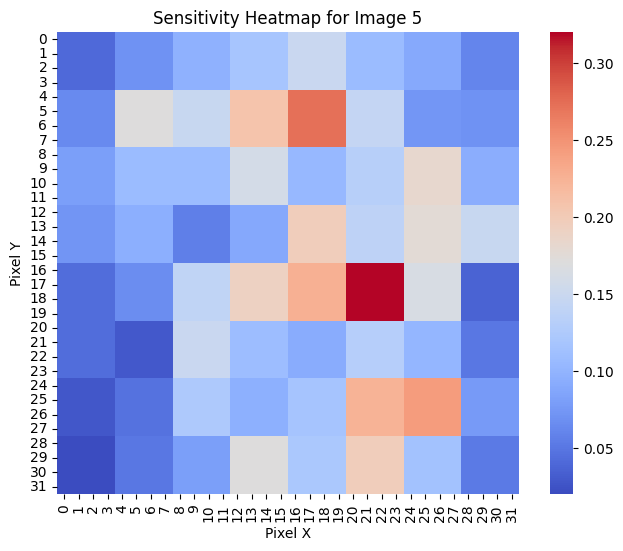}
        \includegraphics[width=0.092\linewidth]{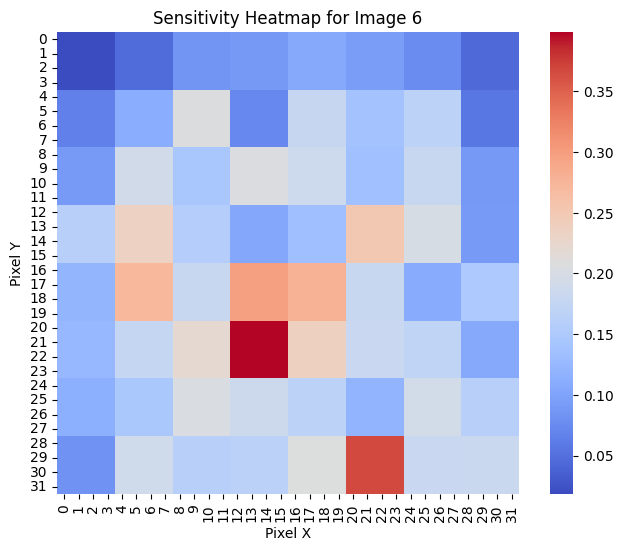}
        \includegraphics[width=0.092\linewidth]{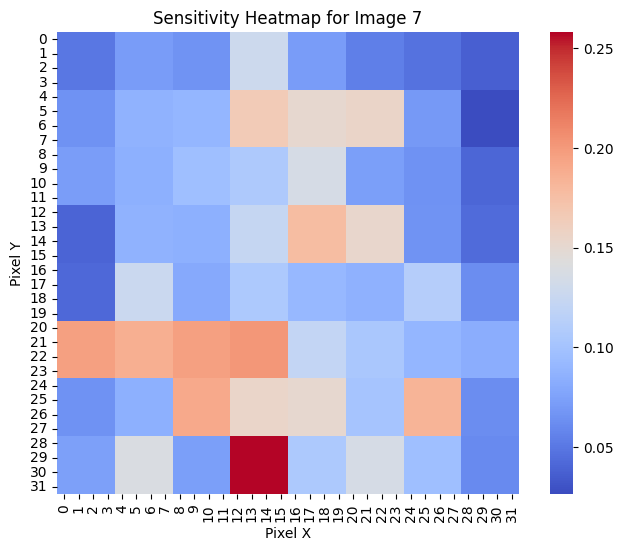}
        \includegraphics[width=0.092\linewidth]{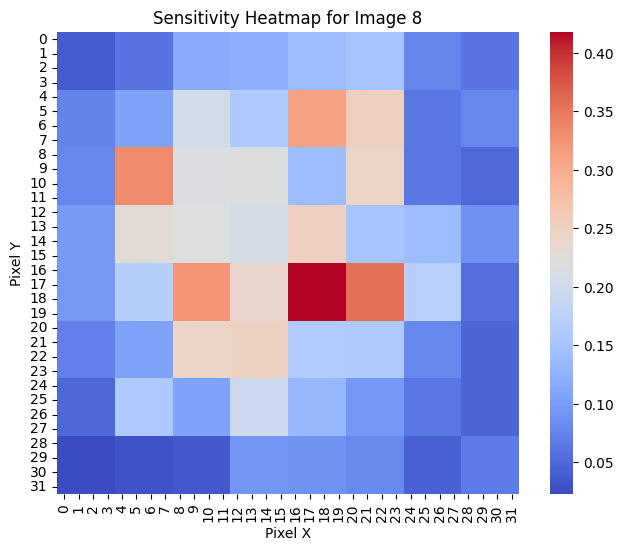}
%    \end{subfigure}
\\
%        \begin{subfigure}{\linewidth}
        \includegraphics[width=0.092\linewidth]{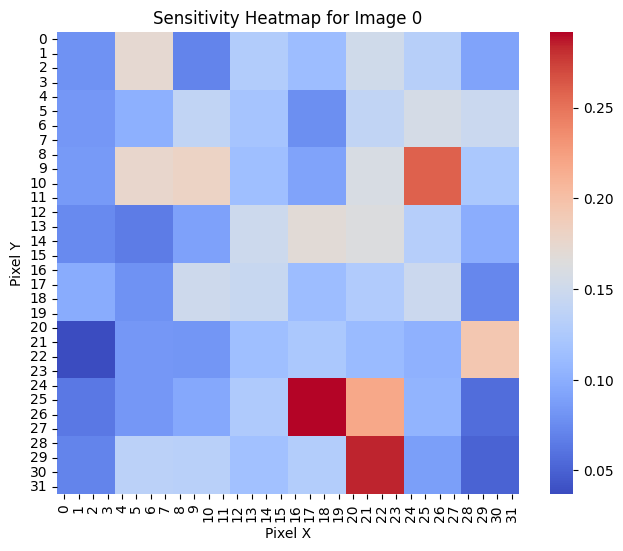}
        \includegraphics[width=0.092\linewidth]{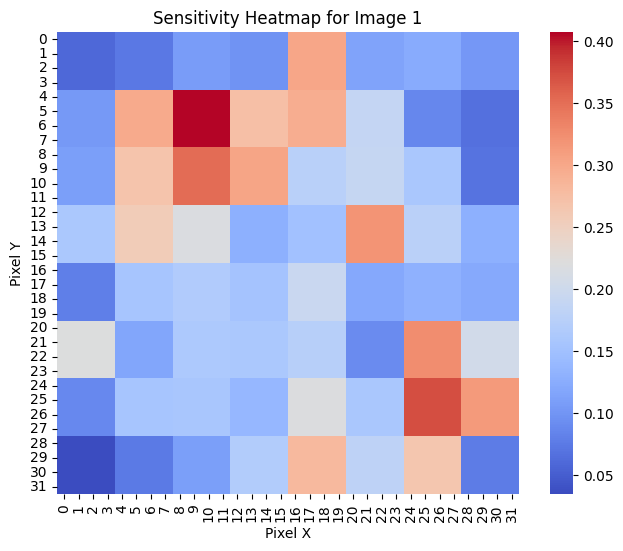}
        \includegraphics[width=0.092\linewidth]{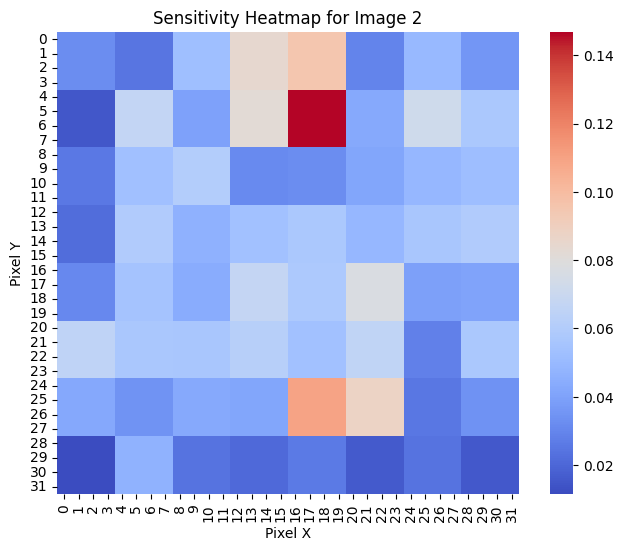}
        \includegraphics[width=0.092\linewidth]{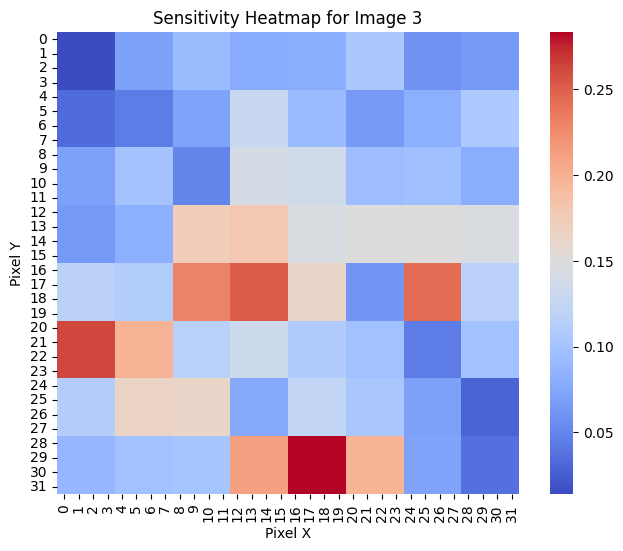}
        \includegraphics[width=0.092\linewidth]{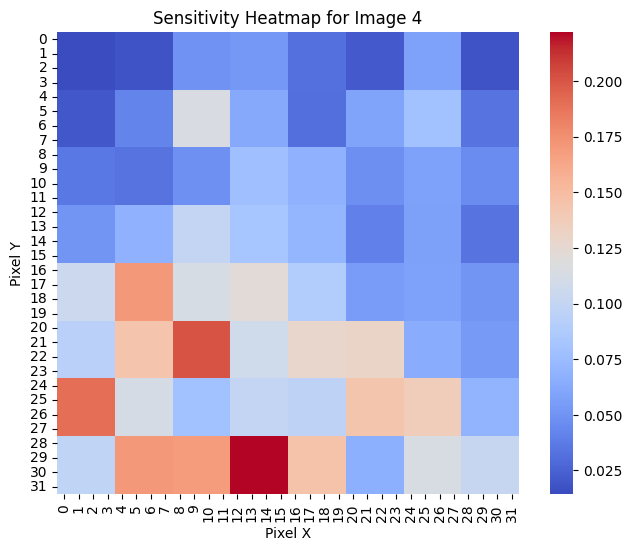}
        \includegraphics[width=0.092\linewidth]{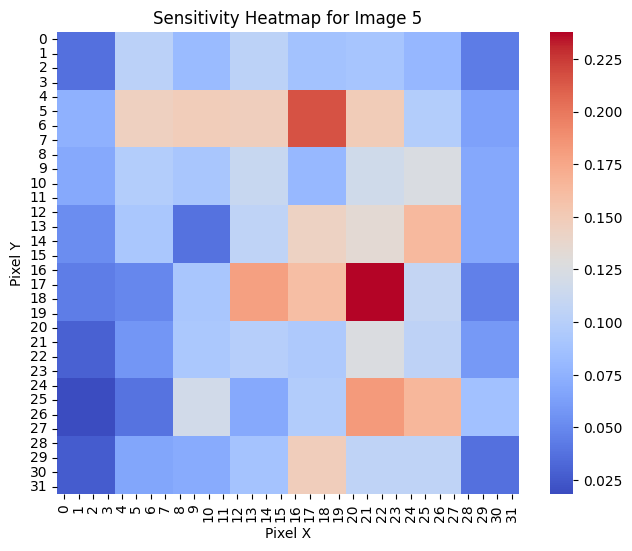}
        \includegraphics[width=0.092\linewidth]{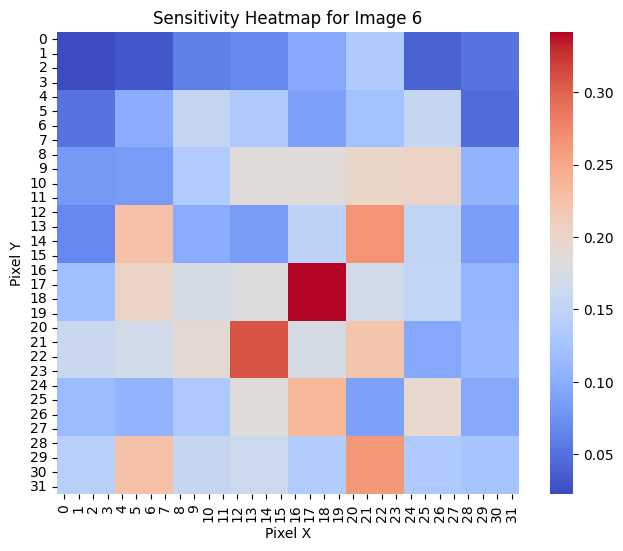}
        \includegraphics[width=0.092\linewidth]{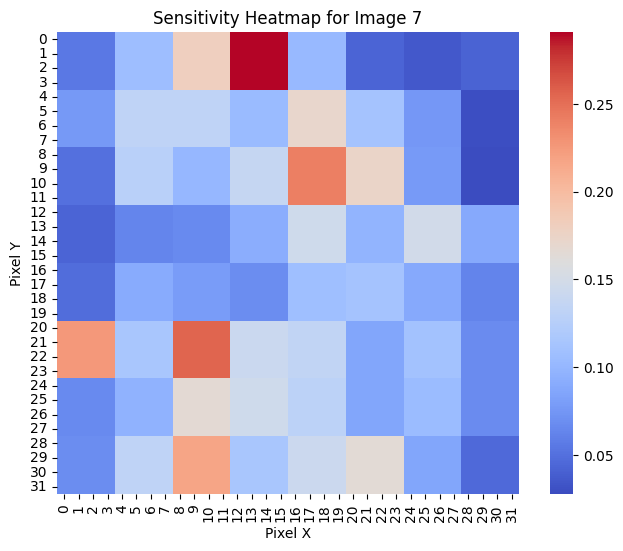}
        \includegraphics[width=0.092\linewidth]{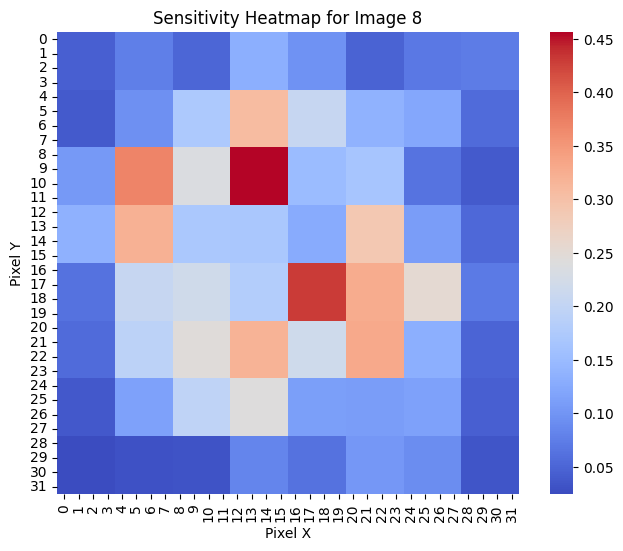}
%    \end{subfigure}
    \vspace{0.2cm} % 图片间距
    \caption{The sensitivity heatmap for the three color channels of block 4. The scale is $10^{-1}$.}
    \label{fig:SensitHeatmapBlock4_RBG_4x4}
\end{figure}

\begin{figure}[htbp]
    \centering   
%    \begin{subfigure}{\linewidth}
        \includegraphics[width=0.092\linewidth]{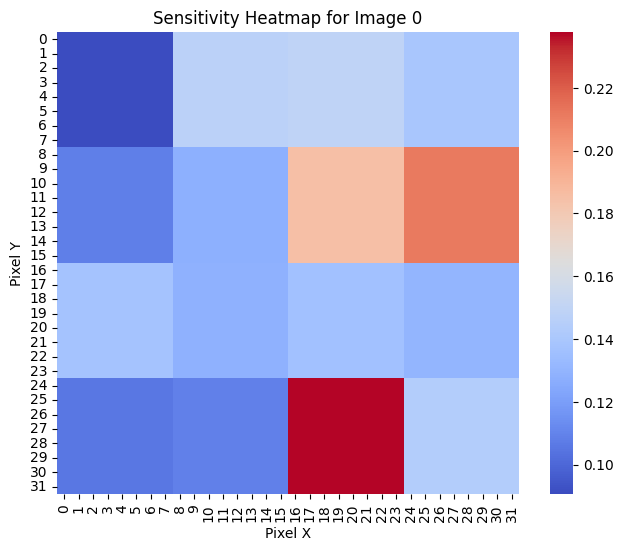}
        \includegraphics[width=0.092\linewidth]{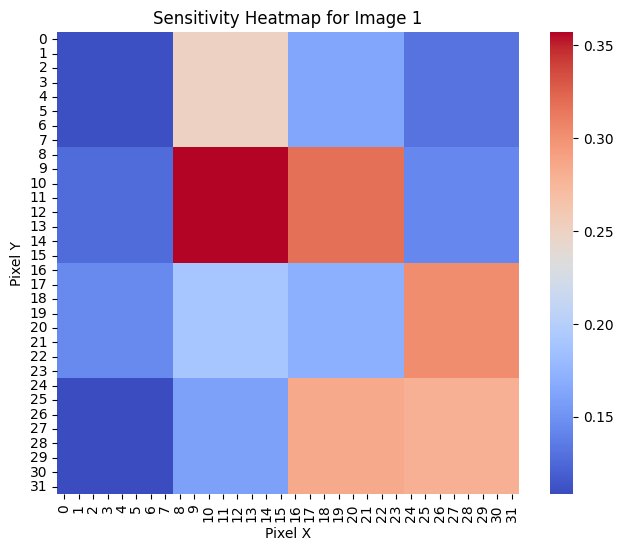}
        \includegraphics[width=0.092\linewidth]{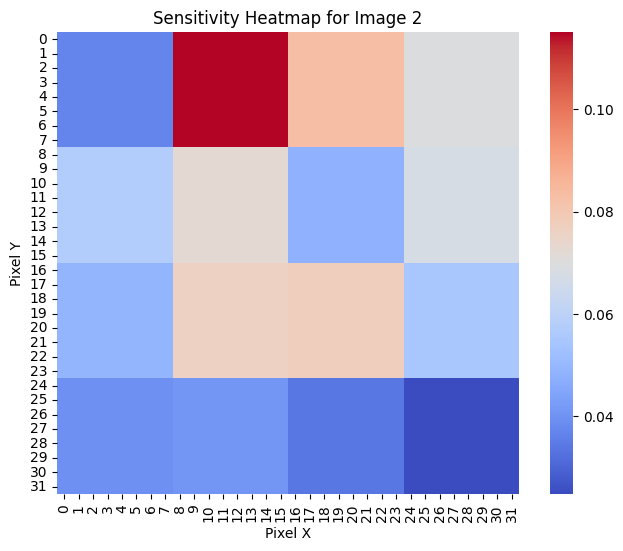}
        \includegraphics[width=0.092\linewidth]{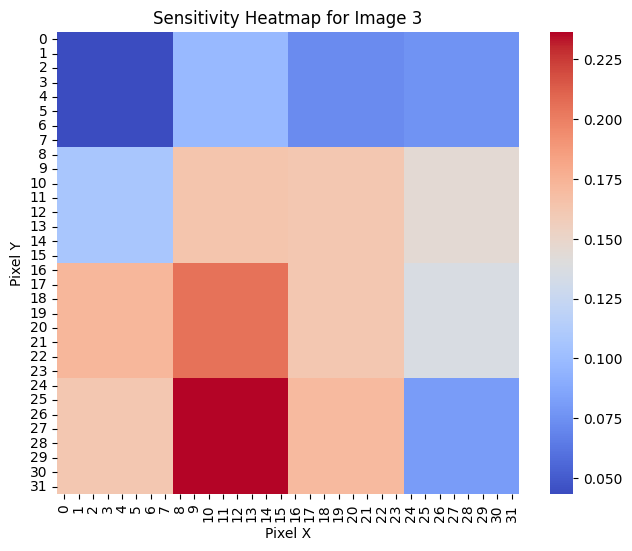}
        \includegraphics[width=0.092\linewidth]{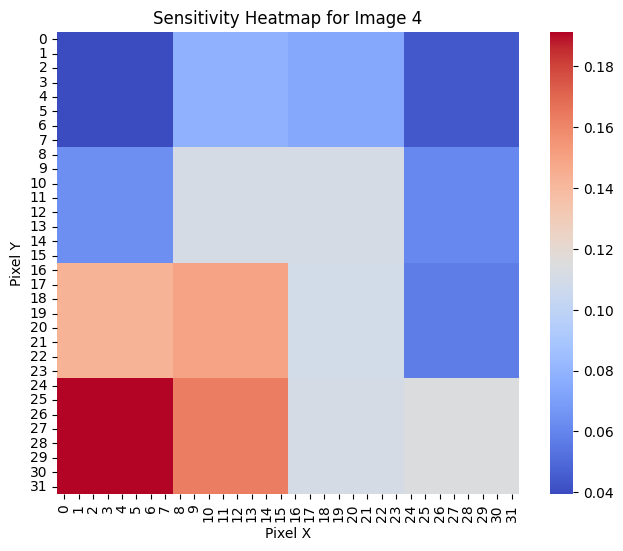}
        \includegraphics[width=0.092\linewidth]{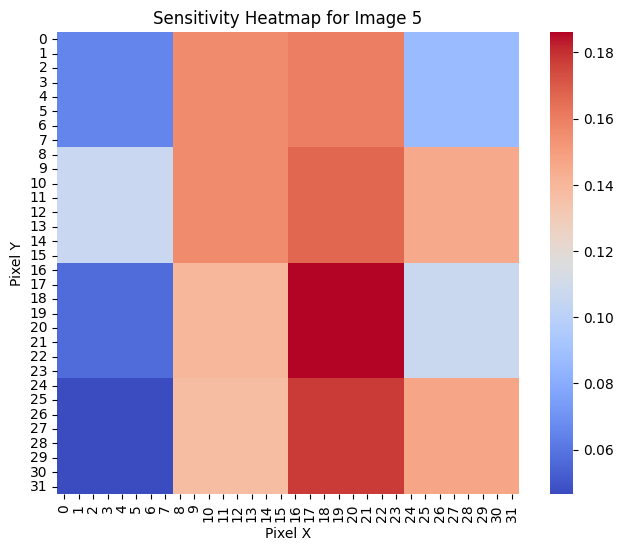}
        \includegraphics[width=0.092\linewidth]{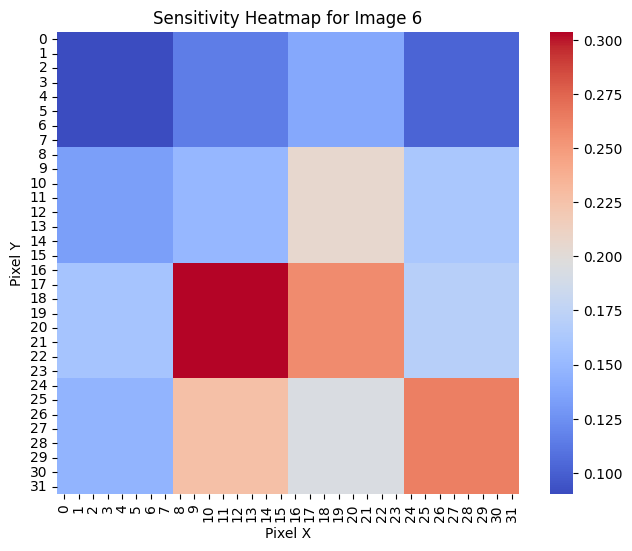}
        \includegraphics[width=0.092\linewidth]{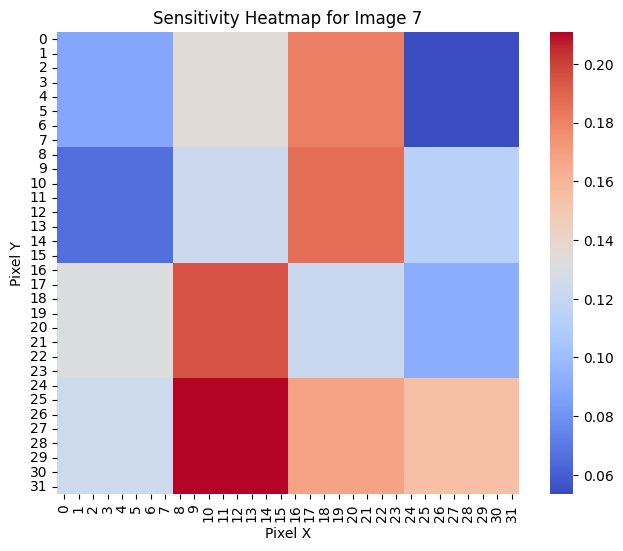}
        \includegraphics[width=0.092\linewidth]{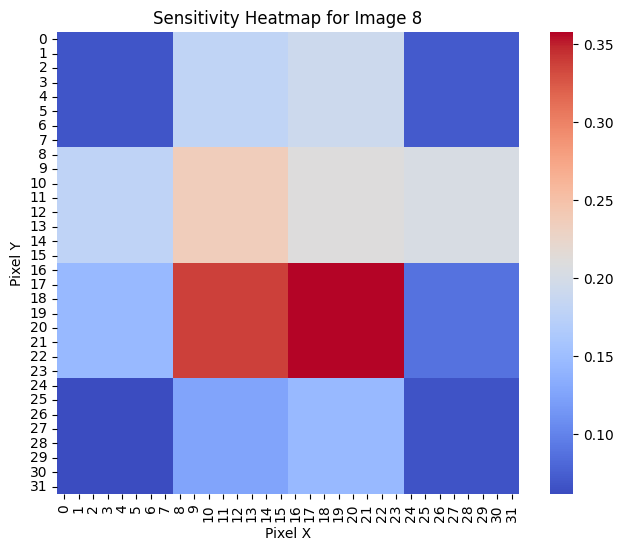}
%    \end{subfigure}
\\
%        \begin{subfigure}{\linewidth}
        \includegraphics[width=0.092\linewidth]{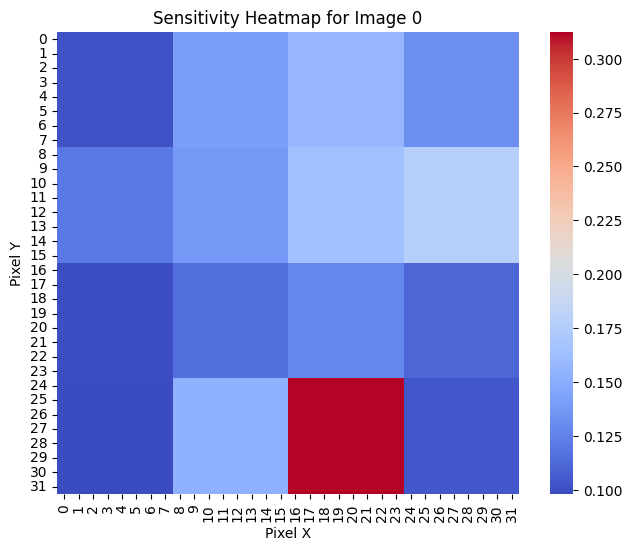}
        \includegraphics[width=0.092\linewidth]{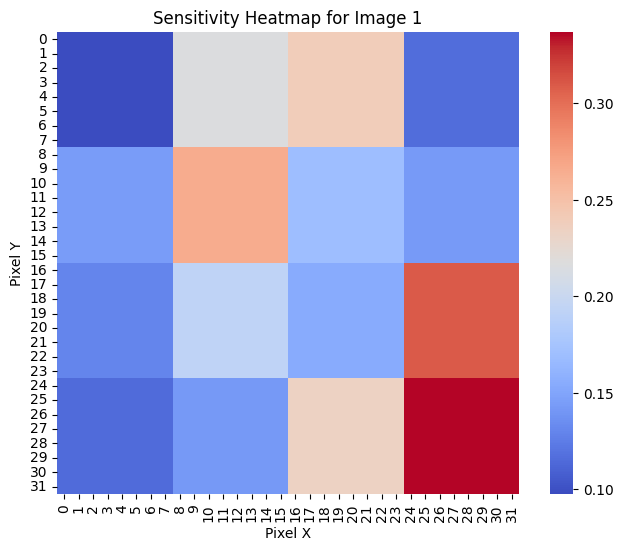}
        \includegraphics[width=0.092\linewidth]{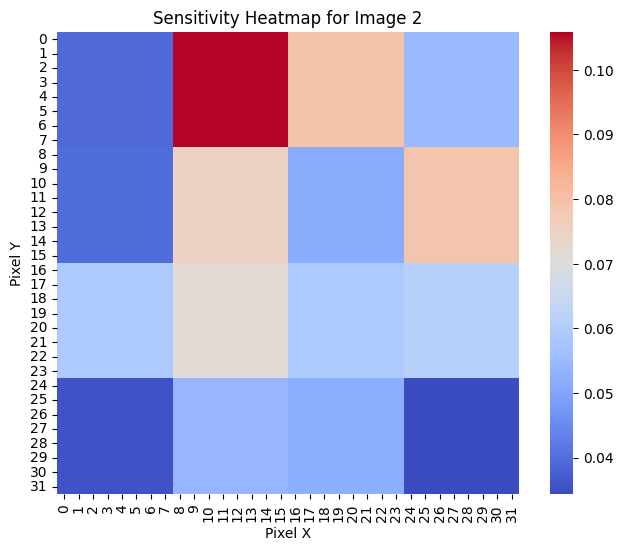}
        \includegraphics[width=0.092\linewidth]{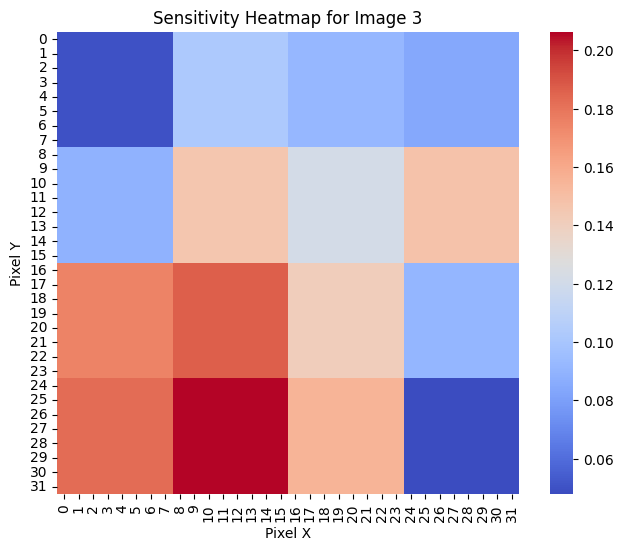}
        \includegraphics[width=0.092\linewidth]{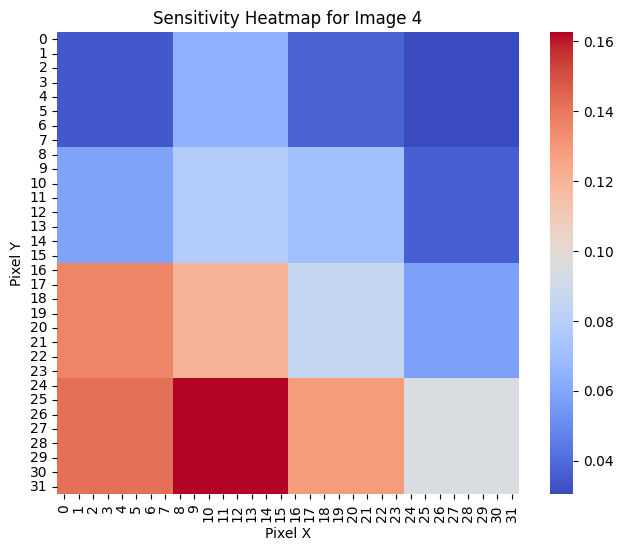}
        \includegraphics[width=0.092\linewidth]{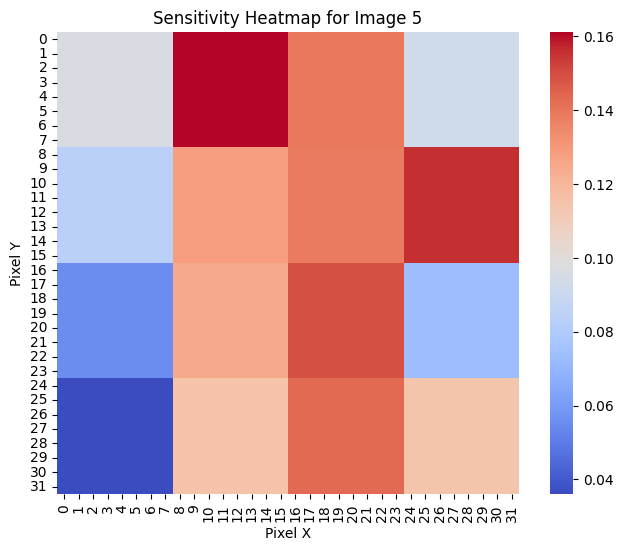}
        \includegraphics[width=0.092\linewidth]{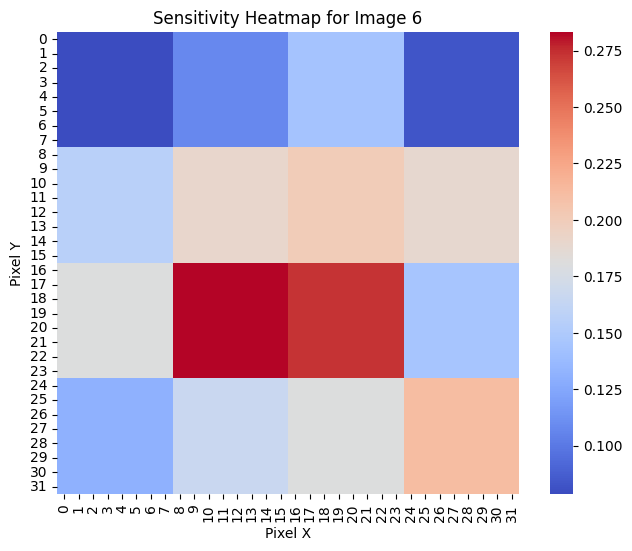}
        \includegraphics[width=0.092\linewidth]{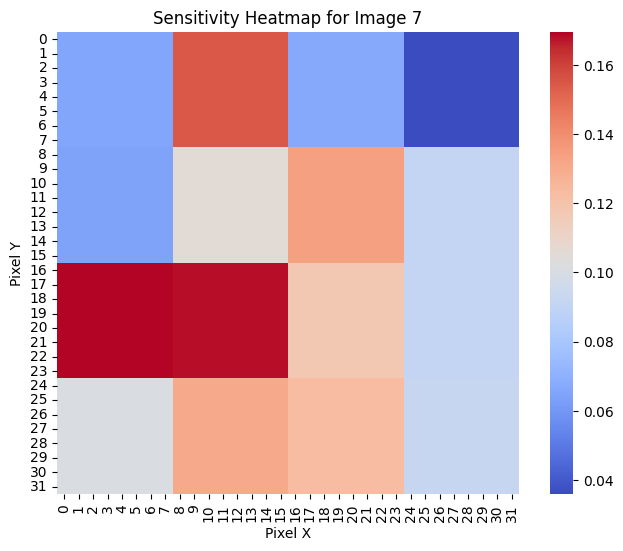}
        \includegraphics[width=0.092\linewidth]{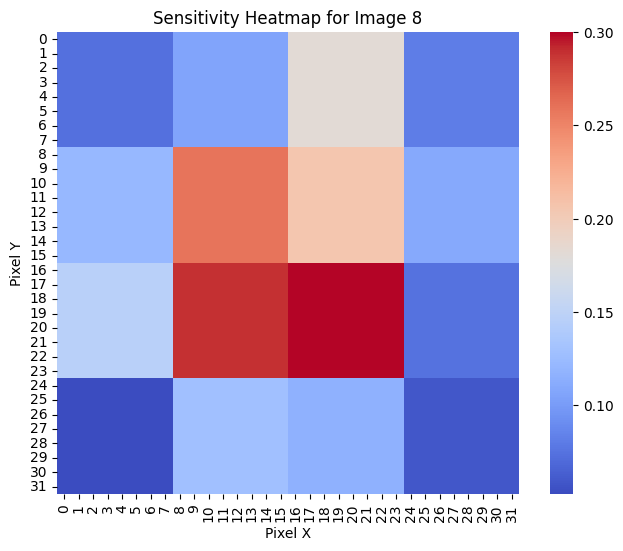}
%    \end{subfigure}
\\
%        \begin{subfigure}{\linewidth}
        \includegraphics[width=0.092\linewidth]{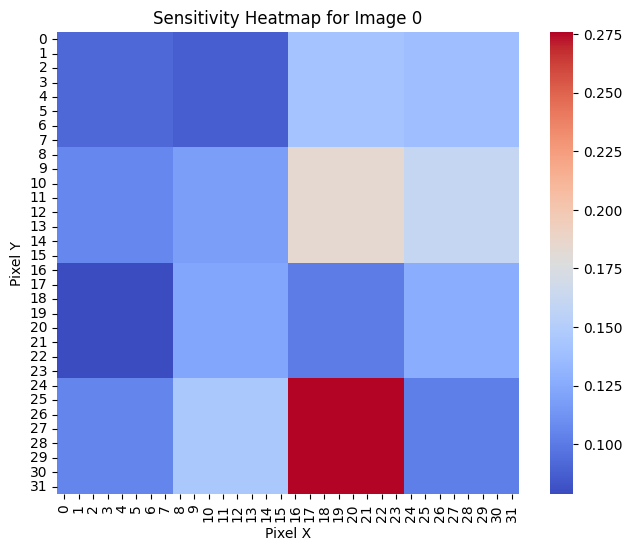}
        \includegraphics[width=0.092\linewidth]{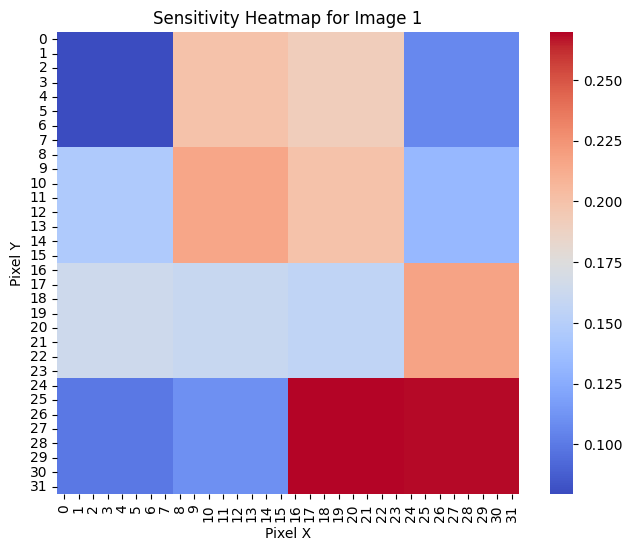}
        \includegraphics[width=0.092\linewidth]{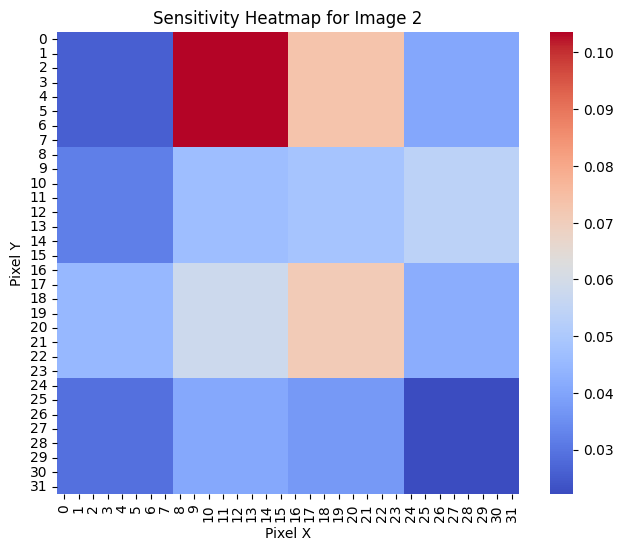}
        \includegraphics[width=0.092\linewidth]{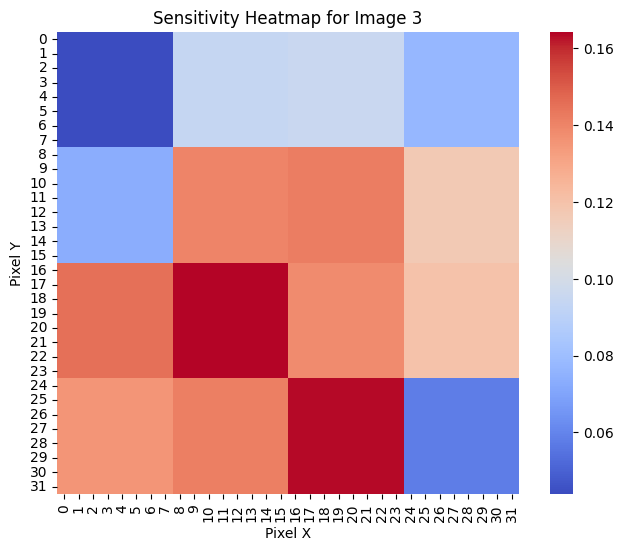}
        \includegraphics[width=0.092\linewidth]{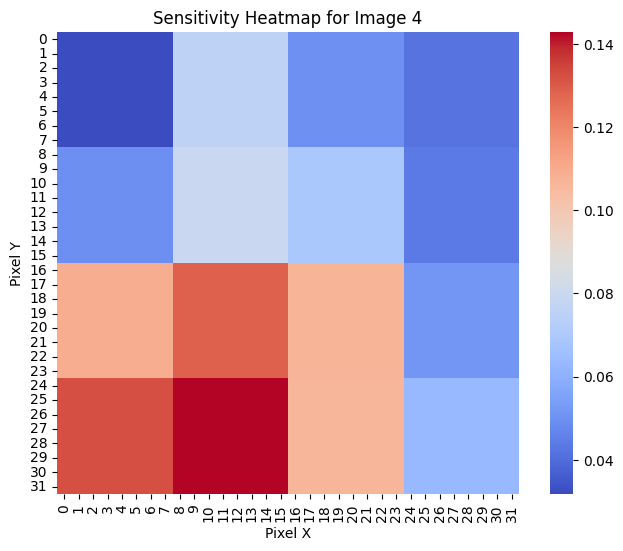}
        \includegraphics[width=0.092\linewidth]{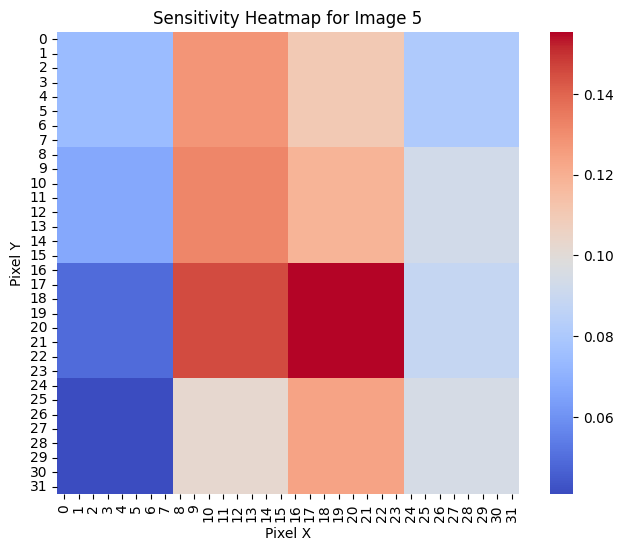}
        \includegraphics[width=0.092\linewidth]{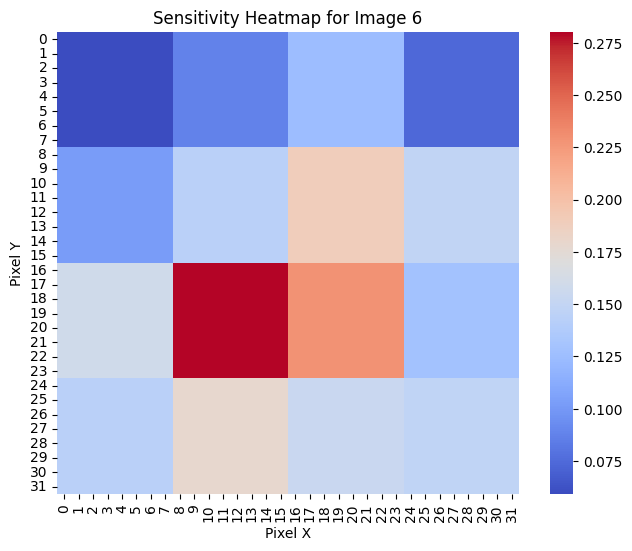}
        \includegraphics[width=0.092\linewidth]{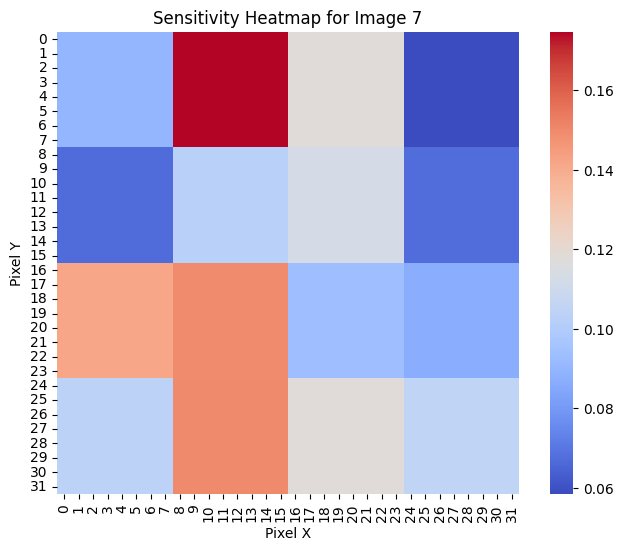}
        \includegraphics[width=0.092\linewidth]{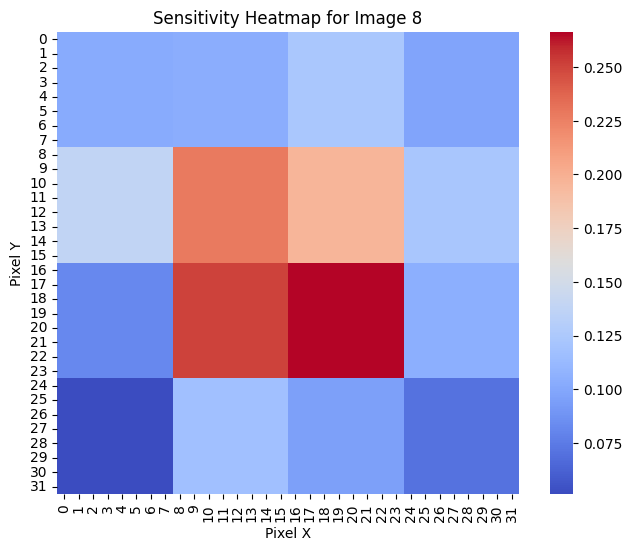}
%    \end{subfigure}
    \vspace{0.2cm} % 图片间距
    \caption{The sensitivity heatmap for the three color channels of block 4 ($8 \times 8$).}
    \label{fig:SensitHeatmapBlock4_RGB8x8}
\end{figure}

\begin{figure}[htbp]
    \centering   
%    \begin{subfigure}{\linewidth}
        \includegraphics[width=0.092\linewidth]{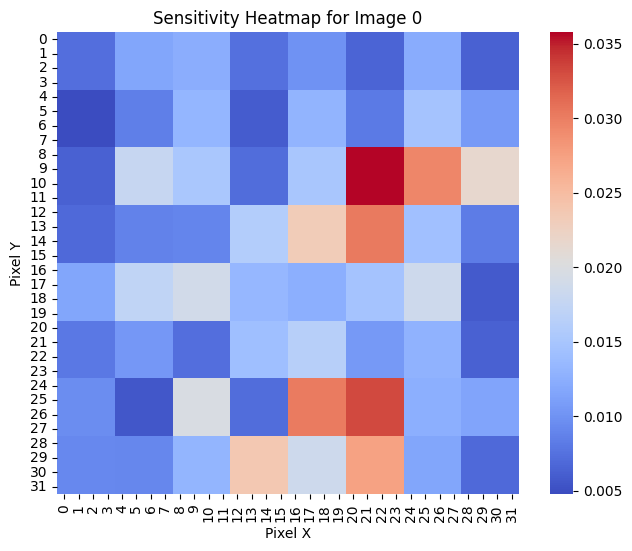}
        \includegraphics[width=0.092\linewidth]{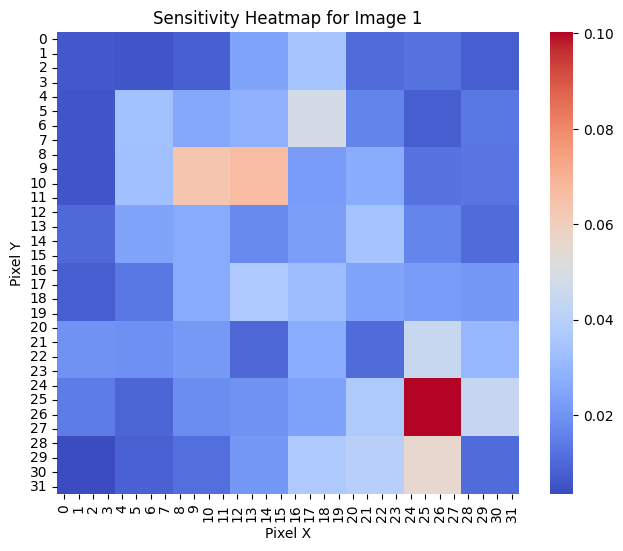}
        \includegraphics[width=0.092\linewidth]{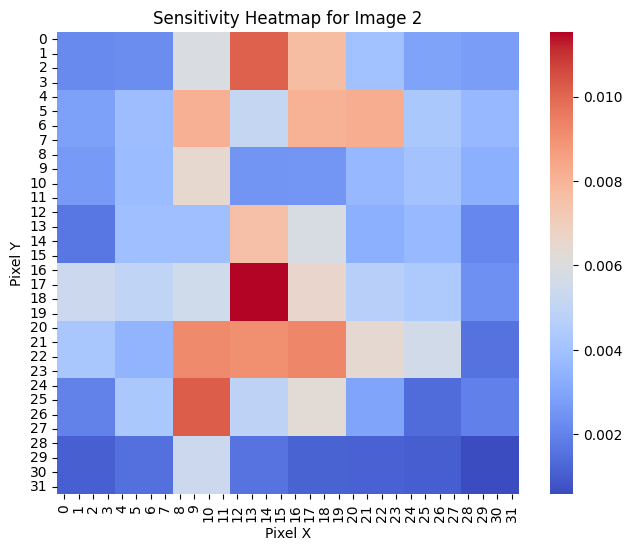}
        \includegraphics[width=0.092\linewidth]{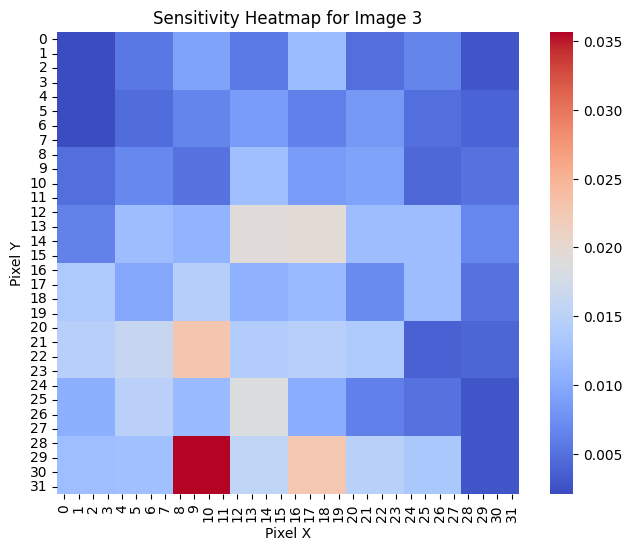}
        \includegraphics[width=0.092\linewidth]{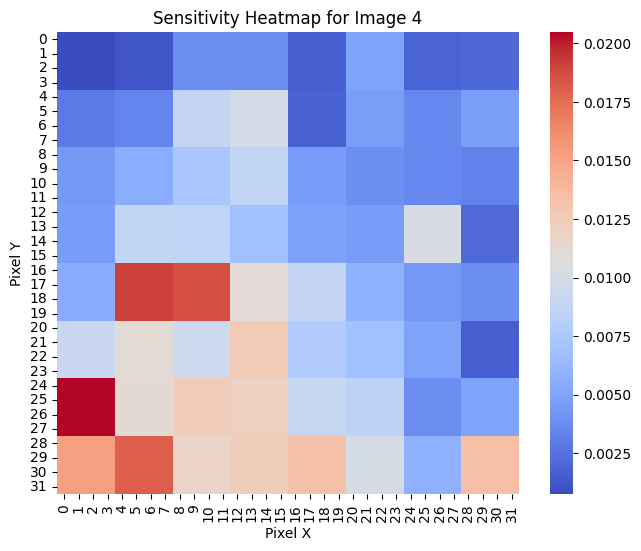}
        \includegraphics[width=0.092\linewidth]{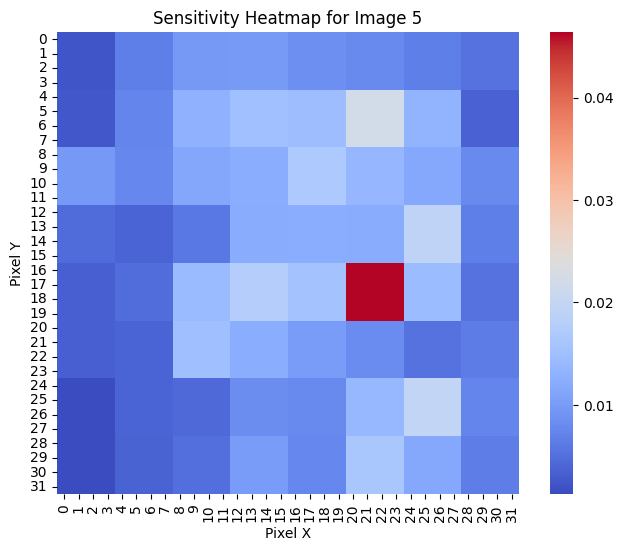}
        \includegraphics[width=0.092\linewidth]{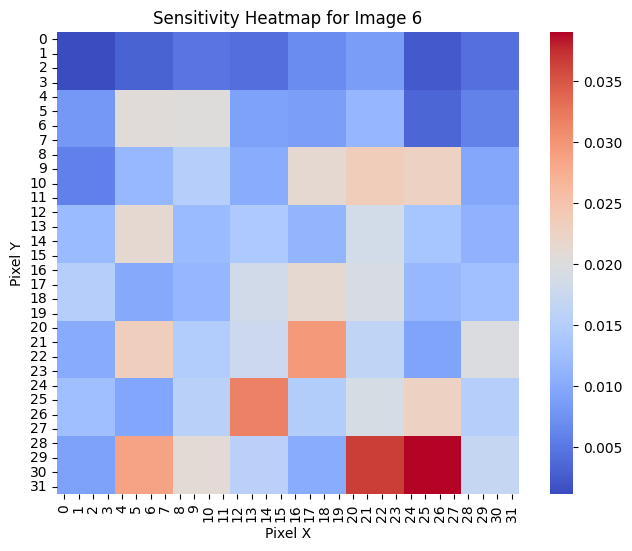}
        \includegraphics[width=0.092\linewidth]{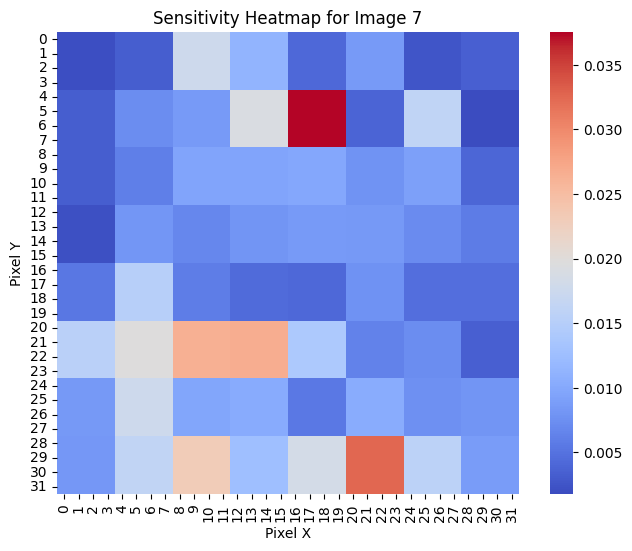}
        \includegraphics[width=0.092\linewidth]{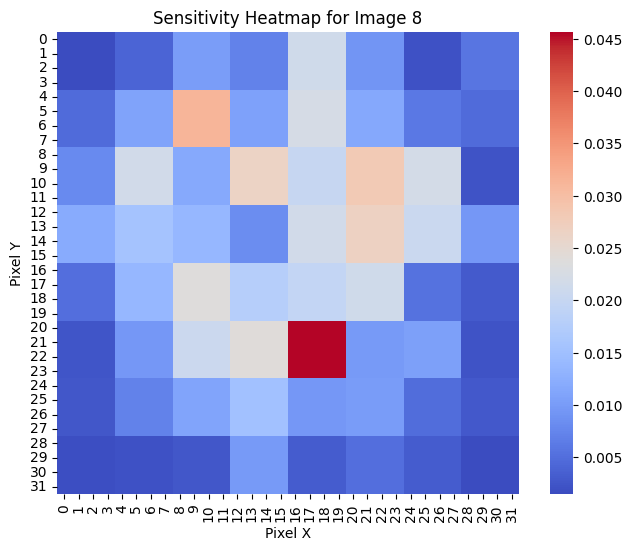}
%    \end{subfigure}
\\
%        \begin{subfigure}{\linewidth}
        \includegraphics[width=0.092\linewidth]{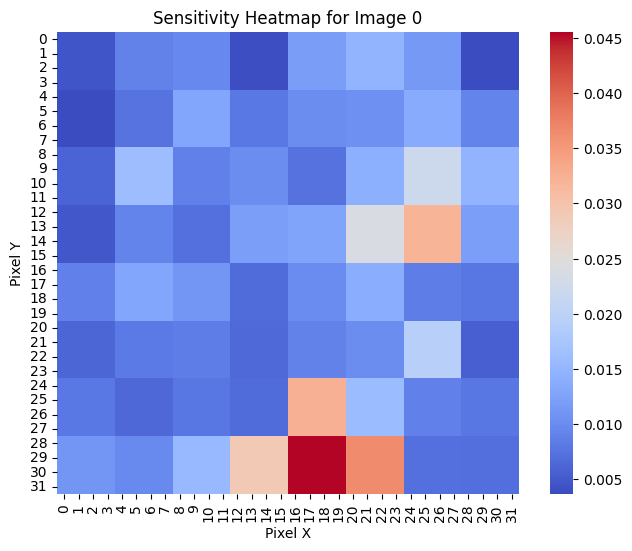}
        \includegraphics[width=0.092\linewidth]{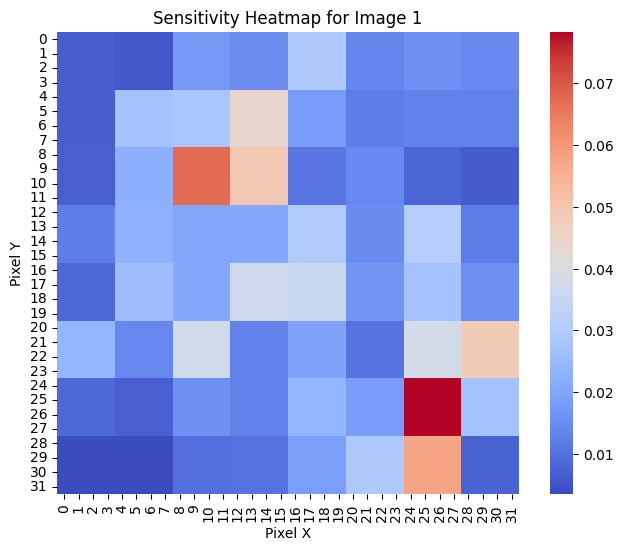}
        \includegraphics[width=0.092\linewidth]{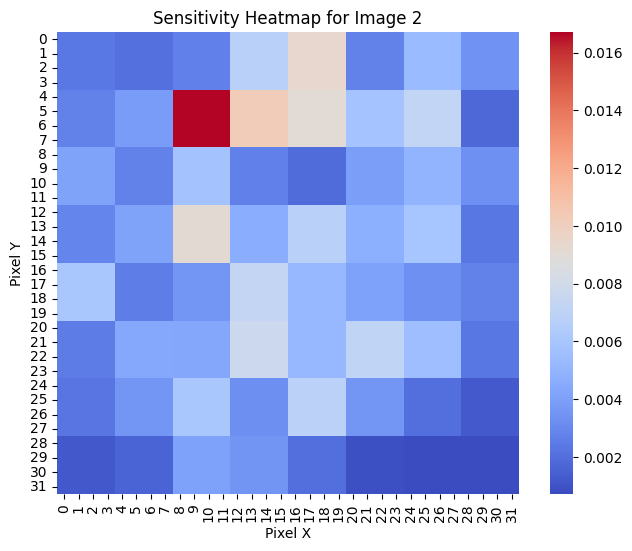}
        \includegraphics[width=0.092\linewidth]{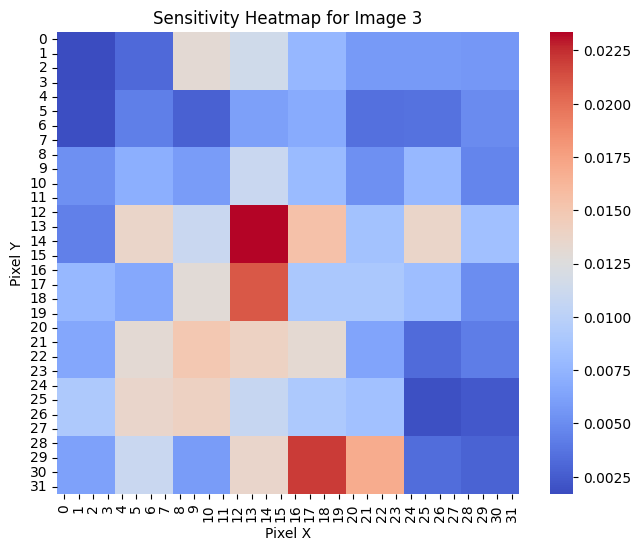}
        \includegraphics[width=0.092\linewidth]{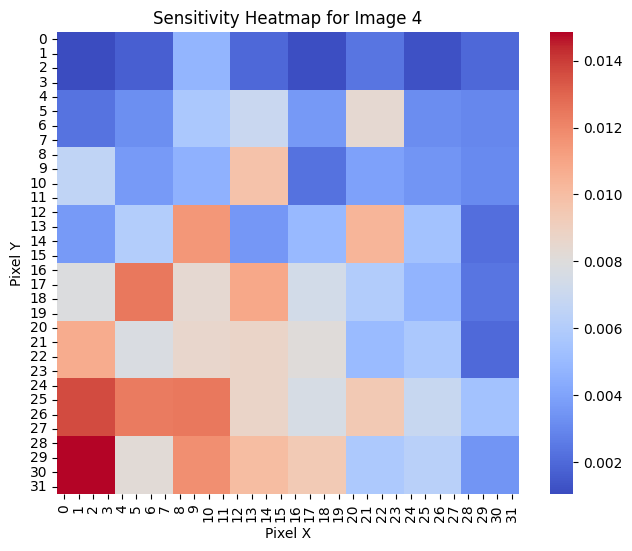}
        \includegraphics[width=0.092\linewidth]{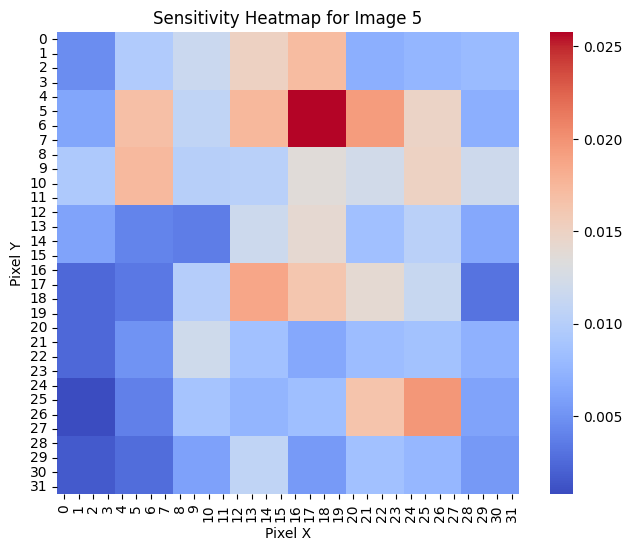}
        \includegraphics[width=0.092\linewidth]{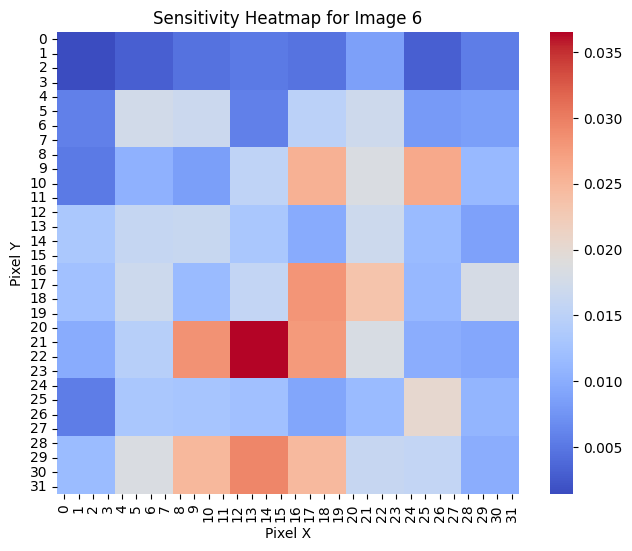}
        \includegraphics[width=0.092\linewidth]{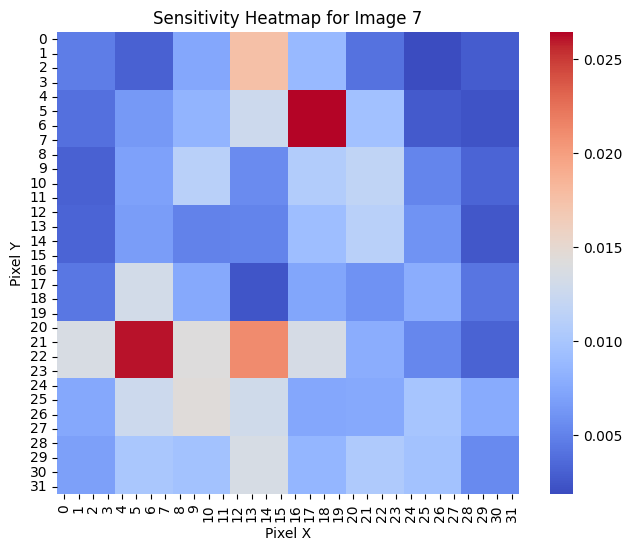}
        \includegraphics[width=0.092\linewidth]{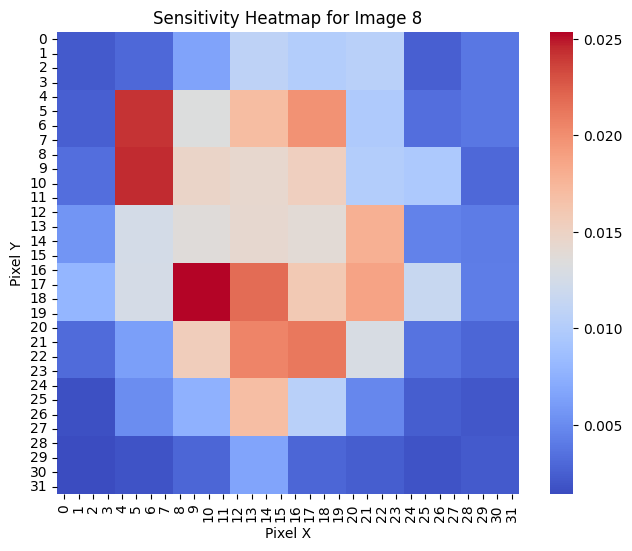}
%    \end{subfigure}
\\
%        \begin{subfigure}{\linewidth}
        \includegraphics[width=0.092\linewidth]{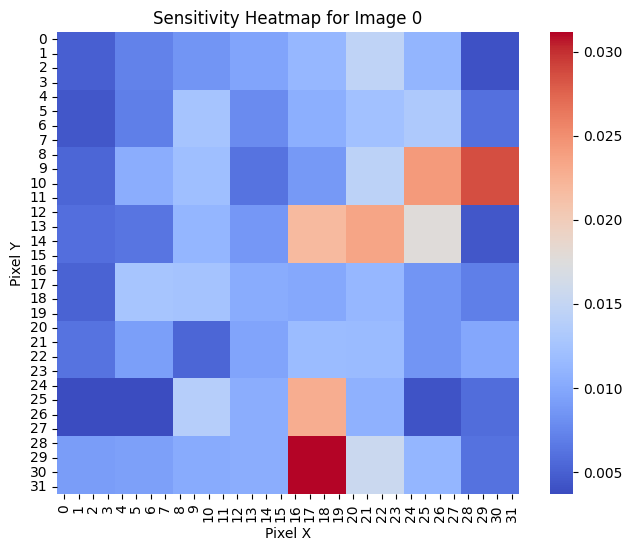}
        \includegraphics[width=0.092\linewidth]{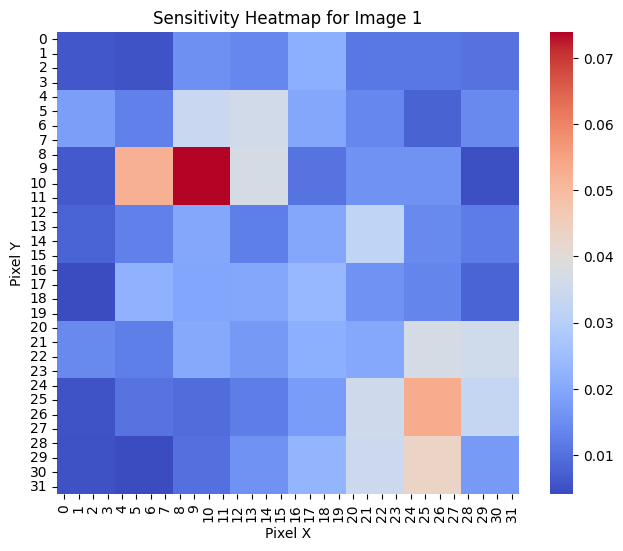}
        \includegraphics[width=0.092\linewidth]{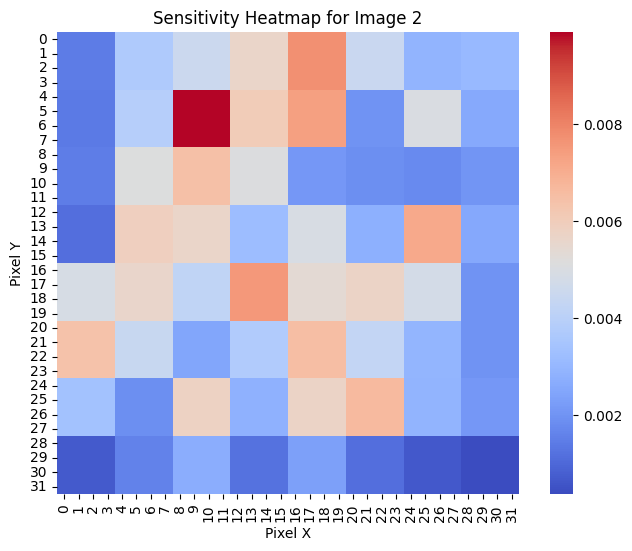}
        \includegraphics[width=0.092\linewidth]{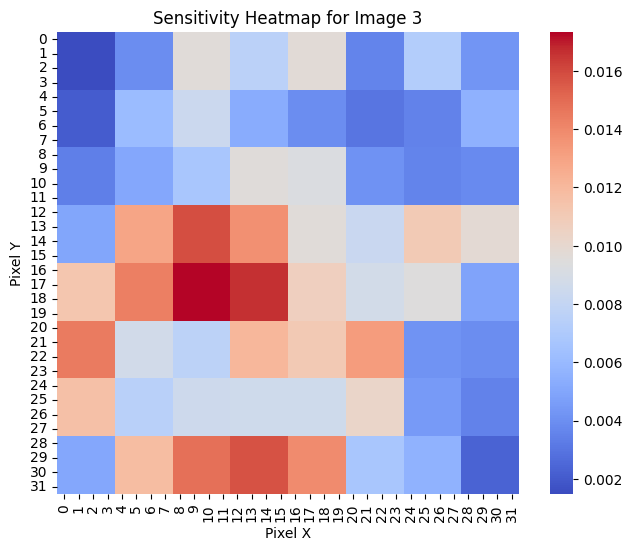}
        \includegraphics[width=0.092\linewidth]{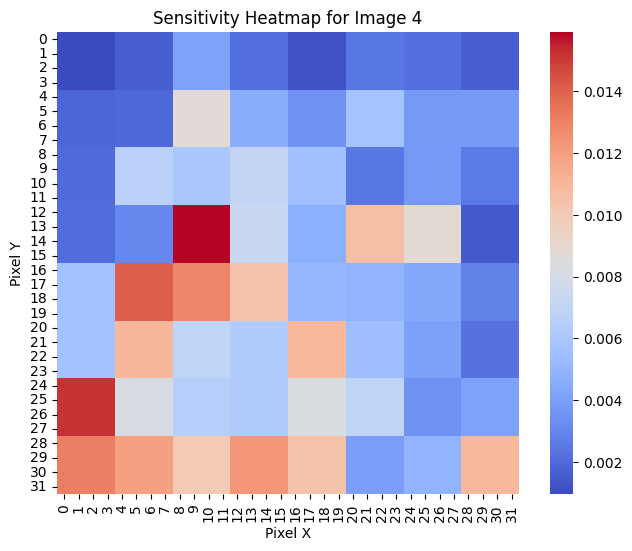}
        \includegraphics[width=0.092\linewidth]{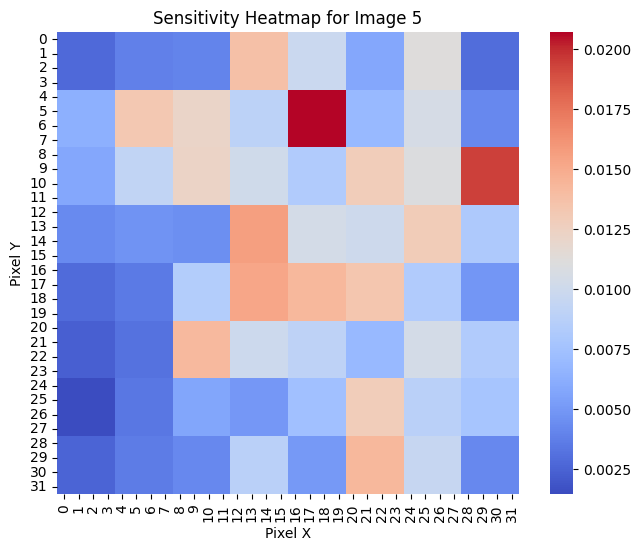}
        \includegraphics[width=0.092\linewidth]{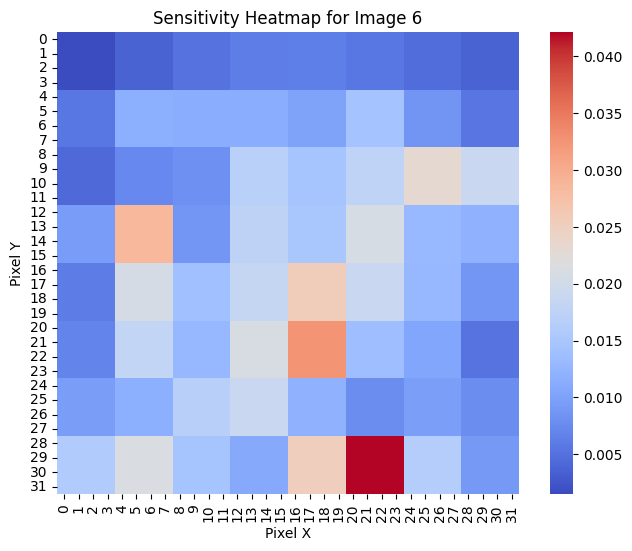}
        \includegraphics[width=0.092\linewidth]{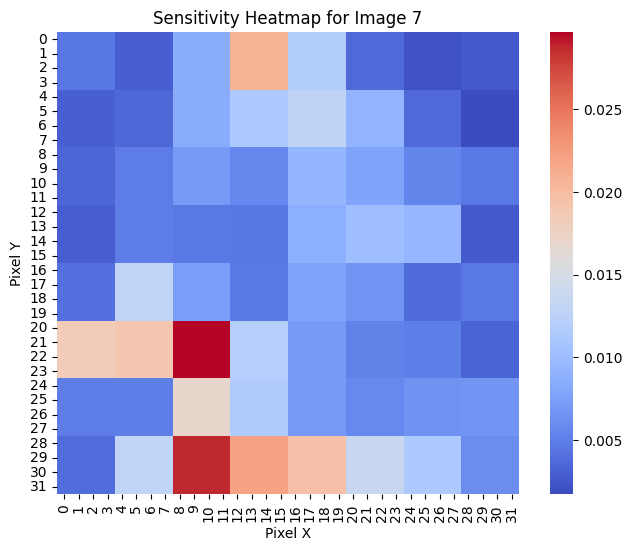}
        \includegraphics[width=0.092\linewidth]{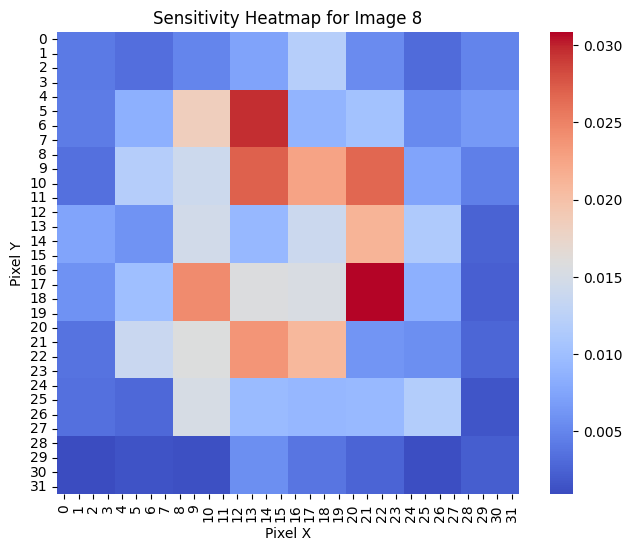}
%    \end{subfigure}
    \vspace{0.2cm} % 图片间距
    \caption{The sensitivity heatmap for the three color channels of block 5. The scale is $10^{-2}$.}
    \label{fig:SensitHeatmapBlock5_RGB}
\end{figure}

For convenience, the sensitivity images of the RGB channels for the fourth block were divided into larger pixel blocks. It is clear that the images for the three color channels are similar, and the influence of color on the fourth block is already quite minimal, with the processing now focused on shape and texture. The results are presented in Figures \ref{fig:SensitHeatmapBlock1} -- \ref{fig:SensitHeatmapBlock5_RGB}.

%\newpage
The chosen method for the local sensitivity analysis and the visualization of results through the heatmaps seems to be the simplest and most in tune with broader neural network visualization methods using attribution.

Attribution (heatmap, spectral relevance analysis) allows us to investigate which part of the input tensor of a neural network contributes to the activation of a specific neuron in the network. For example, in the object detection tasks, the attribution is often understood as a heatmap that shows the contribution of each point in the input image to object detection (for instance, see the review \cite{matveev2021overview} on neural network visualization methods). It is worth noting that the sensitivity analysis of the neural networks is much less discussed compared to their generalization ability \cite{recht2018cifar}.

Comparing different blocks and different color channels within a block we observe that the overall sensitivity gradually decreases with their depth. However, in each block, the significant sensitive areas can still be identified. In the first block, we see many sensitive areas at the edges and corners, but after these features are extracted in the neural network, their sensitivity significantly decreases. Thus, we can confidently assume that in the early stages, the convolutional layers of the VGG-16 model gradually capture low-level features such as color, edges, and corners. Subsequently, the model starts capturing more complex internal shapes and patterns. It can be concluded that using deeper convolutional layers with more parameters in the architecture of a convolutional neural network may be an ineffective solution due to the gradual reduction in the size of sensitive regions.

Despite we obtain the interesting results via visualization of the sensitivity maps, we cannot make many conclusive inferences from them. Therefore, the next section is dedicated to another method of analyzing neural network sensitivity -- the activation maximization method.

\subsection{Activation maximization method}

Activation maximization (AM) visualizes the input preferences for each neuron and helps to determine which input data can maximally activate a specific neuron at a certain level \cite{pospelov2025fast}. Heuristically, one can assume that each neuron is responsible for extracting a certain feature. If the feature extracted by a particular neuron is known, then input data with this feature should strongly activate that neuron. On the other hand, if it is unknown which specific features the neuron extracts, the input data can be iteratively adjusted so that the activation of this neuron becomes maximal, thereby visualizing the features that the neuron has learned to recognize.

Let \( \theta \) be the parameters of the classifier, which maps an image 
\( x \in \mathbb{R}^{H \times W \times C} \) (where \( C \) is the number of color channels, each of which has width \( W \) pixels and height 
\( H \) pixels) to a probability distribution over the output classes. The search for an image \( x \) that maximizes the activation \( a_{ij}(\theta, x) \) of the neuron with index \( i \) in the given layer \( j \) of the classifier network can be formulated as an optimization problem \cite{simonyan2013deep}:
\begin{equation}
x^* = \arg \max_x a_{ij}(\theta, x)
\label{eq:AM}
\end{equation}
Finding the optimal solution using the gradient ascent method:
\begin{equation}
x_{'} = x + \eta \frac{\partial a_{ij}(\theta, x)}{\partial x}
\label{eq:GradStep}
\end{equation}

One of the problems of the aforementioned method is that the resulting visualized inputs often have low distinguishability and differ significantly from real images (since iterations are performed on a randomly initialized input $x_0$, and there are no measures to ensure the found input corresponds to a real image; extreme pixels, structured high-frequency patterns, etc. can lead to high activation values). To find a solution that is as close as possible to a real image, appropriate regularization terms need to be introduced.

For example, prior knowledge about images can help smooth the image or penalize pixels with extreme intensity. Such constraints are typically introduced into the formulation of the Activation Maximization (AM) method as a regularization term $R(x)$:
\begin{equation}
x^* = \arg\max_x \left( a(x) - R(x) \right)
\end{equation}

For example, to promote image smoothing in AM, the regularization term 
$R: \mathbb{R}^{H \times W \times C} \to \mathbb{R}$ can compute the total variation (TV) of the image. Thus, at each update, we follow the gradient to maximize the neuron activation using \eqref{eq:AM}, and \eqref{eq:GradStep} to minimize the total variation loss:
\begin{equation}
x_{t+1} = x_t + \varepsilon_1 \frac{\partial a(x_t)}{\partial x_t} - \varepsilon_2 \frac{\partial R(x_t)}{\partial x_t}.
\end{equation}

However, in practice, we do not always compute the analytical gradient 
$\dfrac{\partial R(x_t)}{\partial x_t}$. Instead, we can define a regularization operator
$$r: \mathbb{R}^{H \times W \times C} \to \mathbb{R}^{H \times W \times C}.$$
For example, one may use a Gaussian blur kernel, and at each step we map $x$ to a more regularized version of itself. In this case, the update step becomes:
\begin{equation}
x_{t+1} = r(x_t) + \varepsilon_1 \frac{\partial a(x_t)}{\partial x_t}.
\end{equation}

\subsection{Application of AM to VGG-16}

\begin{figure}[htbp]
    \centering
    %\begin{subfigure}{1\linewidth}
        \includegraphics[width=0.11\linewidth]{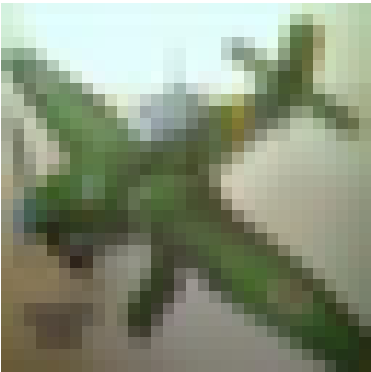}
        \includegraphics[width=0.11\linewidth]{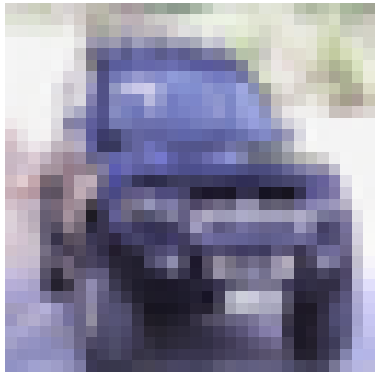}
        \includegraphics[width=0.11\linewidth]{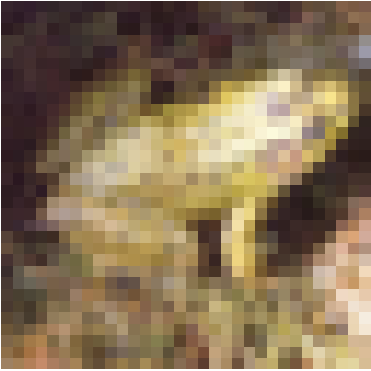}
        \includegraphics[width=0.11\linewidth]{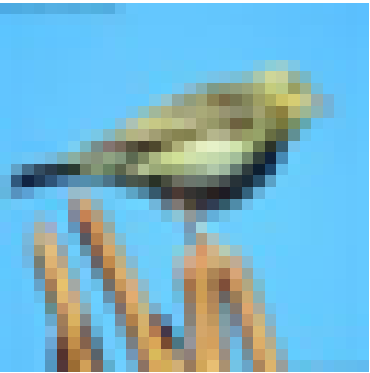}
        \includegraphics[width=0.11\linewidth]{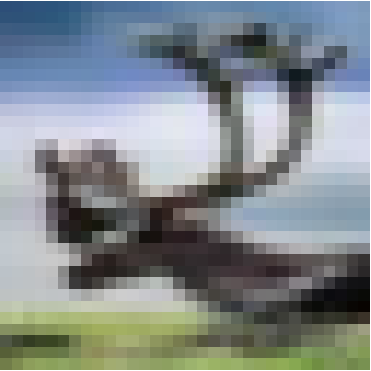}
        \includegraphics[width=0.11\linewidth]{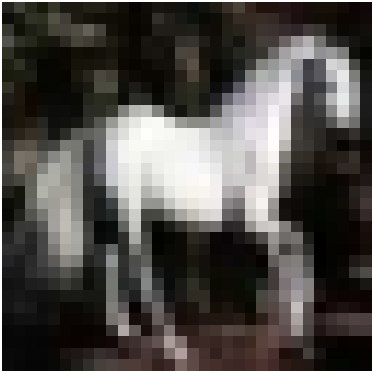}
        \includegraphics[width=0.11\linewidth]{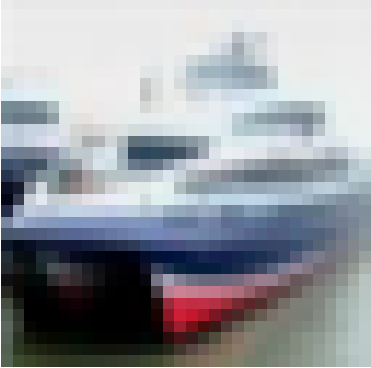}
        \includegraphics[width=0.11\linewidth]{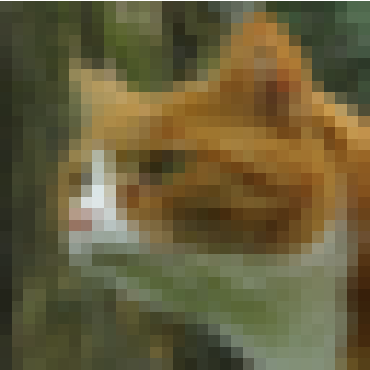}
    \\
    %\begin{subfigure}{1\linewidth}
        \includegraphics[width=0.11\linewidth]{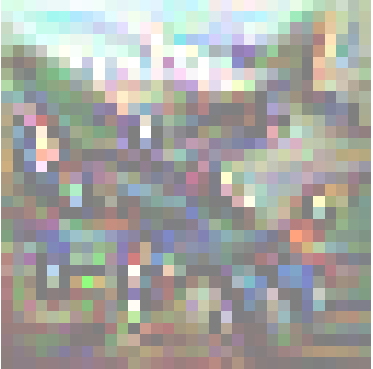}
        \includegraphics[width=0.11\linewidth]{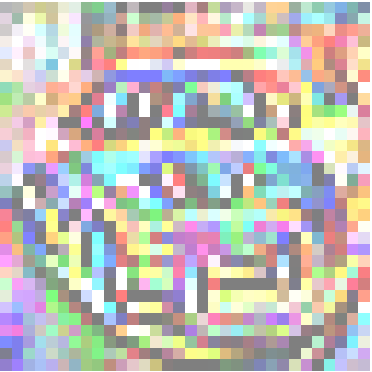}
        \includegraphics[width=0.11\linewidth]{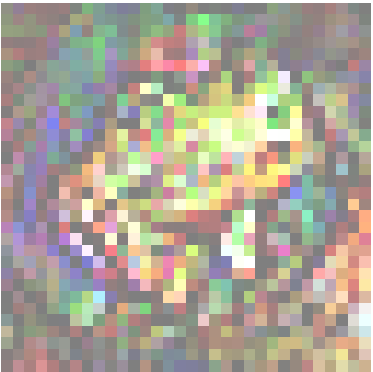}
        \includegraphics[width=0.11\linewidth]{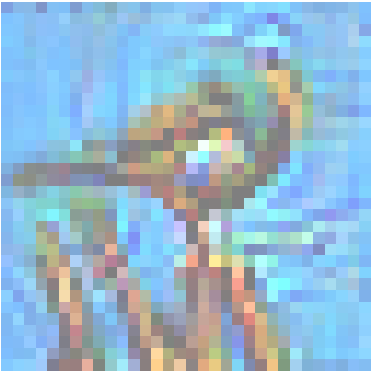}
        \includegraphics[width=0.11\linewidth]{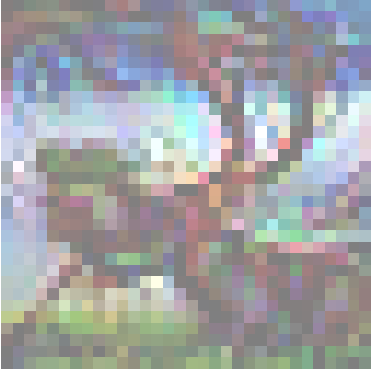}
        \includegraphics[width=0.11\linewidth]{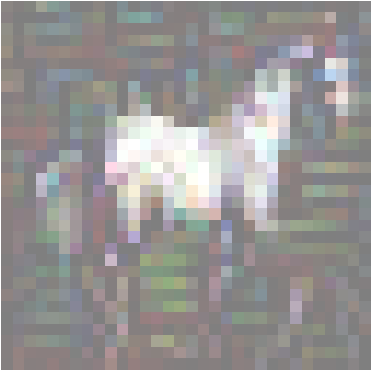}
        \includegraphics[width=0.11\linewidth]{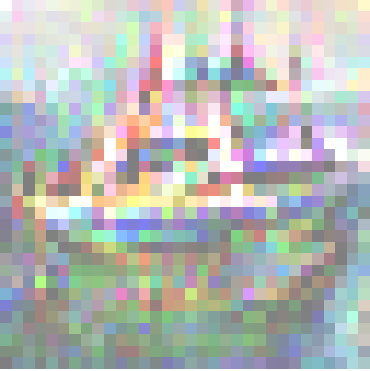}
        \includegraphics[width=0.11\linewidth]{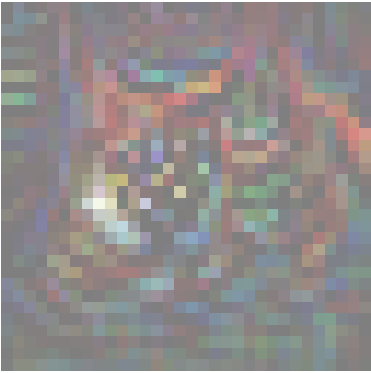}
    %\end{subfigure}
    \caption{The activation maximization pictures with
total variation (TV) for VGG-16.}
    \label{fig:AM}
\end{figure}

We utilize the total variation (TV) regularization method for activation of certain categories from the CIFAR-10 dataset and obtained the following activation maximization results ($\varepsilon_1=0.1$, $\varepsilon_2=0.1$), as shown in Figure \ref{fig:AM}.

It is easy to notice that different classifiers in the VGG-16 model are good at extracting contours and features corresponding to specific categories. For instance, the wings of an airplane, the antlers of a deer, and the ears of a cat — all these features were well extracted.

\begin{figure}[H]
    \centering
%    \begin{subfigure}{1\linewidth}
        \includegraphics[width=0.11\linewidth]{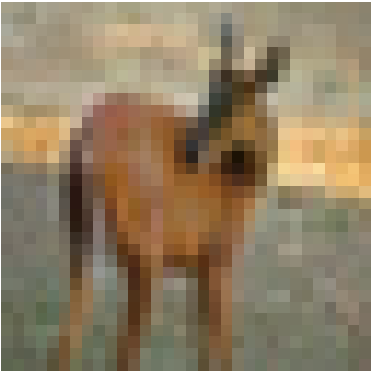}
        \includegraphics[width=0.11\linewidth]{images/AM2/DEER.png}
        \includegraphics[width=0.11\linewidth]{images/AM2/DEER.png}
        \includegraphics[width=0.11\linewidth]{images/AM2/DEER.png}
        \includegraphics[width=0.11\linewidth]{images/AM2/DEER.png}
        \includegraphics[width=0.11\linewidth]{images/AM2/DEER.png}
        \includegraphics[width=0.11\linewidth]{images/AM2/DEER.png}
        \includegraphics[width=0.11\linewidth]{images/AM2/DEER.png}
%    \end{subfigure}
\\
    \vspace{0.2cm} 
%    \begin{subfigure}{1\linewidth}
        \includegraphics[width=0.11\linewidth]{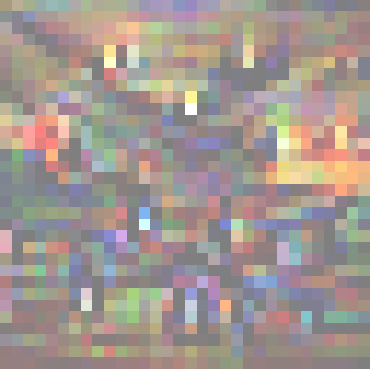}
        \includegraphics[width=0.11\linewidth]{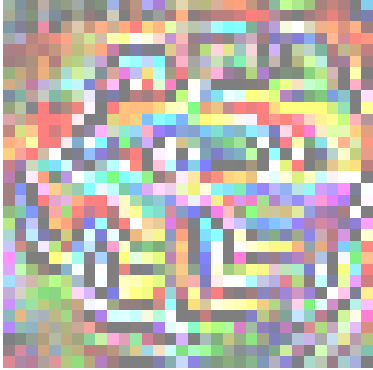}
        \includegraphics[width=0.11\linewidth]{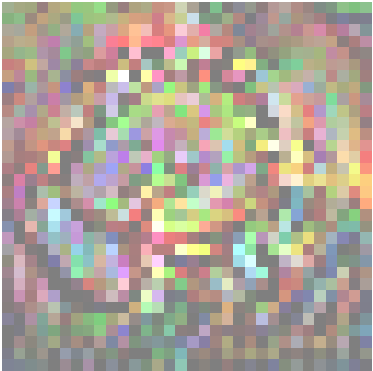}
        \includegraphics[width=0.11\linewidth]{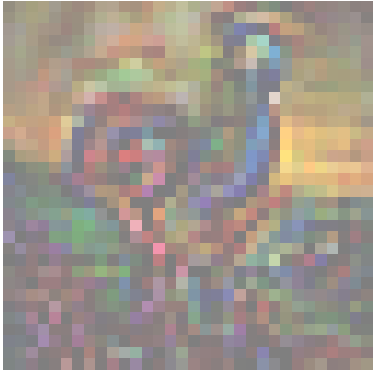}
        \includegraphics[width=0.11\linewidth]{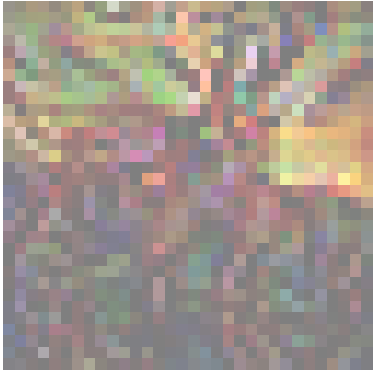}
        \includegraphics[width=0.11\linewidth]{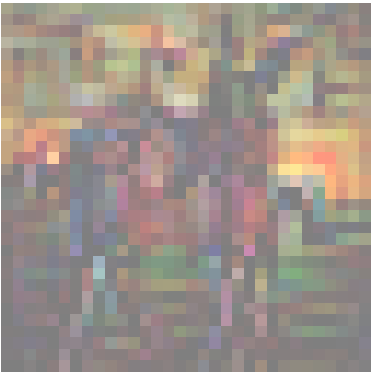}
        \includegraphics[width=0.11\linewidth]{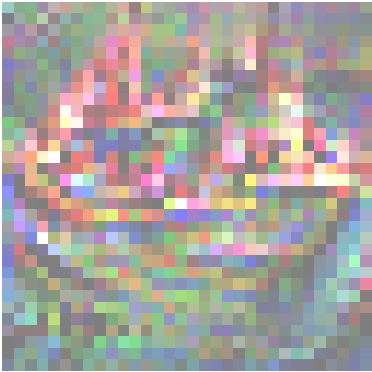}
        \includegraphics[width=0.11\linewidth]{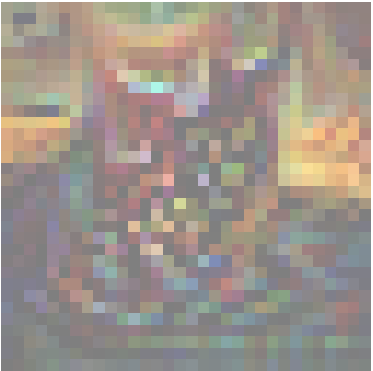}
%    \end{subfigure}
    \caption{The activation images for the features of an airplane, car, frog, bird, deer, horse, ship, and cat. }
    \label{fig:AM_wrong}
\end{figure}

Above, we successfully used parameters corresponding to a specific class of images to activate that class and obtained good results. However, such images may also share characteristics with other classes, so studying the effect of activating these images with parameters related to other features is of great interest. For example, the image of a deer without large antlers was activated using the parameters of an airplane, car, frog, bird, deer, horse, ship, and cat, resulting in new images \ref{fig:AM_wrong}.

We can notice that based on the original image, features corresponding to the activated categories have appeared. Special attention should be given to the fifth and sixth images. These images show that an important characteristic for activating the deer image is the presence of antlers, while the other features are generally similar to those of a horse. Thus, in the absence of antlers, the model may not accurately distinguish between a deer and a horse.

\subsection{Sensitivity analysis for ResNet-18}

Convolutional neural networks are capable of extracting low, medium, and high-level features of an object detection task. The deeper the network, the more diverse features from different levels can be extracted, meaning deeper networks capture more abstract and semantically significant features. However, due to the ``degradation problem'' of neural networks, as the network depth increases, more information is lost.

We investigate the local sensitivity of the VGG-16 model above, and see that as data passes through the model the sensitivity gradually decreases with increasing depth, reaching a scale of $10^{-2}$ in the fifth block. This is an issue for the larger VGG models. In 2015, a special network has been introduced \cite{he2016deep}. Iy employs the residual blocks (or residual units) based on mappings that mix the input signal with the result of the activation function:
$$
\text{Res}(x) = x+ F_{act}(x) 
$$
This solution addressed the inefficiency of very deep convolutional neural networks. Therefore, we consider performing the local sensitivity analysis after training ResNet-18 on the CIFAR-10 dataset. The normalization, standardization parameters, and learning rate are the same as for the previous VGG-16 model, and the model achieved an accuracy of 98.38\% on the training set and 83.04\% on the test set. The software implementation of the corresponding experiments is still available \footnote{See the following link: \url{https://github.com/LUCI1a/CIFAR-10_projects.git}}.

By analogy with VGG16 and using the first 9 images from the airplane category, we obtain local sensitivity heatmaps for layers 1, 2, 3, and 4.

\begin{figure}[htbp]
    \centering  
%    \begin{subfigure}{\textwidth}
        \includegraphics[width=0.092\linewidth]{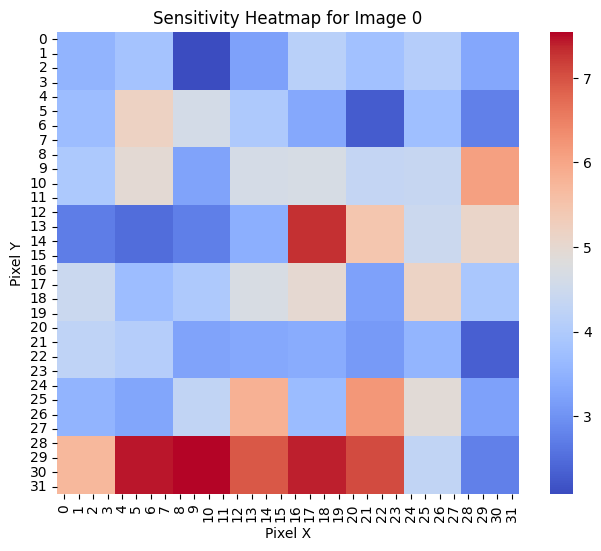}
        \includegraphics[width=0.092\linewidth]{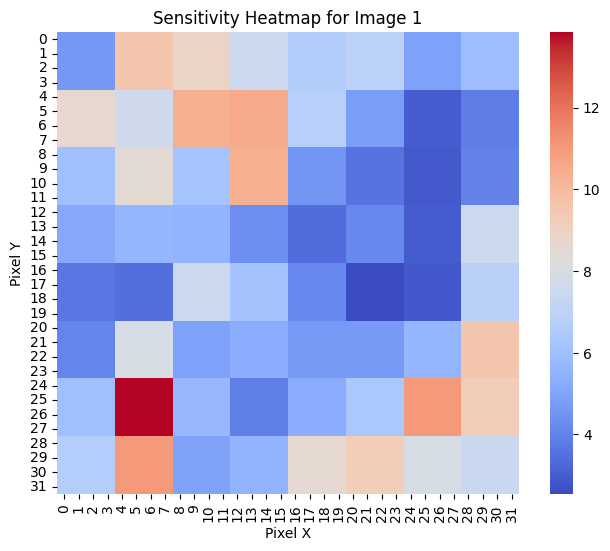}
        \includegraphics[width=0.092\linewidth]{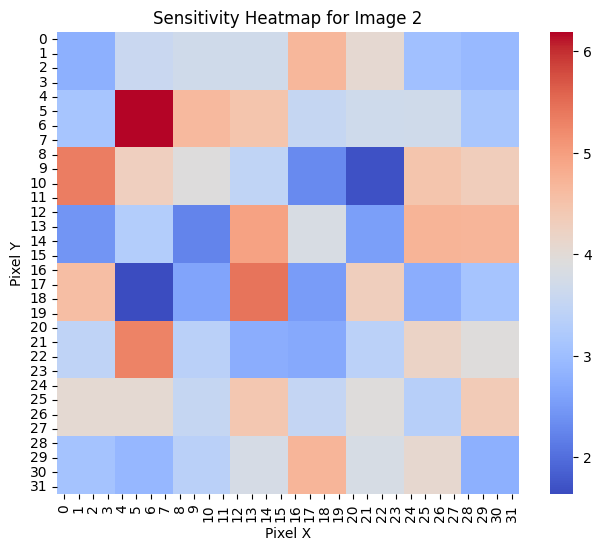}
        \includegraphics[width=0.092\linewidth]{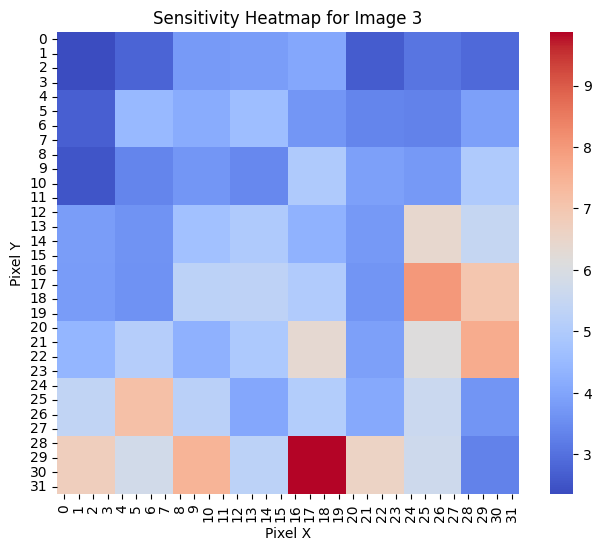}
        \includegraphics[width=0.092\linewidth]{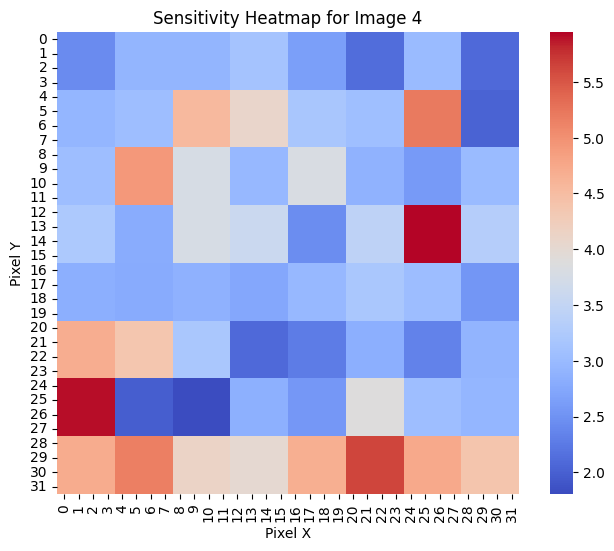}
        \includegraphics[width=0.092\linewidth]{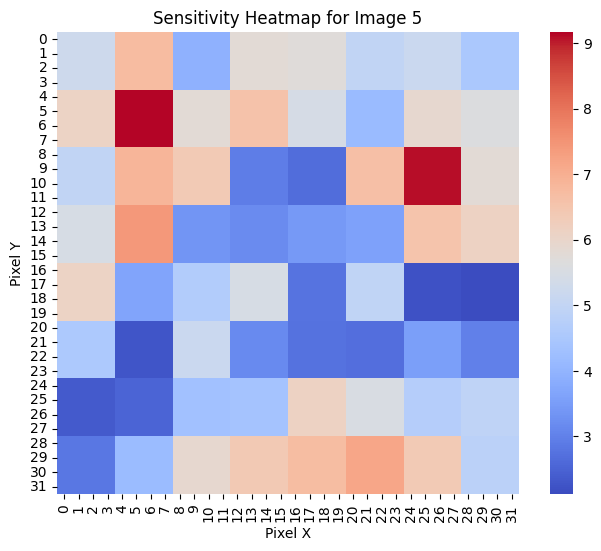}
        \includegraphics[width=0.092\linewidth]{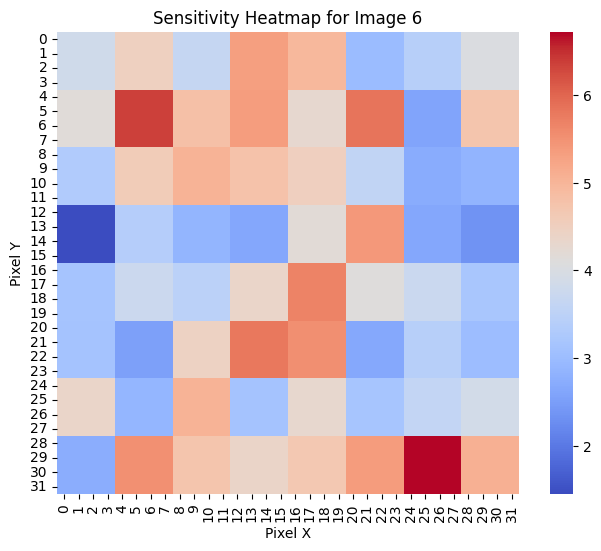}
        \includegraphics[width=0.092\linewidth]{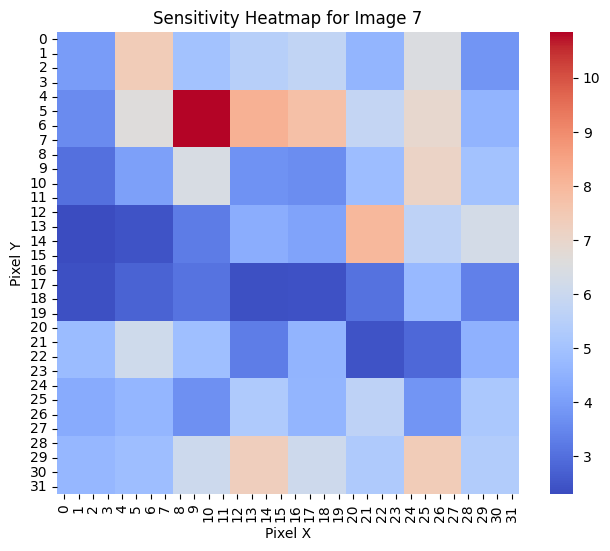}
        \includegraphics[width=0.092\linewidth]{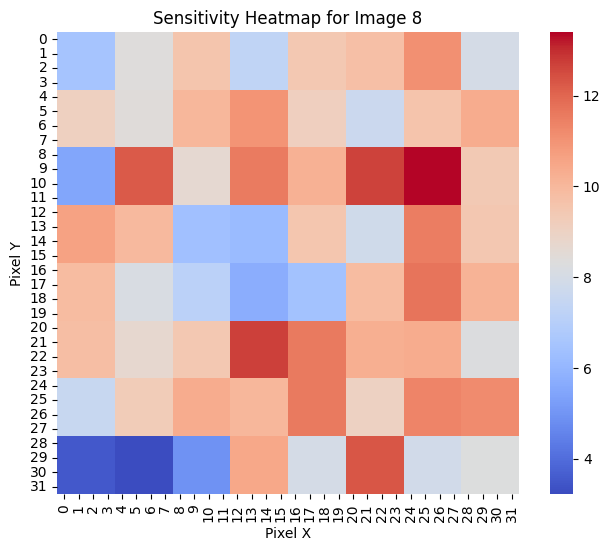}
%    \end{subfigure}
\\
%        \begin{subfigure}{\linewidth}
        \includegraphics[width=0.092\linewidth]{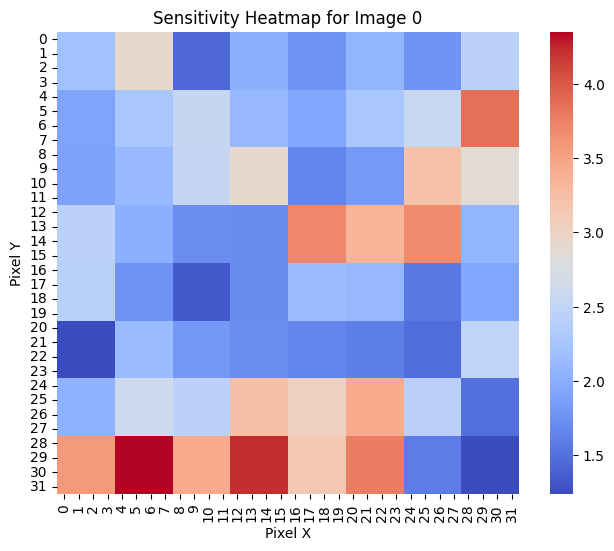}
        \includegraphics[width=0.092\linewidth]{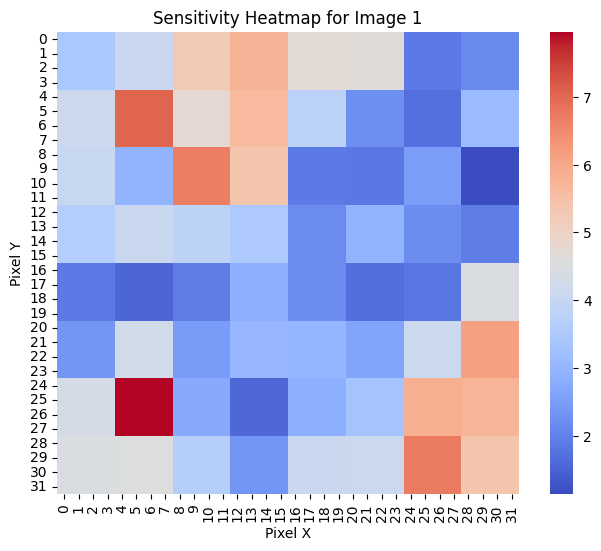}
        \includegraphics[width=0.092\linewidth]{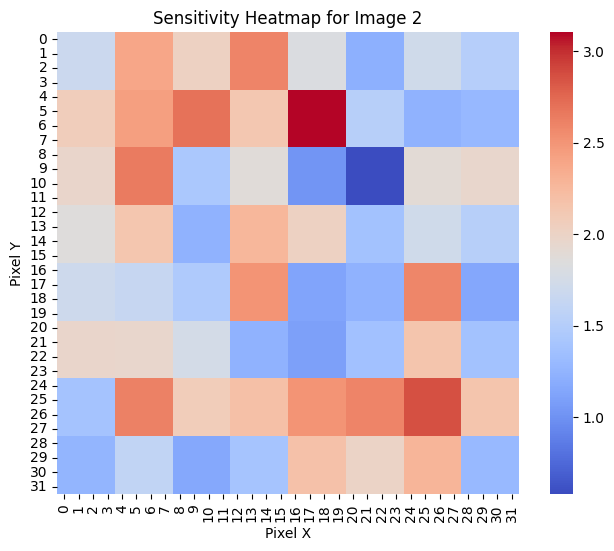}
        \includegraphics[width=0.092\linewidth]{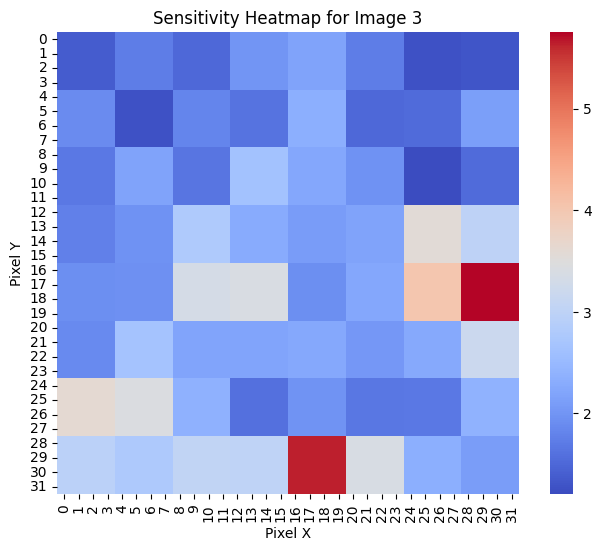}
        \includegraphics[width=0.092\linewidth]{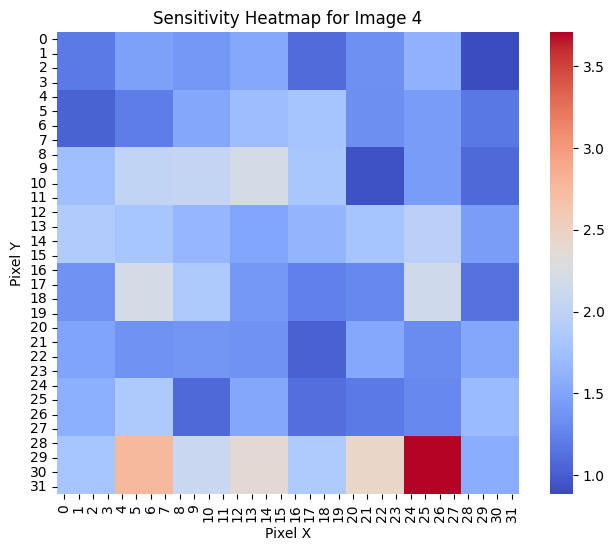}
        \includegraphics[width=0.092\linewidth]{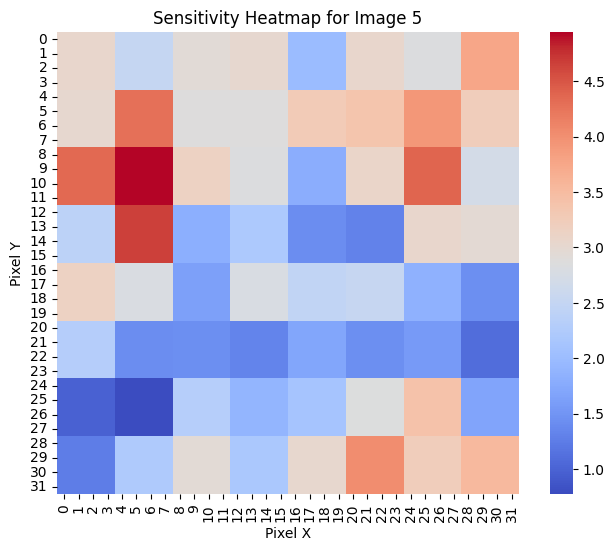}
        \includegraphics[width=0.092\linewidth]{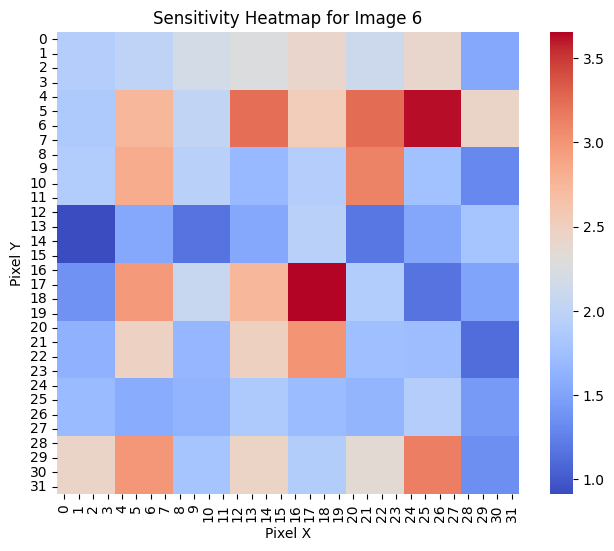}
        \includegraphics[width=0.092\linewidth]{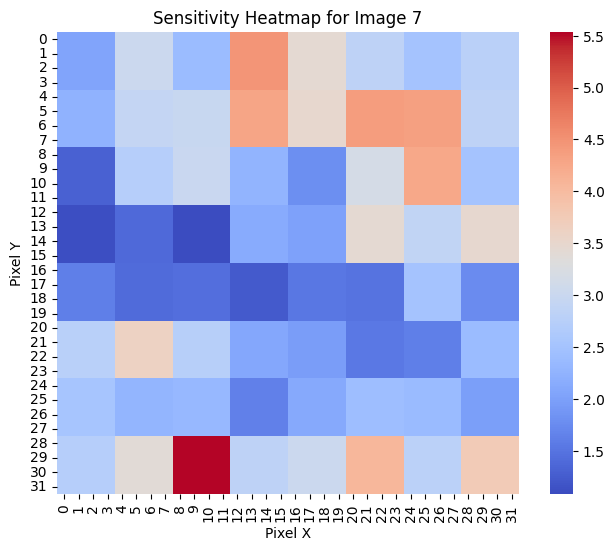}
        \includegraphics[width=0.092\linewidth]{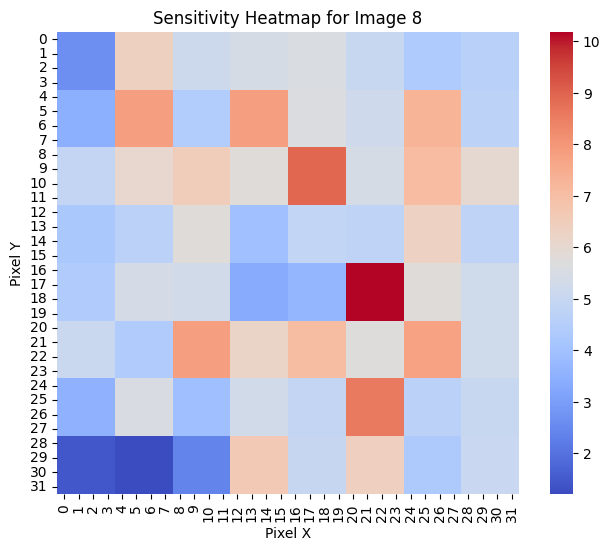}
%    \end{subfigure}
\\
%        \begin{subfigure}{\linewidth}
        \includegraphics[width=0.092\linewidth]{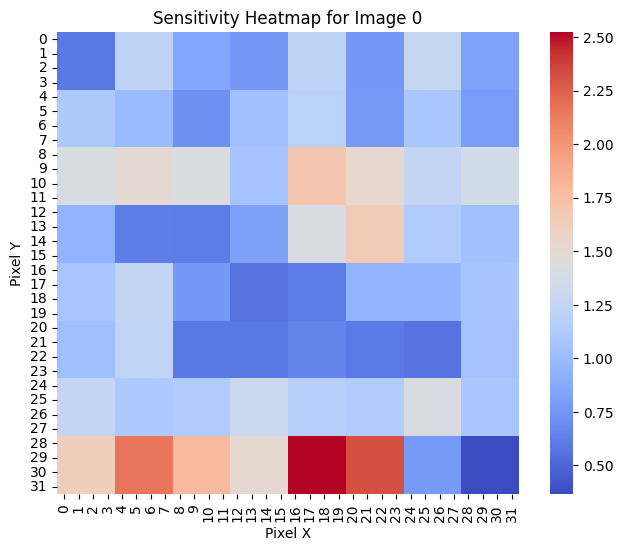}
        \includegraphics[width=0.092\linewidth]{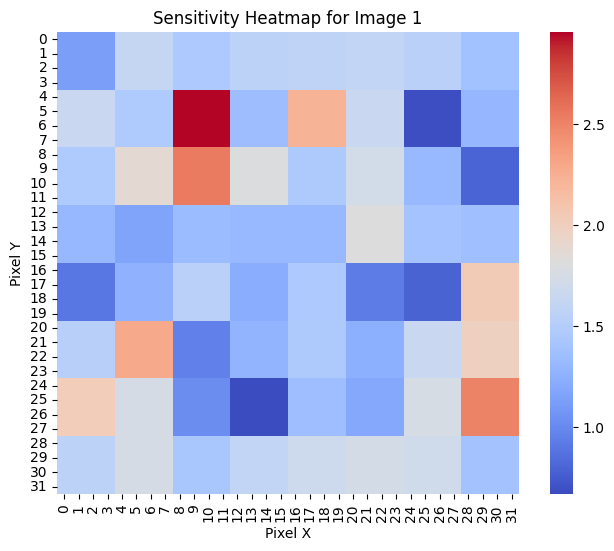}
        \includegraphics[width=0.092\linewidth]{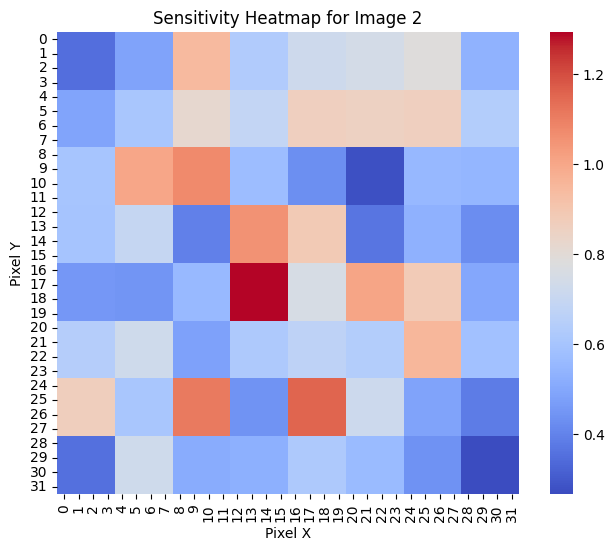}
        \includegraphics[width=0.092\linewidth]{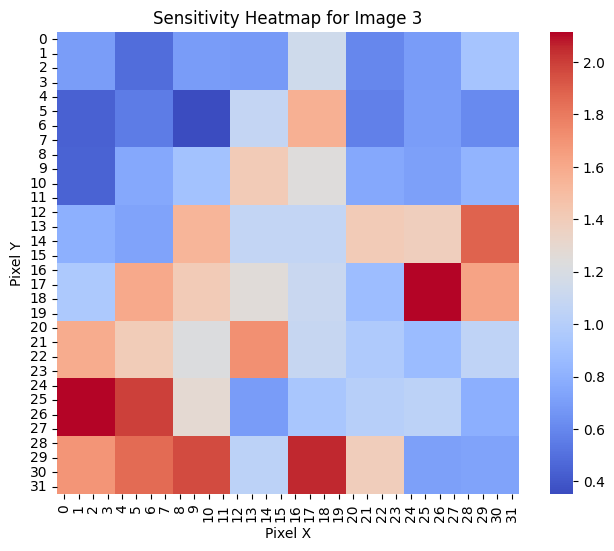}
        \includegraphics[width=0.092\linewidth]{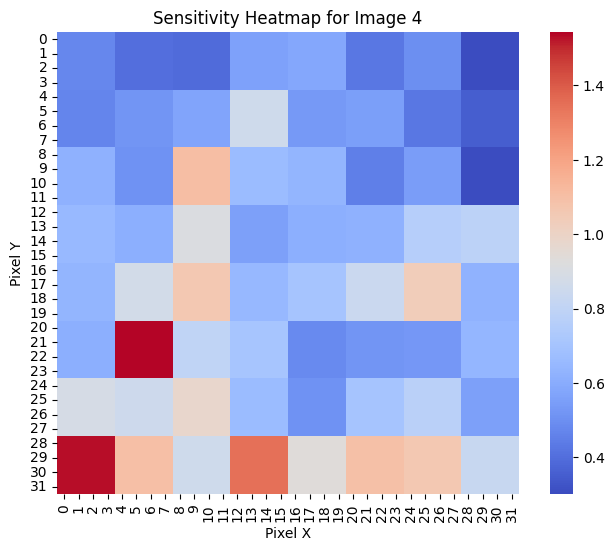}
        \includegraphics[width=0.092\linewidth]{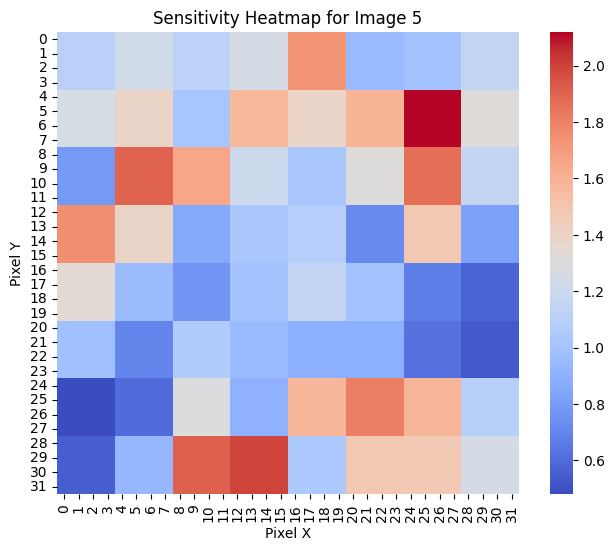}
        \includegraphics[width=0.092\linewidth]{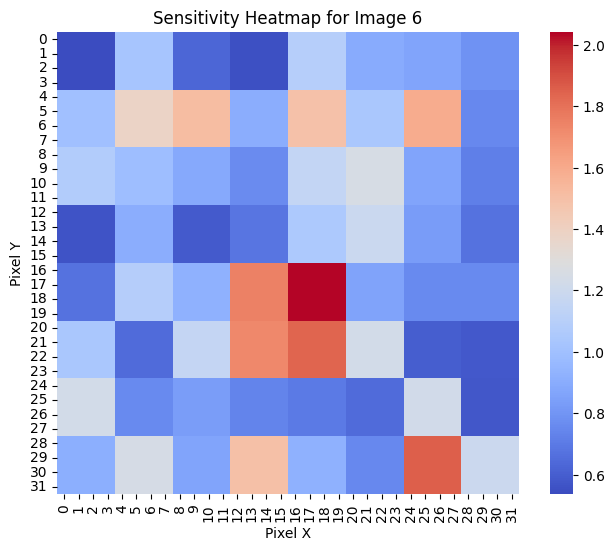}
        \includegraphics[width=0.092\linewidth]{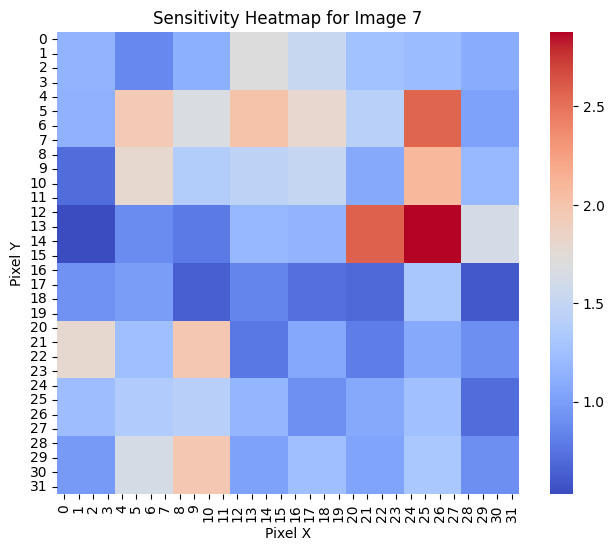}
        \includegraphics[width=0.092\linewidth]{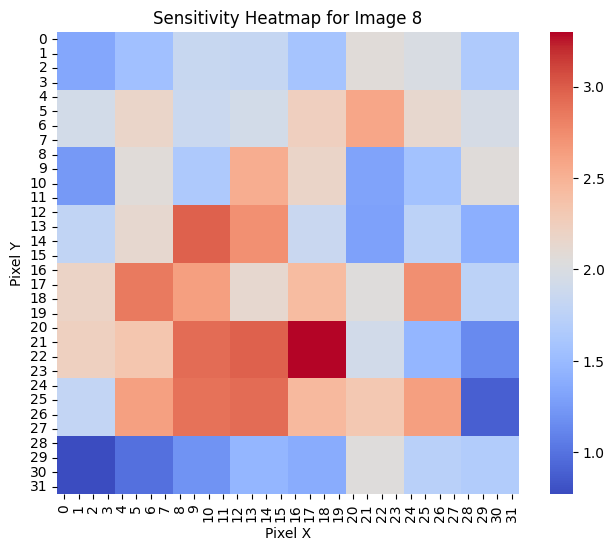}
%    \end{subfigure}
\\
%        \begin{subfigure}{\linewidth}
       \includegraphics[width=0.092\linewidth]{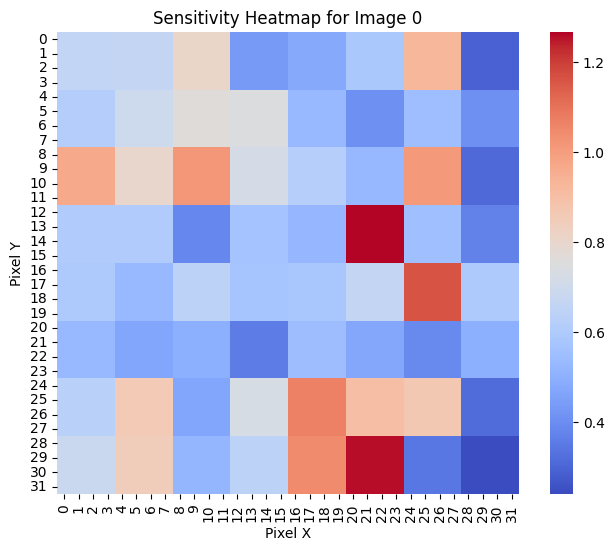}
        \includegraphics[width=0.092\linewidth]{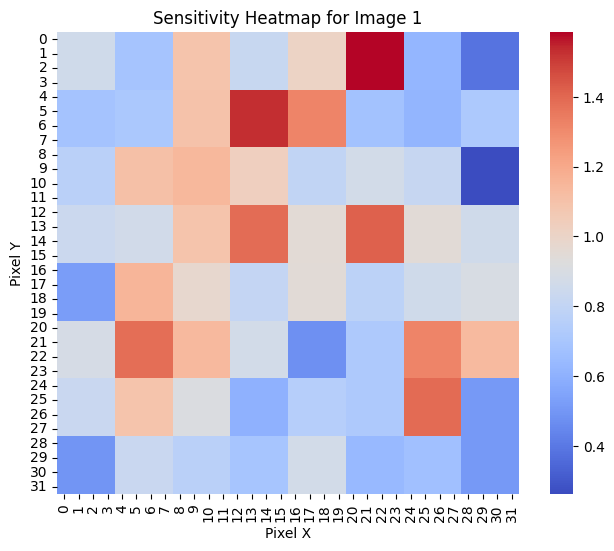}
        \includegraphics[width=0.092\linewidth]{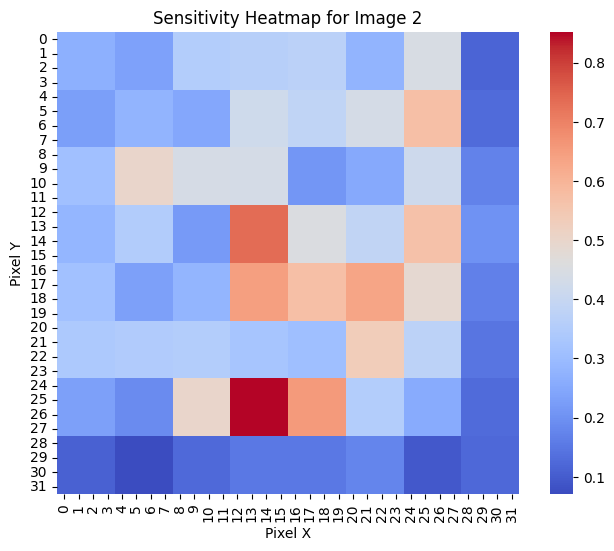}
        \includegraphics[width=0.092\linewidth]{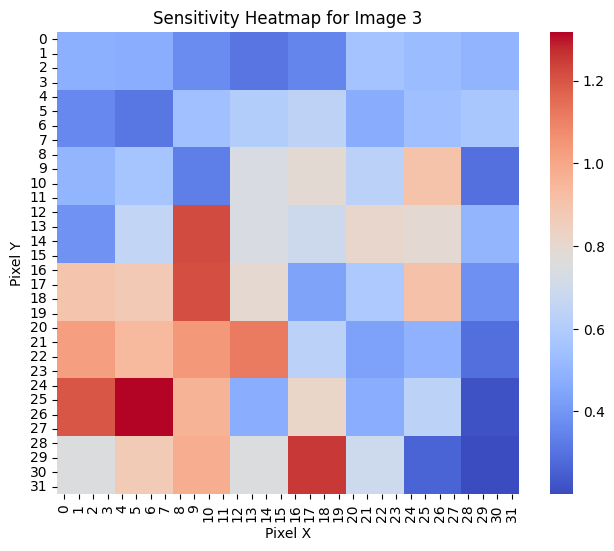}
        \includegraphics[width=0.092\linewidth]{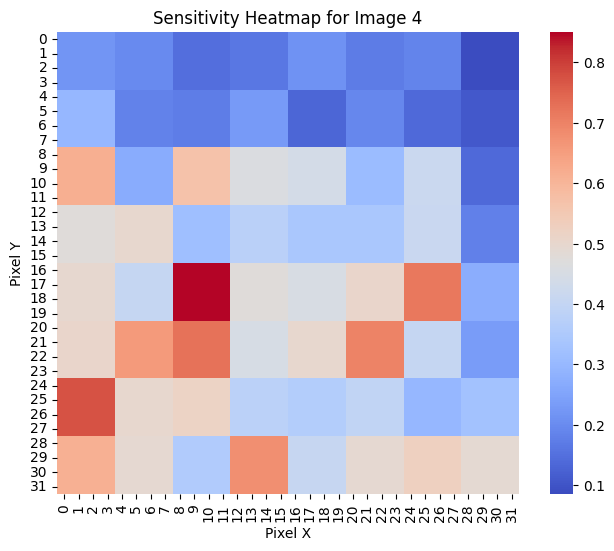}
        \includegraphics[width=0.092\linewidth]{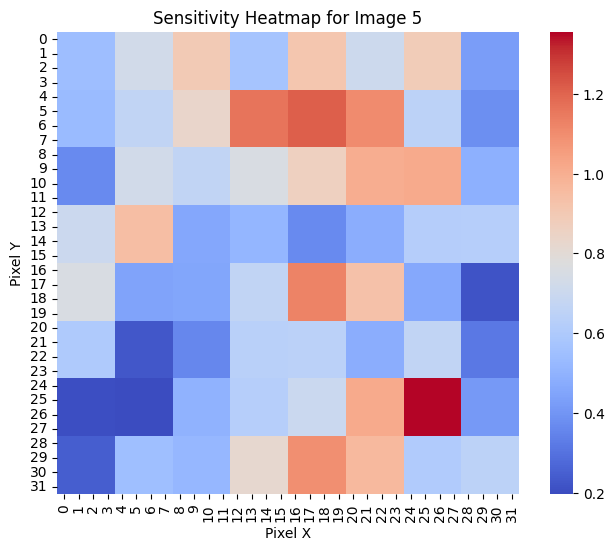}
        \includegraphics[width=0.092\linewidth]{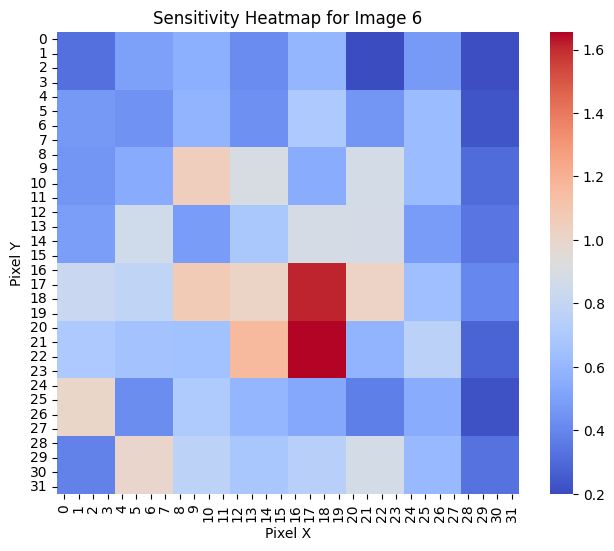}
        \includegraphics[width=0.092\linewidth]{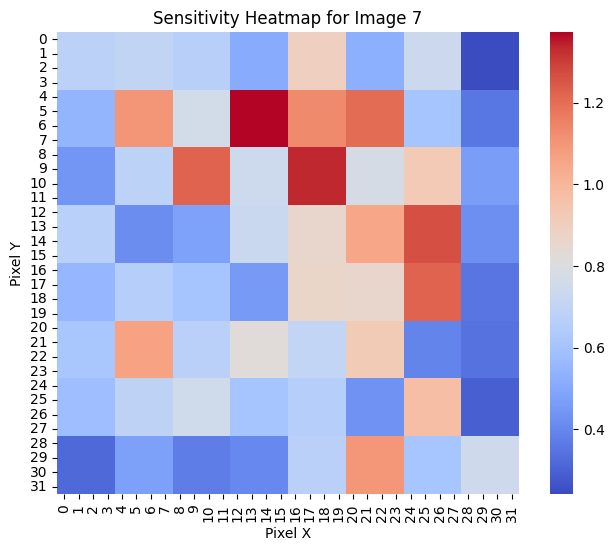}
        \includegraphics[width=0.092\linewidth]{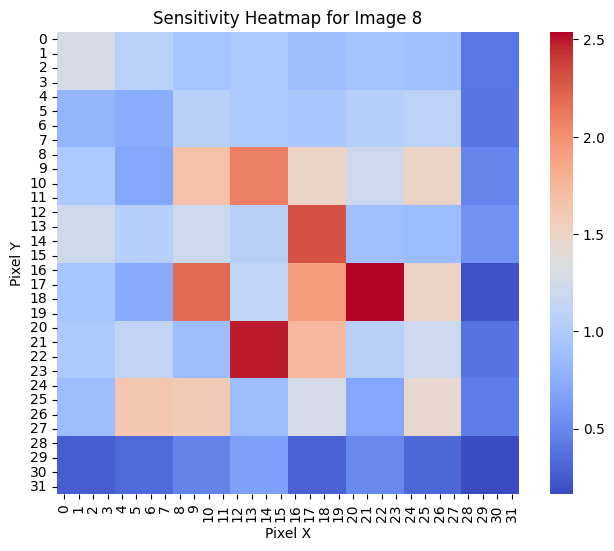}
%    \end{subfigure}
    \vspace{0.2cm} 
    \caption{Sensitivity heatmaps for the red channel of layers 1, 2, 3, and 4 (Top to bottom: layers 1, 2, 3, and 4. The first layer has a scale of $10^{1}$. The scale for layers 2, 3, and 4 is $10^{0}$.}
    \label{fig:ResNet18_red}
\end{figure}

By comparing the sensitivity heatmaps of the red channel of block 5 of the VGG16 model \ref{fig:SensitHeatmapBlock5_RGB}, it is noticeable that in the fourth layer of ResNet18 \ref{fig:ResNet18_red}, the sensitivity level is still $10^{0}$, and it is clear that in the final layers of ResNet18, the sensitivity is much more intense, and the feature capture is more diverse. Using activation maximization (AM), we can observe the ability of our model to capture image features and confirm our assumptions.

\begin{figure}[htbp]
    \centering
    %\begin{subfigure}{1\linewidth}
        \includegraphics[width=0.11\linewidth]{images/AM/plane4.png}
        \includegraphics[width=0.11\linewidth]{images/AM/car2.png}
        \includegraphics[width=0.11\linewidth]{images/AM/frog4.png}
        \includegraphics[width=0.11\linewidth]{images/AM/bird4.png}
        \includegraphics[width=0.11\linewidth]{images/AM/deer5.png}
        \includegraphics[width=0.11\linewidth]{images/AM/horse1.png}
        \includegraphics[width=0.11\linewidth]{images/AM/ship1.png}
        \includegraphics[width=0.11\linewidth]{images/AM/cat10.png}
    %\end{subfigure}
    \\
    %\begin{subfigure}{1\linewidth}
        \includegraphics[width=0.11\linewidth]{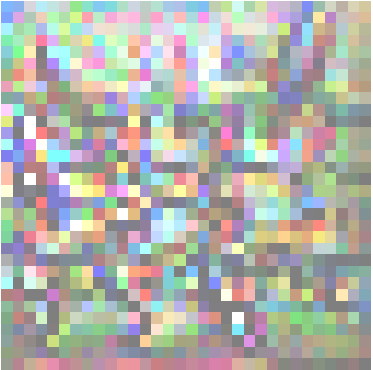}
        \includegraphics[width=0.11\linewidth]{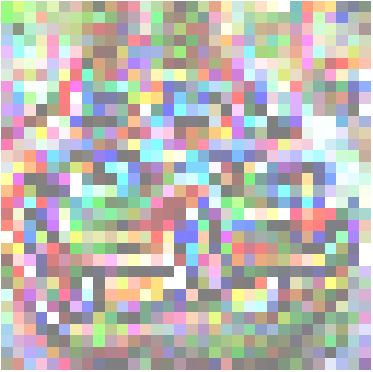}
        \includegraphics[width=0.11\linewidth]{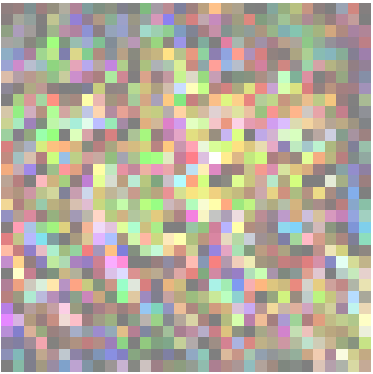}
        \includegraphics[width=0.11\linewidth]{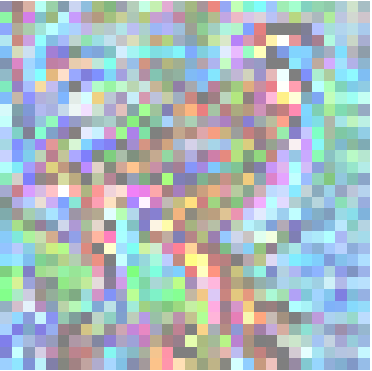}
        \includegraphics[width=0.11\linewidth]{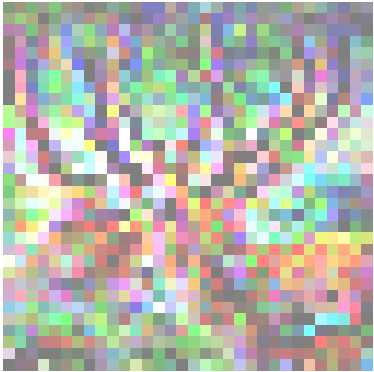}
        \includegraphics[width=0.11\linewidth]{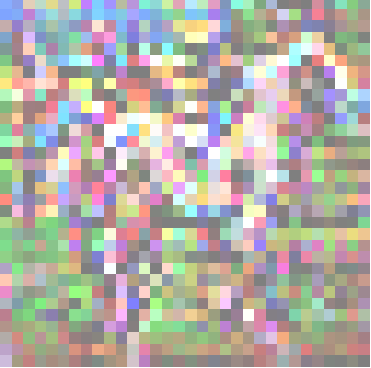}
        \includegraphics[width=0.11\linewidth]{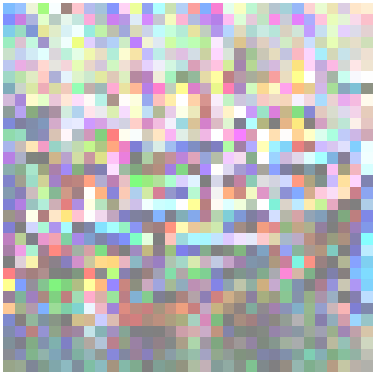}
        \includegraphics[width=0.11\linewidth]{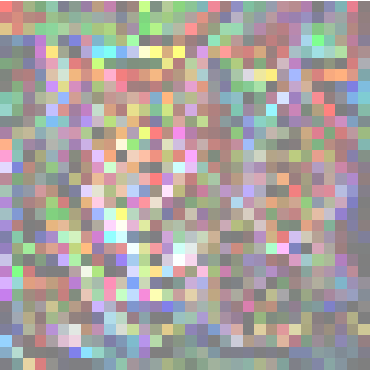}
    %\end{subfigure}
    \caption{Application of AM with TV-regularization for ResNet18.}
    \label{fig:ResNet_AM}
\end{figure}

Although we can still clearly see the contours of the target class features \ref{fig:ResNet_AM}, the activation of the image after using the ResNet-18 model looks more blurred compared to VGG-16 \ref{fig:AM}, with the activation of many details and abstract features, which confirms our hypothesis. Therefore, for tasks with large datasets and more complex features, the deep ResNet network can effectively improve model performance. However, for tasks with small datasets and minimal details, residual networks (ResNet) may cause noise to be perceived as features, which would degrade the model's performance.

\section{Application of VGG-16 and activation maximization to medical data}

During the COVID-19 pandemic, diagnostics relied on fast, safe, and highly sensitive tools. Given the significant practical advantages of lung ultrasound (LUS) compared to other imaging methods, as well as the challenges doctors face when recognizing patterns, we aim to utilize machine learning to assist in diagnosis based on LUS.

We attempt to use the VGG-16 model for classifying lung ultrasound images for COVID-19, bacterial pneumonia, and normal lung conditions \cite{born2021accelerating}. Additionally, we plan to use activation maximization to check how well the model recognizes relevant medical features. Our data is sourced from the repository\footnote{ \url{https://github.com/BorgwardtLab/covid19_ultrasound}}.

\begin{figure}[h]
    \centering
%    \begin{subfigure}{1\linewidth}
        \includegraphics[width=0.3\linewidth]{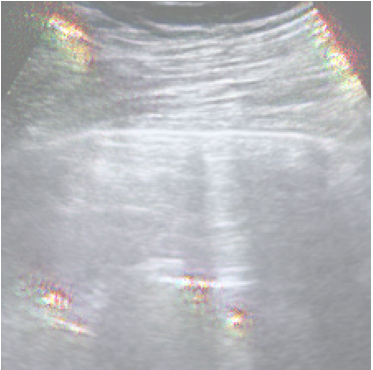}
        \includegraphics[width=0.3\linewidth]{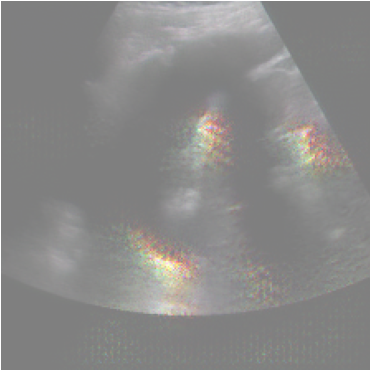}
        \includegraphics[width=0.3\linewidth]{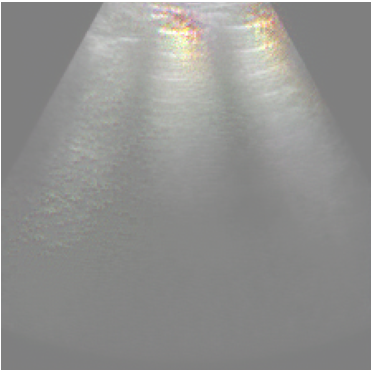}
%    \end{subfigure}
\\
    \vspace{0.2cm} 
%    \begin{subfigure}{1\linewidth}
        \includegraphics[width=0.3\linewidth]{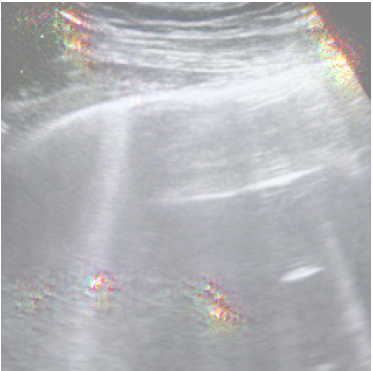}
        \includegraphics[width=0.3\linewidth]{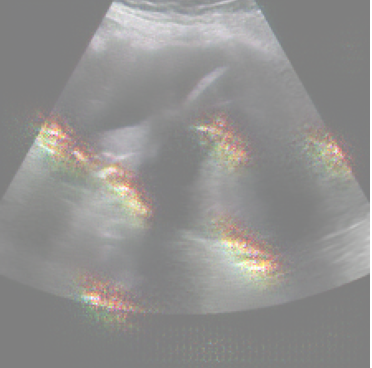}
        \includegraphics[width=0.3\linewidth]{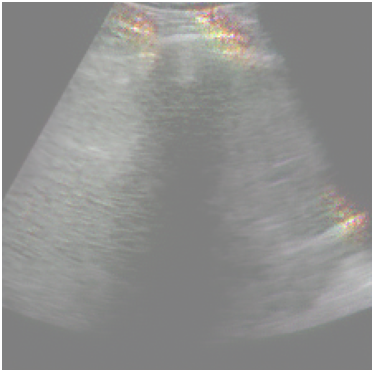}
%    \end{subfigure}
    \caption{AM for COVID-19 (highlights B-lines), bacterial pneumonia (highlights pleural consolidations), and healthy lungs (highlights A-lines).}
    \label{fig:LUS_AM}
\end{figure}

As a result, our model achieved an accuracy of 94\% on the training dataset. After applying activation maximization, we obtained sensitivity heatmap images. In Fig. \ref{fig:LUS_AM}, from left to right, the images represent COVID-19, bacterial pneumonia, and normal lungs.

The classification data and code can be found at the our repository \footnote{\url{https://github.com/LUCI1a/VGG-16_for_LUS.git}}. We have also explored some important features of lung ultrasound \cite{critical2020}:
\begin{itemize}
    \item \textbf{A-lines} — A-lines are linear hyperechoic structures located parallel to the pleural line, arising due to multiple reflections of the sound beam from the pleura as it passes perpendicularly. They appear as straight, clear, and evenly spaced lines, with a gradual decrease in echogenicity as they deepen.    
    \item \textbf{B-lines} — B-lines are hyperechoic linear artifacts that start from the pleural line and extend vertically into the lung, reaching the edge of the ultrasound scanner's screen.    
    \item \textbf{lung consolidation} — Lung consolidation refers to areas of lung tissue that appear "liver-like" on ultrasound images, indicating significant tissue consolidation.
\end{itemize}

\begin{figure}[h]
    \centering
%    \begin{subfigure}{1\linewidth}
        \includegraphics[width=0.3\linewidth]{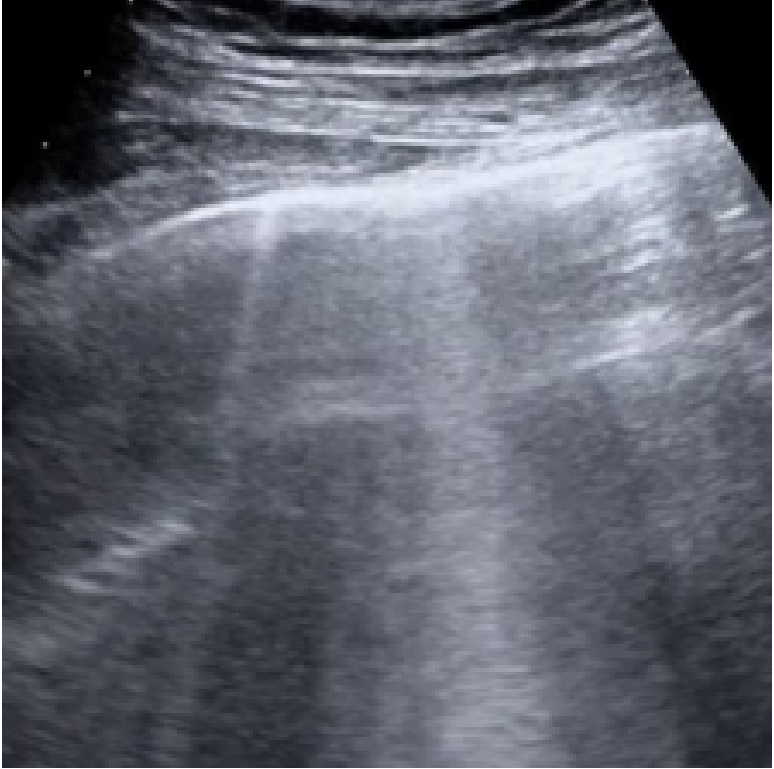}
        \includegraphics[width=0.3\linewidth]{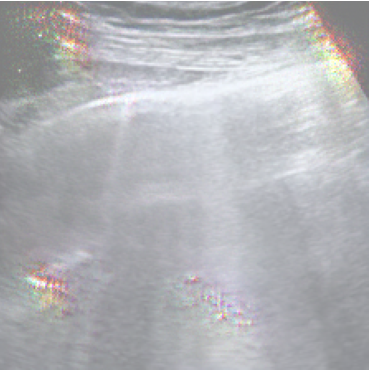}
        \includegraphics[width=0.3\linewidth]{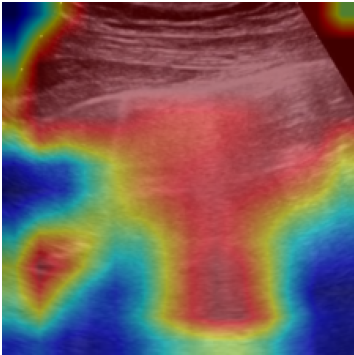}
%    \end{subfigure}
\\
    \vspace{0.2cm} 
%    \begin{subfigure}{1\linewidth}
        \includegraphics[width=0.3\linewidth]{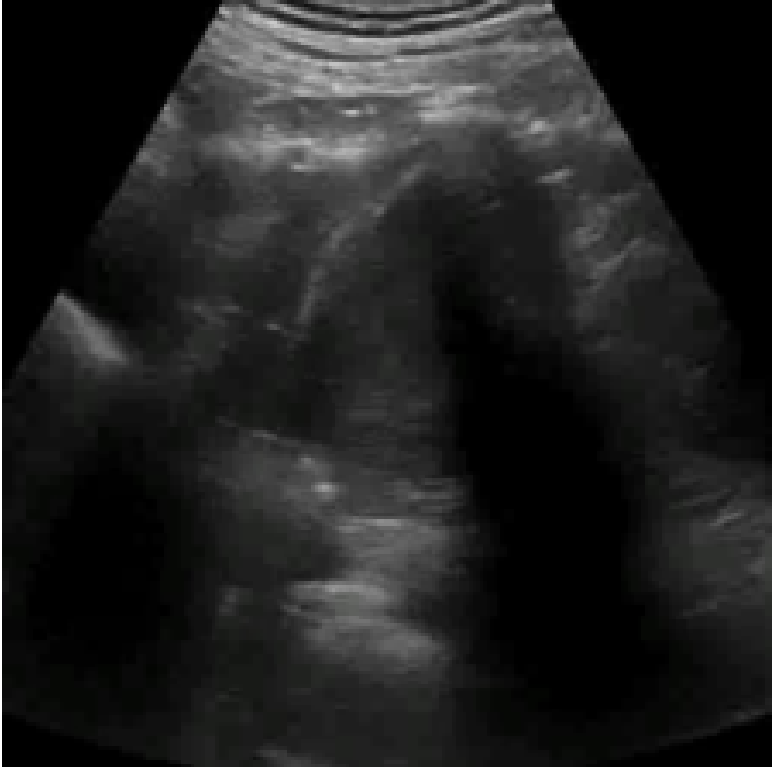}
        \includegraphics[width=0.3\linewidth]{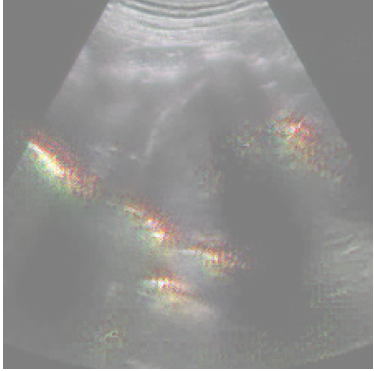}
        \includegraphics[width=0.3\linewidth]{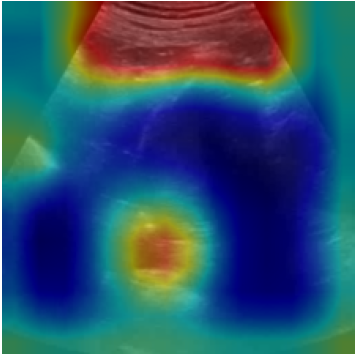}
%    \end{subfigure}
\\
    \vspace{0.2cm} 
%    \begin{subfigure}{1\linewidth}
        \includegraphics[width=0.3\linewidth]{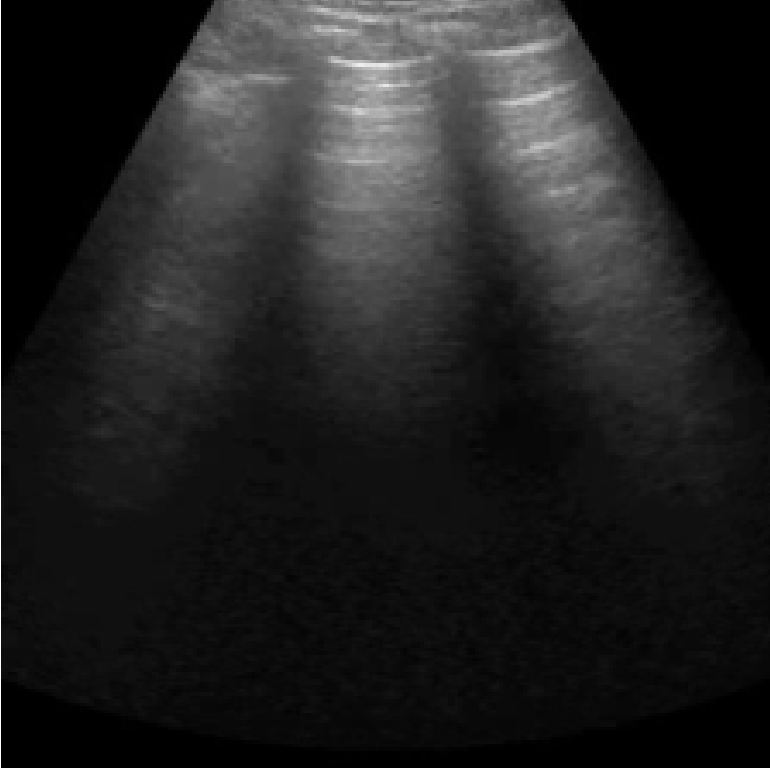}
        \includegraphics[width=0.3\linewidth]{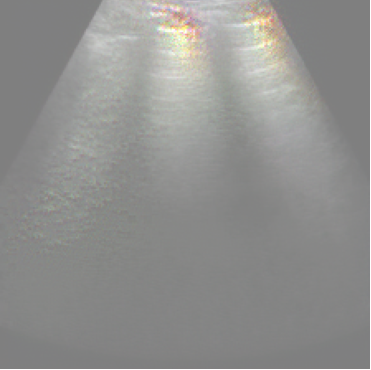}
        \includegraphics[width=0.3\linewidth]{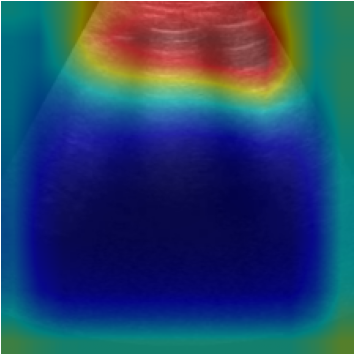}
%    \end{subfigure}
    \caption{The images of lungs with COVID-19, bacterial pneumonia, and healthy lungs are arranged from top to bottom. On the right is the Grad-Cam image.}
    \label{fig:LUS_CAM}
\end{figure}

Observing the results (see Figures \ref{fig:LUS_AM} and  \ref{fig:LUS_CAM}), it can be noted that various types of ultrasound features were well activated. For COVID-19 activation images, the ray-like B-lines exhibit a clear bright glow, but since their boundaries are ray-like and resemble B-lines, they were also activated. This could be one of the reasons for model errors.

For bacterial pneumonia, uneven hyperechoic spots were activated. For healthy lungs, the echogenicity was uniform, and the A-lines, parallel and equidistant, were brightly activated. This aligns with the clinical diagnostic signs observed in \cite{born2021accelerating} and matches the observations made using the Grad-Cam technique \cite{selvaraju2017grad} employed in \cite{born2021accelerating}. However, compared to activation maximization, Grad-Cam \footnote{\url{https://github.com/jacobgil/pytorch-grad-cam.git}} has a broader sensitive area \ref{fig:LUS_CAM}, while the activation maximization method captures the features more effectively.

\section{Conclusion}

In this work, we have explored several methods of the sensitivity analysis for neural network models. At first, we show application of the global sensitivity to a simple feedforward neural network using the Sobol method on a small dataset of clinical diabetes data. The sensitivity analysis identifies three key parameters whose values significantly affect the classifier's performance. The results of the analysis align with a simpler assessment of the dataset using principal component analysis, but the global sensitivity analysis provides a clearer ranking of the input variables by their importance. These observations could be useful for building reliable and interpretable neural network classifiers in the context of tabular data \cite{di2023explainable}.

Increasing the input space dimensionality for neural networks requires a much larger number of quasi-random points for computing Sobol indices. One partial solution to this problem could be using parallel computing \cite{gasanov2020sensitivity} when performing Monte Carlo integration, although computational resources are still limited. Another potential solution could be combining Sobol analysis with low-rank decomposition of the target function in the form of tensor trains \cite{ballester2019sobol}.

Application of the global sensitivity analysis methods to convolutional neural networks like VGG-16 is difficult even for processing the small images from the CIFAR-10 dataset due to the significantly higher dimensionality of the input variables (3072 variables for $32 \times 32 \times 3$ images versus 8 for tabular clinical data). Instead, we present a local sensitivity analysis using heatmaps with small random perturbations of input pixels and employing the activation maximization method.

A gradual decrease in the sensitivity of the hidden layers of the VGG-16 convolutional neural network to local perturbations of input pixels was observed as the depth increased. The presence of residual blocks in the ResNet-18 architecture helps to preserve the sensitivity of deeper layers to local perturbations of input images. At the same time, the activation maximization method clearly demonstrates that the VGG-16 model effectively captures the contours of the classified objects, which appear sharper compared to the results for ResNet-18.

In the case of activation of images from the incorrect classes, we observe that the neural network mistakenly identifies the features (contours) belonging to objects of another class. These results seem to be interesting for analyzing the properties and searching for the structure of adversarial attacks \cite{khrulkov2018art}, which can lead to erroneous performance of neural network classifiers when processing images.

Finally, using the VGG-16 model, we classily the clinical ultrasound images of the lungs, successfully distinguishing the COVID-19, bacterial pneumonia, and normal lungs. The activation maximization method effectively highlights the medical features, although some boundaries might be mistakenly recognized as significant medical features and could lead to misclassification. Therefore, additional preprocessing of image boundaries is a promising solution.

The conducted research demonstrates the potential of a convolutional neural networks in medicine and the possibility of using the activation maximization method to assess their scientific validity.

\section*{Acknowledgements}

Sergey Matveev was supported by Russian Science Foundation project (project No. 25-11-00392).  Authors are grateful to their colleagues from Lomonosov MSU and INM RAS for useful discussions.

\bibliography{paper_en}

\begin{thebibliography}{10}

\bibitem{hopfield1982neural}
J.~J. Hopfield, ``Neural networks and physical systems with emergent collective
  computational abilities.,'' {\em Proceedings of the national academy of
  sciences}, vol.~79, no.~8, pp.~2554--2558, 1982.

\bibitem{jumper2021highly}
J.~Jumper, R.~Evans, A.~Pritzel, T.~Green, M.~Figurnov, O.~Ronneberger,
  K.~Tunyasuvunakool, R.~Bates, A.~{\v{Z}}{\'\i}dek, A.~Potapenko, {\em
  et~al.}, ``Highly accurate protein structure prediction with {A}lpha{F}old,''
  {\em nature}, vol.~596, no.~7873, pp.~583--589, 2021.

\bibitem{antsiferova2022video}
A.~Antsiferova, S.~Lavrushkin, M.~Smirnov, A.~Gushchin, D.~Vatolin, and
  D.~Kulikov, ``Video compression dataset and benchmark of learning-based
  video-quality metrics,'' {\em Advances in Neural Information Processing
  Systems}, vol.~35, pp.~13814--13825, 2022.

\bibitem{liu2019neural}
H.~Liu, T.~Chen, M.~Lu, Q.~Shen, and Z.~Ma, ``Neural video compression using
  spatio-temporal priors,'' {\em arXiv preprint arXiv:1902.07383}, 2019.

\bibitem{vaswani2017attention}
A.~Vaswani, ``Attention is all you need,'' {\em Advances in Neural Information
  Processing Systems}, 2017.

\bibitem{kharyuk2018employing}
P.~Kharyuk, D.~Nazarenko, I.~Oseledets, I.~Rodin, O.~Shpigun, A.~Tsitsilin, and
  M.~Lavrentyev, ``Employing fingerprinting of medicinal plants by means of
  lc-ms and machine learning for species identification task,'' {\em Scientific
  reports}, vol.~8, no.~1, p.~17053, 2018.

\bibitem{matveev2021overview}
S.~A. Matveev, I.~V. Oseledets, E.~S. Ponomarev, and A.~V. Chertkov, ``Overview
  of visualization methods for artificial neural networks,'' {\em Computational
  Mathematics and Mathematical Physics}, vol.~61, no.~5, pp.~887--899, 2021.

\bibitem{ramsauer2020hopfield}
H.~Ramsauer, B.~Sch{\"a}fl, J.~Lehner, P.~Seidl, M.~Widrich, T.~Adler,
  L.~Gruber, M.~Holzleitner, M.~Pavlovi{\'c}, G.~K. Sandve, {\em et~al.},
  ``Hopfield networks is all you need,'' {\em arXiv preprint arXiv:2008.02217},
  2020.

\bibitem{khrulkov2018art}
V.~Khrulkov and I.~Oseledets, ``Art of singular vectors and universal
  adversarial perturbations,'' in {\em Proceedings of the IEEE Conference on
  Computer Vision and Pattern Recognition}, pp.~8562--8570, 2018.

\bibitem{mopuri2017fast}
K.~R. Mopuri, U.~Garg, and R.~V. Babu, ``Fast feature fool: A data independent
  approach to universal adversarial perturbations,'' {\em arXiv preprint
  arXiv:1707.05572}, 2017.

\bibitem{angelov2021explainable}
P.~P. Angelov, E.~A. Soares, R.~Jiang, N.~I. Arnold, and P.~M. Atkinson,
  ``Explainable artificial intelligence: an analytical review,'' {\em Wiley
  Interdisciplinary Reviews: Data Mining and Knowledge Discovery}, vol.~11,
  no.~5, p.~e1424, 2021.

\bibitem{kharyuk2025exploring}
P.~Kharyuk, S.~Matveev, and I.~Oseledets, ``Exploring specialization and
  sensitivity of convolutional neural networks in the context of simultaneous
  image augmentations,'' {\em arXiv preprint arXiv:2503.03283}, 2025.

\bibitem{smith1988using}
J.~W. Smith, J.~E. Everhart, W.~Dickson, W.~C. Knowler, and R.~S. Johannes,
  ``Using the {A}{D}{A}{P} learning algorithm to forecast the onset of diabetes
  mellitus,'' in {\em Proceedings of the annual symposium on computer
  application in medical care}, p.~261, American Medical Informatics
  Association, 1988.

\bibitem{akturk2020diabetes}
M.~Akturk, ``Diabetes dataset,'' {\em Kaggle. com}, 2020.

\bibitem{spall2005introduction}
J.~C. Spall, {\em Introduction to stochastic search and optimization:
  estimation, simulation, and control}.
\newblock John Wiley \& Sons, 2005.

\bibitem{sobol2001global}
I.~M. Sobol, ``Global sensitivity indices for nonlinear mathematical models and
  their {M}onte {C}arlo estimates,'' {\em Mathematics and computers in
  simulation}, vol.~55, no.~1-3, pp.~271--280, 2001.

\bibitem{rosolem2012fully}
R.~Rosolem, H.~V. Gupta, W.~J. Shuttleworth, X.~Zeng, and L.~G.~G.
  de~Gon{\c{c}}alves, ``A fully multiple-criteria implementation of the sobol
  method for parameter sensitivity analysis,'' {\em Journal of Geophysical
  Research: Atmospheres}, vol.~117, no.~D7, 2012.

\bibitem{herman2017salib}
J.~Herman and W.~Usher, ``{S}{A}{L}ib: An open-source python library for
  sensitivity analysis,'' {\em Journal of Open Source Software}, vol.~2, no.~9,
  p.~97, 2017.

\bibitem{sobol1998quasi}
I.~M. Sobol, ``On quasi-monte carlo integrations,'' {\em Mathematics and
  computers in simulation}, vol.~47, no.~2-5, pp.~103--112, 1998.

\bibitem{gasanov2020sensitivity}
M.~Gasanov, A.~Petrovskaia, A.~Nikitin, S.~Matveev, P.~Tregubova, M.~Pukalchik,
  and I.~Oseledets, ``Sensitivity analysis of soil parameters in crop model
  supported with high-throughput computing,'' in {\em International Conference
  on Computational Science}, pp.~731--741, Springer, 2020.

\bibitem{sysoev2019sensitivity}
A.~Sysoev, A.~Ciurlia, R.~Sheglevatych, and S.~Blyumin, ``Sensitivity analysis
  of neural network models: {A}pplying methods of analysis of finite
  fluctuations,'' {\em Periodica polytechnica Electrical engineering and
  computer science}, vol.~63, no.~4, pp.~306--311, 2019.

\bibitem{zhang2020sensitivity}
C.~Zhang, Y.~Zhang, J.~Zhao, and J.~Luo, ``Sensitivity analysis of weather
  factors affecting pv module output power based on artificial neural network
  and sobol algorithm,'' in {\em 2020 IEEE/IAS Industrial and Commercial Power
  System Asia (I\&CPS Asia)}, pp.~246--250, IEEE, 2020.

\bibitem{simonyan2014very}
K.~Simonyan, ``Very deep convolutional networks for large-scale image
  recognition,'' {\em arXiv preprint arXiv:1409.1556}, 2014.

\bibitem{alex2009learning}
K.~Alex, ``Learning multiple layers of features from tiny images,'' {\em
  https://www. cs. toronto. edu/kriz/learning-features-2009-TR. pdf}, 2009.

\bibitem{JSSv102i07}
J.~Pizarroso, J.~Portela, and A.~Muñoz, ``Neural{S}ens: {S}ensitivity
  {A}nalysis of {N}eural {N}etworks,'' {\em Journal of Statistical Software},
  vol.~102, no.~7, p.~1–36, 2022.

\bibitem{recht2018cifar}
B.~Recht, R.~Roelofs, L.~Schmidt, and V.~Shankar, ``Do cifar-10 classifiers
  generalize to cifar-10?,'' {\em arXiv preprint arXiv:1806.00451}, 2018.

\bibitem{pospelov2025fast}
N.~Pospelov, A.~Chertkov, M.~Beketov, I.~Oseledets, and K.~Anokhin, ``Fast
  gradient-free activation maximization for neurons in spiking neural
  networks,'' {\em Neurocomputing}, vol.~618, p.~129070, 2025.

\bibitem{simonyan2013deep}
K.~Simonyan, A.~Vedaldi, and A.~Zisserman, ``Deep inside convolutional
  networks: {V}isualising image classification models and saliency maps,'' {\em
  arXiv preprint arXiv:1312.6034}, 2013.

\bibitem{he2016deep}
K.~He, X.~Zhang, S.~Ren, and J.~Sun, ``Deep residual learning for image
  recognition,'' in {\em Proceedings of the IEEE conference on computer vision
  and pattern recognition}, pp.~770--778, 2016.

\bibitem{born2021accelerating}
J.~Born, N.~Wiedemann, M.~Cossio, C.~Buhre, G.~Br{\"a}ndle, K.~Leidermann,
  J.~Goulet, A.~Aujayeb, M.~Moor, B.~Rieck, {\em et~al.}, ``Accelerating
  detection of lung pathologies with explainable ultrasound image analysis,''
  {\em Applied Sciences}, vol.~11, no.~2, p.~672, 2021.

\bibitem{critical2020}
C.~U.~G. of~Ultrasound Medicine~Committee, C.~M.~E. Association, {\em et~al.},
  ``Gan ran xing fei yan chao sheng zhen duan zhuan jia jian yi,'' {\em Zhong
  Hua Yi Xue Chao Sheng Za Zhi(Dian Zi Ban)}, vol.~17, no.~03, p.~244, 2020.

\bibitem{selvaraju2017grad}
R.~R. Selvaraju, M.~Cogswell, A.~Das, R.~Vedantam, D.~Parikh, and D.~Batra,
  ``Grad-cam: Visual explanations from deep networks via gradient-based
  localization,'' in {\em Proceedings of the IEEE international conference on
  computer vision}, pp.~618--626, 2017.

\bibitem{di2023explainable}
F.~Di~Martino and F.~Delmastro, ``Explainable {A}{I} for clinical and remote
  health applications: a survey on tabular and time series data,'' {\em
  Artificial Intelligence Review}, vol.~56, no.~6, pp.~5261--5315, 2023.

\bibitem{ballester2019sobol}
R.~Ballester-Ripoll, E.~G. Paredes, and R.~Pajarola, ``Sobol tensor trains for
  global sensitivity analysis,'' {\em Reliability Engineering \& System
  Safety}, vol.~183, pp.~311--322, 2019.

\end{thebibliography}
\bibliographystyle{ieeetr}

\end{document}